\documentclass[review]{elsarticle}

\usepackage{hyperref}

\journal{Elsevier}

\usepackage{color}
\usepackage{amsmath}
\usepackage{graphicx}
\usepackage{subfigure}
\usepackage{rotating}
\usepackage{multirow, makecell}
\usepackage{geometry}
\usepackage{appendix}
\usepackage{booktabs}

\biboptions{sort&compress}              %参考文献折叠

\begin{document}

% \graphicspath{{figures/}} 
	
\begin{frontmatter}
		
\title{Numerical stability analysis of shock-capturing methods for strong shocks II: high-order finite-volume schemes}
		
\author[mymainaddress]{Weijie Ren}

\author[mymainaddress]{Wenjia Xie\corref{mycorrespondingauthor}}
\cortext[mycorrespondingauthor]{Corresponding author}
\ead{xiewenjia@nudt.edu.cn}

\author[mymainaddress]{Ye Zhang}

\author[mymainaddress]{Hang Yu}

\author[mymainaddress]{Zhengyu Tian}

\address[mymainaddress]{College of Aerospace Science and Engineering, National University of Defense Technology, Hunan 410073, China}

\begin{abstract}
	The shock instability problem commonly arises in flow simulations involving strong shocks, particularly when employing high-order schemes, limiting their applications in hypersonic flow simulations. This study focuses on exploring the numerical characteristics and underlying mechanisms of shock instabilities in fifth-order finite-volume WENO schemes. To this end, for the first time, we have established the matrix stability analysis method for the fifth-order scheme. By predicting the evolution of perturbation errors in the exponential growth stage, this method provides quantitative insights into the behavior of shock-capturing and helps elucidate the mechanisms that cause shock instabilities. Results reveal that even dissipative solvers also suffer from shock instabilities when the spatial accuracy is increased to fifth-order. Further investigation indicates that this is due to the excessively high spatial accuracy of the WENO scheme near the numerical shock structure. Moreover, the shock instability problem of fifth-order schemes is demonstrated to be a multidimensional coupling problem. To stably capture strong shocks, it is crucial to have sufficient dissipation on transverse faces and ensure at least two points within the numerical shock structure in the direction perpendicular to the shock. The source location of instability is also clarified by the matrix stability analysis method, revealing that the instability arises from the numerical shock structure. Additionally, stability analysis demonstrates that local characteristic decomposition helps mitigate shock instabilities in high-order schemes, although the instability still persists. These conclusions pave the way for a better understanding of the shock instability in fifth-order schemes and provide guidance for the development of more reliable high-order shock-capturing methods for compressible flows with high Mach numbers.
\end{abstract}

		
\end{frontmatter}
	
\section{Introduction}
When conducting simulations of supersonic or hypersonic flows containing strong shocks, the shock instability problem becomes a challenge that must be addressed and resolved. A comprehensive classification of various shock instabilities can be found in \cite{Quirk_Contribution_1994}. Among them, the carbuncle phenomenon is one of the most famous. The carbuncle phenomenon is first observed by Peery and Imlay \cite{Perry_Bluntbody_1988} when they compute the supersonic flow field around a blunt-body using the Roe solver. It is conventionally referred to as a spurious solution in blunt-body calculations, characterized by the growth of a protuberance ahead of the bow shock along the stagnation line \cite{Quirk_Contribution_1994}. Moreover, the carbuncle phenomenon is generally considered a concrete manifestation of numerical shock instabilities.

It has been demonstrated that the shock instability problem is a widespread issue. Besides affecting the compressible Euler equations of gas dynamics, the shock instability problem also plagues the approximation of other hyperbolic systems of conservation laws, such as the shallow water equations \cite{Navas-Montilla_Overcoming_2017,Navas-Montilla_Improved_2019,Kemm_Note_2014,Bader_Carbuncle_2014} and the MHD (magnetohydrodynamics) equations \cite{Baty_Robust_2023,Wang_Robust_2023,Kitamura_SLAU2_2020,Lee_Solution_2013,Hanawa_Improving_2008}. When the shock instability problem occurs, simulating shocks will result in unphysical solutions. Consequently, results that heavily rely on shocks become suspicious, such as hypersonic heating \cite{Kitamura_Evaluation_2008,Kitamura_Evaluation_2009,Kitamura_Evaluation_2010,Kitamura_Further_2013,Nastac_Improved_2022}, shock-disturbance interaction \cite{Chuvakhov_ShockCapturing_2021}, and detonation wave simulation \cite{Papalexandris_Numerical_2000,Choi_Celllike_2007,Teng_Numerical_2014}. What exacerbates the situation is that various shock-capturing methods encounter the shock instability problem, including Godunov-type schemes \cite{Quirk_Contribution_1994,Pandolfi_Numerical_2001}, kinetic schemes \cite{Ohwada_Remedy_2013}, and LBM (lattice Boltzmann method) \cite{Chen_Rotated_2023a,Esfahanian_Improvement_2015,Chen_Development_2023}.

Given that the shock instability problem will ruin all efforts to accurately resolve strong shocks and cast doubt on simulation results, extensive research has been conducted to investigate this issue. It has been well demonstrated that the occurrence of the shock instability problem is highly sensitive to some factors, including Riemann solvers \cite{Quirk_Contribution_1994,Pandolfi_Numerical_2001,Dumbser_Matrix_2004,Xie_Numerical_2017}, computational grid \cite{Henderson_Grid_2007,Ohwada_Remedy_2013}, inflow Mach number \cite{Dumbser_Matrix_2004}, numerical shock structure \cite{Dumbser_Matrix_2004,Xie_Numerical_2017,DanielZaide_Shock_2011}, and the order of spatial accuracy \cite{Tu_Evaluation_2014,Ren_Numerical_2023}. For a comprehensive literature review on this topic, readers can refer to these references and the cited works therein. However, the majority of these studies are almost based on first- or second-order schemes. Compared with low-order schemes, high-order schemes offer significant advantages, such as higher accuracy and the ability to capture flow details more effectively. Consequently, they have broad application prospects in the numerical simulations of supersonic/hypersonic flows. Therefore, it is necessary to conduct research on the shock instability problem for high-order schemes. Kemm \cite{Kemm_Heuristical_2018} investigates the carbuncle phenomenon through numerical experiments and heuristic considerations. He thinks that increasing the spatial order offers an alternative approach to stabilize the shock position by introducing more degrees of freedom to remodel the Rankine-Hugoniot condition at a captured shock. However, Tu et al. \cite{Tu_Evaluation_2014} perform a series of numerical experiments and have a different perspective. They point out that high-order schemes are at a higher risk of shock instabilities. The same conclusion is further supported by Jiang et al.\cite{Jiang_Effective_2017}, who find through numerical experiments that higher-order reconstruction leads to more severe shock anomalies, especially when the spatial accuracy is enhanced from third-order to fifth-order. Also, they find that higher-order schemes equipped with characteristic reconstruction will yield more stable results, although instabilities remain. To cure shock anomalies in high-order cases, Ohwada et al. \cite{Ohwada_Simple_2018a} propose a hybrid scheme. The main ingredients of this scheme involve variants of the standard numerical flux, MUSCL, and Lagrange’s polynomial. In smooth regions, the scheme achieves a fifth-order spatial accuracy, while in the vicinity of shocks, it operates at second-order. This hybrid scheme aims to balance robustness and low dissipation, guaranteeing stable computations and accurate heat flux results. Rodionov \cite{Rodionov_Artificial_2019,Rodionov_Simplified_2021} cures the shock instability in high-order schemes by introducing artificial viscosity. The artificial viscosity terms are in the form of right-hand sides of the Navier-Stokes equations and aim to smear shear waves within the shock layer only. Thus, this method is general and can be implemented to cure the shock instability of high-order schemes.

The characteristics and mechanisms of shock instabilities in high-order schemes require further investigation, as they form the foundation for developing robust high-order schemes. However, the research presented above indicates the limited focus on analyzing the shock instability of high-order schemes. This limitation is mainly attributed to the lack of effective analytical tools. Most studies investigating shock instability in high-order schemes rely heavily on numerical experiments. While numerical experiments can accurately predict whether shock instability problems will occur in actual simulations, they are not the most suitable tool for exploring the underlying mechanisms since they are time-consuming and involve the coupling of multiple factors that influence shock instabilities. It should be noted that Dumbser et al.\cite{Dumbser_Matrix_2004} propose the matrix stability analysis method. This method develops a model for predicting the evolution of perturbation errors and links the ability of schemes to stably capture shocks with the eigenvalues of the stability matrix. By utilizing this method, researchers can investigate the influence of factors such as Riemann solvers, computational grids, features of the numerical shock, and boundary conditions on shock instability separately. Additionally, the matrix stability analysis method is easy to implement by programming. As a result, it has been widely applied to investigate the characteristics and mechanisms of shock instabilities \cite{Dumbser_Matrix_2004,Shen_Stability_2014,Chen_Mechanism_2018,Chauvat_Shock_2005} and to evaluate the robustness of novel shock-stable schemes \cite{Chen_MechanismDerived_2018,Chen_Rotated_2023a,Chen_Novel_2018,Chen_Lowdiffusion_2021,Hu_Development_2023,Hu_Shockstable_2022,Sun_Effective_2022}.

However, the matrix stability analysis method proposed by Dumbser et al.\cite{Dumbser_Matrix_2004} is originally developed for first-order schemes, thus limiting its applicability in investigating the shock instability problem for higher-order schemes. In our previous work \cite{Ren_Numerical_2023}, we extend this method to analyze the shock instabilities of second-order schemes, specifically the second-order MUSCL approach. Based on the previous work of \cite{Dumbser_Matrix_2004} and our own extension \cite{Ren_Numerical_2023}, in the current study, the matrix stability analysis method of the fifth-order WENO scheme is established and employed to investigate the shock instabilities for high-order schemes. The subsequent sections of this paper are organized as follows: In section \ref{section 2}, the Euler equations of compressible flows and finite volume discretization are presented. The matrix stability analysis method for the fifth-order WENO scheme is established and validated in section \ref{section 3}. In section \ref{section 4}, this method is employed to investigate the shock instability of fifth-order schemes. Moreover, section \ref{section 5} contains conclusions and a discussion on the development of robust high-order schemes.

% Although the shock instability problem is first observed in Godunov-type schemes, it has been demonstrated that many shock-capturing methods for compressible flows will encounter the shock instability problem, such as the kinetic scheme and LBM

% various nuermical methods for compressible flows encounter shock instability: Godunov-type schemes, Kinetic scheme, GRP, LBM ,i.e. shock-capturing methods almost all fail.

% various governing models: typical problem for hyperbolic system, i.e., Euler, MHD, high temperature, shaollow water, and so on.

% background: hypersonic heating, detonation wave simulation, shock-shock interaction, shock-boudnary layer interaction, shock-disturbance interaction.DNS

% The severity of instability for high-order methods: even van Leer, Rusanov encounter shock instability (Ohawa).

% The overview of Research on high order instability, artificial viscosity, low Mach number, Jiang, Tu, Ohawa, shallow water ....Reza Zangeneh,

% The power of Matrix, i.e., used for Godunov-type schemes, even LBM

\section{Governing equations and finite-volume discretization}\label{section 2}

\subsection{The Euler equations}\label{subsection 2.1}

In the current work, the two-dimensional Euler equations in integral form are taken as the governing equations, which can be written as
\begin{equation}\label{eq Euler equations}
  \frac{\partial}{\partial t} \int_{\Omega} \mathbf{U} {\rm{d}}{\Omega} + \oint_{\partial \Omega} {\mathbf{F}} {\rm{d}}S = 0,
\end{equation}
where $\Omega$ denotes the control volume and $\partial \Omega$ is its boundary. $\mathbf{U}$ and $\mathbf{F}$ are the vectors of conservative variables and flux components normal to the boundary $\partial \Omega$, which are given respectively by
\begin{equation}\label{eq conservative variables and flux components}
	\mathbf{U} = \left[ {\begin{array}{{c}}
	\rho \\
	{\rho u}\\
	{\rho v}\\
	{\rho e}
	\end{array}} \right],
	\quad
	{\mathbf{F}} = \left[ {\begin{array}{{c}}
	\rho q \\
	{\rho q u + p n_x}\\
	{\rho q v + p n_y}\\	
	{\left(\rho e+p\right)q}
	\end{array}} \right].
\end{equation}
In (\ref{eq conservative variables and flux components}), $\rho$, $e$, and $p$ represent density, specific total energy, and pressure, respectively. The directed velocity, $q=un_x+vn_y$, is the component of velocity acting in the $\mathbf{n}$ direction, where $\mathbf{n}={\left[ {{n_x},\;{n_y}} \right]^T}$ is the outward unit vector normal to the surface element ${\rm{d}}S$ and $\left( {u,v} \right)$ is the flow velocity. In order to make the equations closed, the equation of state must be introduced, which can be expressed as
\begin{equation}\label{eq equation of state}
	p = \left(\gamma-1\right) \rho \left[e-\frac{1}{2}\left(u^2+v^2\right)\right],
\end{equation}
where $\gamma$ is the specific heat ratio.

\subsection{The high-order finite-volume WENO schemes}\label{subsection 2.2}
\begin{figure}[htbp]
	\centering
	\includegraphics[width=0.6\textwidth]{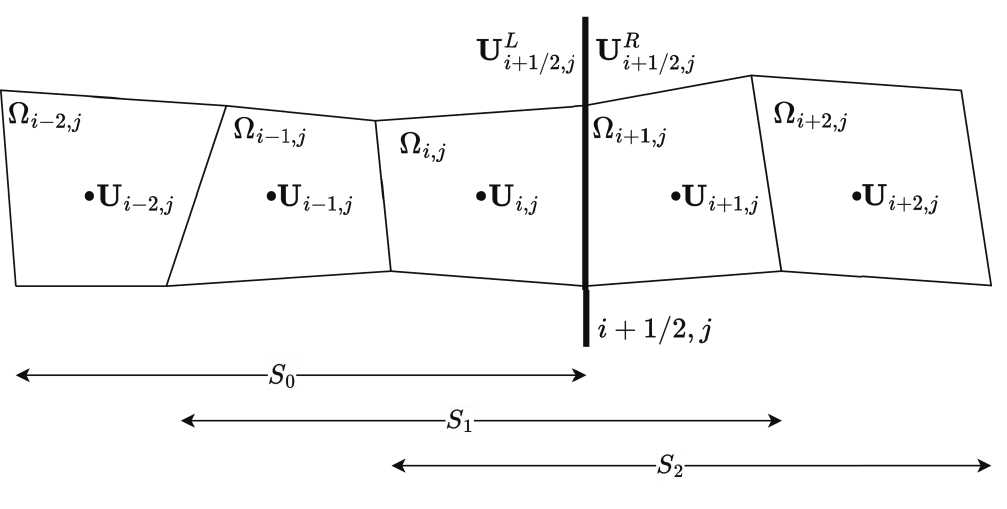}
	\caption{Candidate stencils $ S_k $ in the WENO method.}
	\label{fig stencil of WENO5}
\end{figure}

We consider discretizing (\ref{eq Euler equations}) with the cell-centered finite-volume method over a two-dimensional domain subdivided into some structured quadrilateral cells. The semi-discrete finite volume scheme can be written as
\begin{equation}\label{eq discrete Euler equations}
  \frac{\mathrm{d} {\mathbf{U}}_{i,j}}{\mathrm{~d} t}=-\frac{1}{\left|\Omega_{i,j}\right|} \sum_{k=1}^{4} \mathcal{L}_{k} \mathbf{F}_{k},
\end{equation}
where $ {\mathbf{U}}_{i,j} $ denotes the average of \textbf{U} on $ \Omega _{i,j} $. $ \left|\Omega_{i,j}\right| $ is the volume of $ \Omega _{i,j} $ and $\mathcal{L}_{k}$ stands for the length of the cell interface. $\mathbf{F}_{k}$ is the numerical flux that is supposed to be constant along the individual cell interface. In the current study, approximate Riemann solvers are used to determine $\mathbf{F}_{k}$ at each cell interface. As shown in Fig.\ref{fig stencil of WENO5}, considering the interface between $ \Omega_{i,j} $ and $ \Omega_{i+1,j} $, the numerical flux can be expressed in the following form
\begin{equation}\label{eq numerical flux}
  \mathbf{F}_{i+1/2,j}=\mathbf{F}\left(\mathbf{U}_{i+1/2,j}^L,\mathbf{U}_{i+1/2,j}^R \right),
\end{equation}
where $ \mathbf{U}_{i+1 / 2,j}^{L} $ and $ \mathbf{U}_{i+1 / 2,j}^{R} $ are the conservative variables on the left and right sides of the interface.

In the current work, we devote our efforts to investigating the shock instability problem of fifth-order finite-volume WENO schemes. Detailed information on the WENO method can be found in the references \cite{Jiang_Efficient_1996,Shu_Essentially_2020}, and only a brief description of the reconstruction process is given here. As shown in Fig.\ref{fig stencil of WENO5}, to obtain $ \mathbf{U}_{i+1/2}^L $, the fifth-order WENO method uses a global five-point stencil which is subdivided into three three-point substencils. $ \mathbf{U}_{i+1/2}^L $ can be written as
\begin{equation}\label{eq WENO method}
  \mathbf{U}_{i+1/2,j}^L = \boldsymbol{\omega}_{i+1/2,j}^{L,0}\mathbf{U}_{i+1 / 2,j}^{L,0}+\boldsymbol{\omega}_{i+1/2,j}^{L,1}\mathbf{U}_{i+1 / 2,j}^{L,1}+\boldsymbol{\omega}_{i+1/2,j}^{L,2}\mathbf{U}_{i+1 / 2,j}^{L,2},
\end{equation}
where $ \mathbf{U}_{i+1 / 2,j}^{L,m} \enspace (m={0,1,2}) $ is the third-order quadratic polynomial in substencil $ S_m $, which is given by
\begin{equation}
  \left\{\begin{aligned}
    \mathbf{U}_{i+1 / 2,j}^{L,0}&= \frac{1}{3} \mathbf{U}_{i-2,j}-\frac{7}{6} \mathbf{U}_{i-1,j}+\frac{11}{6}\mathbf{U}_{i,j} \\
    \mathbf{U}_{i+1 / 2,j}^{L,1}&=-\frac{1}{6} \mathbf{U}_{i-1,j}+\frac{5}{6} \mathbf{U}_{i,j}  +\frac{1}{3} \mathbf{U}_{i+1,j} \\
    \mathbf{U}_{i+1 / 2,j}^{L,2}&= \frac{1}{3} \mathbf{U}_{i,j}  +\frac{5}{6} \mathbf{U}_{i+1,j}-\frac{1}{6} \mathbf{U}_{i+2,j}
  \end{aligned}\right..
\end{equation}
In (\ref{eq WENO method}), $ \boldsymbol{\omega}_{i+1/2,j}^{L,m} \enspace (m={0,1,2}) $ is the diagonal matrix of the nonlinear weights, which can be defined as
\begin{equation}
  \boldsymbol{\omega} _{i + 1/2,j}^{L,m} = \text{diag}\left[{\left( {\omega _{i + 1/2,j}^{L,m}} \right)}_\rho,{\left( {\omega _{i + 1/2,j}^{L,m}} \right)}_{\rho u},{\left( {\omega _{i + 1/2,j}^{L,m}} \right)}_{\rho v},{\left( {\omega _{i + 1/2,j}^{L,m}} \right)}_{\rho e}\right] ,\quad m={0,1,2},
\end{equation}
where the subscripts indicate that the nonlinear weights are calculated by different conservative variables. Since the nonlinear weights have the same calculation procedure for different variables, for simplicity, we only present the calculation of the nonlinear weights corresponding to $ \rho $, which are defined as
\begin{equation}
  {\left( {\omega _{i + 1/2,j}^{L,m}} \right)}_\rho=\frac{\alpha_m }{\sum_{l=0}^{2}\alpha_l}, \quad \alpha_m=\frac{d_m}{(\beta_m+\epsilon)^2},\quad m={0,1,2}.
\end{equation}
The linear weights $ d_m $ are
\begin{equation}
  d_0=\frac{1}{10},d_0=\frac{3}{5},d_0=\frac{3}{10}.
\end{equation}
The parameter $ \epsilon = 10^{-15} $ is used to avoid the division by zero in the denominator. $ \beta_m $ are called the smoothness indicators and given by
\begin{equation}
  \left\{\begin{aligned}
    \beta_0&=\frac{13}{12}\left(\rho_{i-2,j}-2 \rho_{i-1,j}+\rho_{i,j}\right)^2+\frac{1}{4}\left(\rho_{i-2,j}-4 \rho_{i-1,j}+3 \rho_i\right)^2 \\
    \beta_1&=\frac{13}{12}\left(\rho_{i-1,j}-2 \rho_{i,j}+\rho_{i+1,j}\right)^2+\frac{1}{4}\left(\rho_{i-1,j}-\rho_{i+1,j}\right)^2 \\
    \beta_2&=\frac{13}{12}\left(\rho_{i,j}-2 \rho_{i+1,j}+\rho_{i+2,j}\right)^2+\frac{1}{4}\left(3 \rho_{i,j}-4 \rho_{i+1,j}+\rho_{i+2,j}\right)^2
    \end{aligned}\right..
\end{equation}

The procedure described above in this section is the WENO-JS method for reconstructing the conservative variables. Based on the idea of WENO-JS, Borges et al.\cite{Borges_Improved_2008} propose the WENO-Z method to improve the order of accuracy at critical points. For the WENO-Z method, the smoothness indicators are given by
\begin{equation}
  \beta_m^z=\left(\frac{\beta_m+\epsilon}{\beta_m+\tau_5+\epsilon}\right),\quad m={0,1,2},
\end{equation}
where $ \tau_5=|\beta_0-\beta_2| $. So, the new nonlinear weights are
\begin{equation}
  \omega_m^z=\frac{a_m^z}{\sum_{l=0}^2\alpha_l^z},\quad a_m^z=\frac{d_m}{\beta_m^z}=d_k\left(1+\frac{\tau_5}{\beta_m+\epsilon}\right),\quad m={0,1,2}.
\end{equation}
In this paper, the WENO-Z method is employed to reconstruct the variables on the left and right sides of the interface. And we only introduce the reconstruction of $ \mathbf{U}_{i+1/2,j}^L $ in this section. The reconstruction of $ \mathbf{U}_{i+1/2,j}^R $ can also be obtained by symmetry.

% \section{The planar steady shock instability for high-order methods}

\section{A matrix stability analysis method for the WENO scheme}\label{section 3}
The matrix stability analysis method has been demonstrated to be an efficient tool for investigating the shock instability problem \cite{Dumbser_Matrix_2004,Shen_Stability_2014,Simon_Cure_2018}. However, the matrix stability analysis method developed in \cite{Dumbser_Matrix_2004} only applies to the first-order scheme, limiting its application to the higher-order schemes. In \cite{Ren_Numerical_2023}, the matrix stability analysis method for the finite-volume MUSCL approach is proposed, and some unique stability characteristics for second-order shock-capturing methods are revealed quantitatively for the first time. Compared with its low-order counterparts, the establishment of the matrix stability analysis method for high-order shock-capturing schemes is more challenging and needs to be dealt with carefully. In the following sections, the implementation of the stability analysis method for fifth-order finite-volume schemes is presented in detail.

\subsection{Derivation of the matrix stability method for fifth-order WENO schemes}\label{subsection 3.2}

The basic idea of the matrix stability analysis method is to study the nature of the temporal evolution of random perturbations, which are introduced into the steady mean values of the flow field at the initial time. This section provides a brief introduction to the main progress of the matrix stability analysis method in \cite{Dumbser_Matrix_2004}, and then develop the matrix stability analysis method for fifth-order WENO schemes based on their characteristics.

For the steady field, we assume that
\begin{equation}\label{eq variables decomposition}
	\mathbf{U}_{i,j}=\bar{\mathbf{U}}_{i,j} + \delta\mathbf{U}_{i,j},
\end{equation}
where $ \bar{\mathbf{U}}_{i,j} $ denotes the steady mean value and $ \delta \mathbf{U}_{i,j} $ is the small numerical random perturbation. Then the flux function of the face $ i+1/2,j $ can be linearized around the steady mean value as follows
\begin{equation}   
	\begin{aligned}
		\mathbf{F}_{i+1/2,j} &= \mathbf{F}_{i+1/2,j}\left(\mathbf{U}_{i,j},\mathbf{U}_{i+1,j}\right)\\
		&= \mathbf{F}_{i+1/2,j}\left(\bar{{\mathbf{U}}}_{i+1/2,j}^{L},\bar{{\mathbf{U}}}_{i+1/2,j}^{R}\right)+\frac{\partial\mathbf{F}_{i+1/2,j}}{\partial\mathbf{U}_{i,j}} \delta \mathbf{U}_{i,j}+\frac{\partial\mathbf{F}_{i+1/2,j}}{\partial\mathbf{U}_{i+1,j}} \delta \mathbf{U}_{i+1,j}
	\end{aligned}.
	\label{eq first-order flux linearized}
\end{equation}
Then substituting (\ref{eq variables decomposition}) and (\ref{eq first-order flux linearized}) into (\ref{eq discrete Euler equations}) and taking into account that the mean-field is steady, the linear error evolution model can be obtained
\begin{equation}   
	\begin{aligned}\label{eq first-order error evolution}
		\frac{\mathrm{d} \delta \mathbf{U}_{i, j}}{\mathrm{dt}} = &-\left[\sigma_{i+1/2,j}\frac{\partial\mathbf{F}_{i+1/2,j}}{\partial\mathbf{U}_{i,j}}+\sigma_{i,j+1/2}\frac{\partial\mathbf{F}_{i,j+1/2}}{\partial\mathbf{U}_{i,j}}+\sigma_{i-1/2,j}\frac{\partial\mathbf{F}_{i-1/2,j}}{\partial\mathbf{U}_{i,j}}+\sigma_{i,j-1/2}\frac{\partial\mathbf{F}_{i,j-1/2}}{\partial\mathbf{U}_{i,j}} \right]\delta \mathbf{U}_{i,j}\\
		= &-\sigma_{i+1/2,j}\frac{\partial\mathbf{F}_{i+1/2,j}}{\partial\mathbf{U}_{i+1,j}}\delta\mathbf{U}_{i+1,j}-\sigma_{i,j+1/2}\frac{\partial\mathbf{F}_{i,j+1/2}}{\partial\mathbf{U}_{i,j+1}}\delta\mathbf{U}_{i,j+1}\\
		&-\sigma_{i-1/2,j}\frac{\partial\mathbf{F}_{i-1/2,j}}{\partial\mathbf{U}_{i-1,j}}\delta\mathbf{U}_{i-1,j}-\sigma_{i,j-1/2}\frac{\partial\mathbf{F}_{i,j-1/2}}{\partial\mathbf{U}_{i,j-1}}\delta\mathbf{U}_{i,j-1}
	  \end{aligned},
\end{equation}
where 
\begin{equation}   
	\sigma_{i\pm1/2,j\pm1/2} = \frac{\mathcal{L}_{i\pm1/2,j\pm1/2}}{\Omega_{i,j}}.
\end{equation}
(\ref{eq first-order error evolution}) is the linear error evolution model proposed in \cite{Dumbser_Matrix_2004}. By comparing (\ref{eq numerical flux}) and (\ref{eq first-order flux linearized}), it can be found that (\ref{eq first-order error evolution}) is only applicable to first-order schemes, since it assumes 
\begin{equation}   
	\mathbf{U}_{i+1/2,j}^L = \mathbf{U}_{i,j} \quad \text{and} \quad \mathbf{U}_{i+1/2,j}^R = \mathbf{U}_{i+1,j}.
\end{equation}
As a result, in order to estabilsh the linear error evolution model for high-order schemes, the influence of reconstruction must be considered.

\begin{figure}[htbp]
	\centering
	\includegraphics[width=0.99\textwidth]{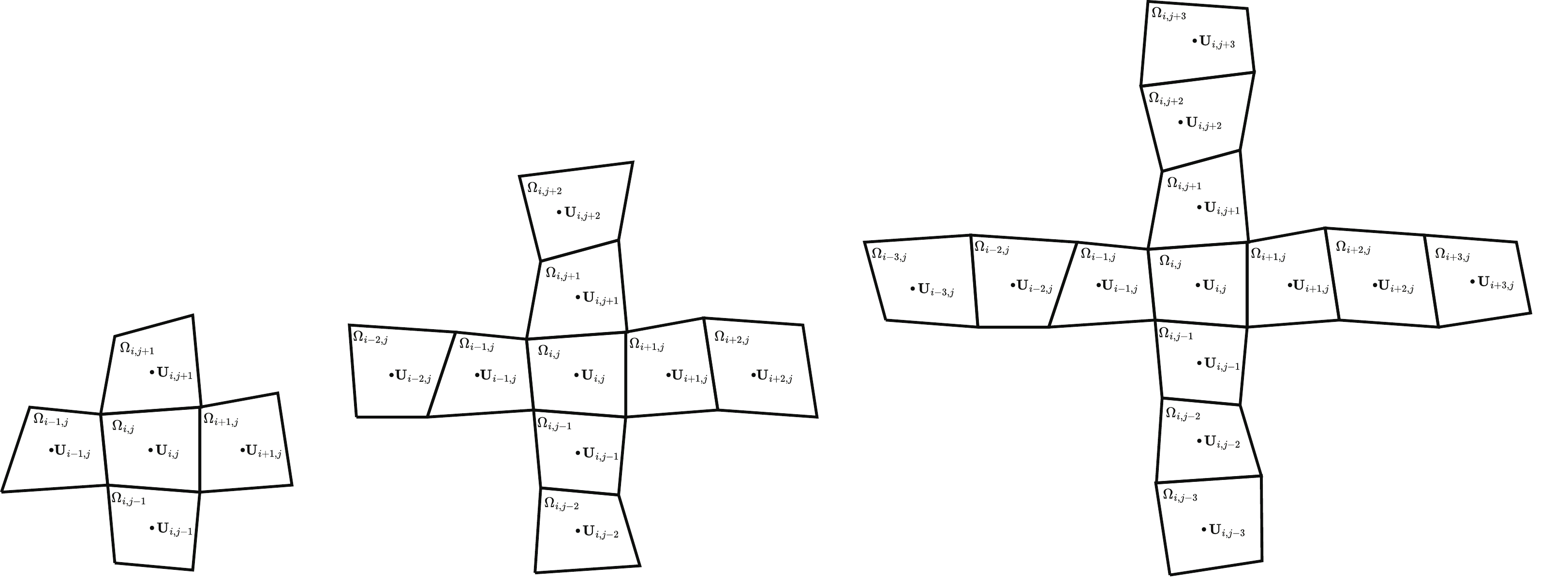}
	\caption{The cells that influence the evolution of $ \delta\mathbf{U}_{i,j} $ for different spatial accuracy.(From the left to the right is for first-order scheme, second-order schemes, and fifth-order schemes, respectively.)}
	\label{fig stencil compare}
\end{figure}

Substituting (\ref{eq variables decomposition}) into the WENO reconstruction of $ \mathbf{U}_{i+1/2,j}^L $ and $ \mathbf{U}_{i+1/2,j}^R $, it can be gotten that
\begin{equation}  \label{eq fifth-order variables decomposition}
  \begin{aligned}
    \delta \mathbf{U}_{i+1/2,j}^L &=\frac{1}{3} \boldsymbol{\omega}_{i+1/2,j}^{L,0} \delta\mathbf{U}_{i-2,j}-\frac{1}{6}\left(7\boldsymbol{\omega}_{i+1/2,j}^{L,0}+\boldsymbol{\omega}_{i+1/2,j}^{L,1}\right)\delta\mathbf{U}_{i-1,j}\\
    &+\frac{1}{6}\left(11\boldsymbol{\omega}_{i+1/2,j}^{L,0}+5\boldsymbol{\omega}_{i+1/2,j}^{L,1}+2\boldsymbol{\omega}_{i+1/2,j}^{L,2}\right)\delta\mathbf{U}_{i,j}\\
    &+\frac{1}{6}\left(2\boldsymbol{\omega}_{i+1/2,j}^{L,1}+5\boldsymbol{\omega}_{i+1/2,j}^{L,2} \right)\delta\mathbf{U}_{i+1,j}-\frac{1}{6}\boldsymbol{\omega}_{i+1/2,j}^{L,2}\delta\mathbf{U}_{i+2,j},\\
    \delta \mathbf{U}_{i+1/2}^R&=\frac{1}{3} \boldsymbol{\omega}_{i+1/2,j}^{R,0} \delta\mathbf{U}_{i+3,j}-\frac{1}{6}\left(7\boldsymbol{\omega}_{i+1/2,j}^{R,0}+\boldsymbol{\omega}_{i+1/2,j}^{R,1}\right)\delta\mathbf{U}_{i+2,j}\\
    &+\frac{1}{6}\left(11\boldsymbol{\omega}_{i+1/2,j}^{R,0}+5\boldsymbol{\omega}_{i+1/2,j}^{R,1}+2\boldsymbol{\omega}_{i+1/2,j}^{R,2}\right)\delta\mathbf{U}_{i+1,j}\\
    &+\frac{1}{6}\left(2\boldsymbol{\omega}_{i+1/2,j}^{R,1}+5\boldsymbol{\omega}_{i+1/2,j}^{R,2} \right)\delta\mathbf{U}_{i,j}-\frac{1}{6}\boldsymbol{\omega}_{i+1/2,j}^{R,2}\delta\mathbf{U}_{i-1,j}.
  \end{aligned}
\end{equation}
Note that in (\ref{eq fifth-order variables decomposition}), we assume that the perturbations will not influence the computation of the nonlinear weights $ \boldsymbol{\omega}_{i+1/2,j}^{L/R,m} $. As shown in (\ref{eq numerical flux}), the numerical flux is the function of $ \mathbf{U}_{i+1/2,j}^L $ and $ \mathbf{U}_{i+1/2,j}^R $. Therefore, $ \mathbf{F}_{i+1/2,j} $ can be linearized around the steady mean value as
\begin{equation}
  \begin{aligned}\label{eq simplify second-order flux linearized}
    \mathbf{F}_{i+1/2,j}& =\mathbf{F}_{i+1/2,j}\left(\bar{{\mathbf{U}}}_{i+1/2,j}^{L},\bar{{\mathbf{U}}}_{i+1/2,j}^{R}\right)  \\
    &+\boldsymbol{\alpha}_{i+1/2,j}^{L}\delta\mathbf{U}_{i-2,j}+\boldsymbol{\beta}_{i+1/2,j}^{L}\delta\mathbf{U}_{i-1,j}+\boldsymbol{\chi}_{i+1/2,j}^{L}\delta\mathbf{U}_{i,j} \\
    &+\boldsymbol{\chi}_{i+1/2,j}^{R}\delta\mathbf{U}_{i+1,j}+\boldsymbol{\beta}_{i+1/2,j}^{R}\delta\mathbf{U}_{i+2,j}+\boldsymbol{\alpha}_{i+1/2,j}^{R}\delta\mathbf{U}_{i+3,j}
    \end{aligned},
\end{equation}
where
\begin{equation}\label{eq alpha beta chi}
  \begin{aligned}
    \boldsymbol{\alpha}_{i+1/2,j}^{L/R}&=\frac{1}{3} \frac{\partial\mathbf{F}_{i+1/2,j}}{\partial\mathbf{U}_{i+1/2,j}^{L/R}} \boldsymbol{\omega}_{i+1/2,j}^{L/R,0},\\
    \boldsymbol{\beta}_{i+1/2,j}^{L/R}&=-\frac{1}{6}\left[\frac{\partial\mathbf{F}_{i+1/2,j}}{\partial\mathbf{U}_{i+1/2,j}^{L/R}}\left(7\boldsymbol{\omega}_{i+1/2,j}^{L/R,0}+\boldsymbol{\omega}_{i+1/2,j}^{L/R,1}\right)+\frac{\partial\mathbf{F}_{i+1/2,j}}{\partial\mathbf{U}_{i+1/2,j}^{R/L}}\boldsymbol{\omega}_{i+1/2,j}^{R/L,2}\right],\\
    \boldsymbol{\chi}_{i+1/2,j}^{L/R}&=\frac{1}{6}\left[\frac{\partial\mathbf{F}_{i+1/2,j}}{\partial\mathbf{U}_{i+1/2,j}^{L/R}}\left(11\boldsymbol{\omega}_{i+1/2,j}^{L/R,0}+5\boldsymbol{\omega}_{i+1/2,j}^{L/R,1}+2\boldsymbol{\omega}_{i+1/2,j}^{L/R,2}\right)\right.\\
    &\quad\;+\left.\frac{\partial\mathbf{F}_{i+1/2,j}}{\partial\mathbf{U}_{i+1/2,j}^{R/L}}\left(2\boldsymbol{\omega}_{i+1/2,j}^{R/L,1}+5\boldsymbol{\omega}_{i+1/2,j}^{R/L,2} \right)\right].
  \end{aligned}
\end{equation}
Substituting (\ref{eq variables decomposition}) and (\ref{eq simplify second-order flux linearized}) into (\ref{eq discrete Euler equations}), we can get the linear error evolution model of $ \Omega_{i,j} $
\begin{equation}\label{eq linear error evolution model for conservative variables}
  \begin{aligned}
    \frac{\mathrm{d} \delta \mathbf{U}_{i, j}}{\mathrm{dt}}=&\left(\boldsymbol{\zeta}_{i+1/2,j}^{L}+\boldsymbol{\zeta}_{i,j+1/2}^{L}+\boldsymbol{\zeta}_{i-1/2,j}^{R}+\boldsymbol{\zeta}_{i,j-1/2}^{R}\right)\delta\mathbf{U}_{i,j}\\
    +&\left(\boldsymbol{\zeta}_{i+1/2,j}^R+\boldsymbol{\eta}_{i-1/2,j}^R\right)\delta\mathbf{U}_{i+1,j}
    +\left(\boldsymbol{\zeta}_{i,j+1/2}^R+\boldsymbol{\eta}_{i,j-1/2}^R\right)\delta\mathbf{U}_{i,j+1}\\
    +&\left(\boldsymbol{\zeta}_{i-1/2,j}^L+\boldsymbol{\eta}_{i+1/2,j}^L\right)\delta\mathbf{U}_{i-1,j}
    +\left(\boldsymbol{\zeta}_{i,j-1/2}^L+\boldsymbol{\eta}_{i,j+1/2}^L\right)\delta\mathbf{U}_{i,j-1}\\
    +&\left(\boldsymbol{\eta}_{i+1/2,j}^R+\boldsymbol{\theta}_{i-1/2,j}^R\right)\delta\mathbf{U}_{i+2,j}
    +\left(\boldsymbol{\eta}_{i,j+1/2}^R+\boldsymbol{\theta}_{i,j-1/2}^R\right)\delta\mathbf{U}_{i,j+2}\\
    +&\left(\boldsymbol{\eta}_{i-1/2,j}^L+\boldsymbol{\theta}_{i+1/2,j}^L\right)\delta\mathbf{U}_{i-2,j}
    +\left(\boldsymbol{\eta}_{i,j-1/2}^L+\boldsymbol{\theta}_{i,j+1/2}^L\right)\delta\mathbf{U}_{i,j-2}\\
    +&\boldsymbol{\theta}_{i+1/2,j}^R\delta\mathbf{U}_{i+3,j}+\boldsymbol{\theta}_{i,j+1/2}^R\mathbf{U}_{i,j+3}+\boldsymbol{\theta}_{i-1/2,j}^L\mathbf{U}_{i-3,j}+\boldsymbol{\theta}_{i,j-1/2}^L\delta\mathbf{U}_{i,j-3}
  \end{aligned},
\end{equation}
where
\begin{equation}
  \begin{aligned}\label{eq zeta eta theta}
    \boldsymbol{\zeta}_{i\pm1/2,j\pm1/2}^{L/R} &= -\frac{\mathcal{L}_{i\pm1/2,j\pm1/2}}{\Omega_{i,j}}\boldsymbol{\chi}_{i\pm1/2,j\pm1/2}^{L/R}\\
    \boldsymbol{\eta}_{i\pm1/2,j\pm1/2}^{L/R} &= -\frac{\mathcal{L}_{i\pm1/2,j\pm1/2}}{\Omega_{i,j}}\boldsymbol{\beta}_{i\pm1/2,j\pm1/2}^{L/R}\\
    \boldsymbol{\theta}_{i\pm1/2,j\pm1/2}^{L/R} &= -\frac{\mathcal{L}_{i\pm1/2,j\pm1/2}}{\Omega_{i,j}}\boldsymbol{\alpha}_{i\pm1/2,j\pm1/2}^{L/R}\\
  \end{aligned}.
\end{equation}
As shown in Fig.\ref{fig stencil compare}, it can be found from (\ref{eq linear error evolution model for conservative variables}) that the evolution of perturbations within $ \Omega_{i,j} $ in fifth-order schemes is influenced by the perturbation errors in more cells than the first- \cite{Dumbser_Matrix_2004} and second-order \cite{Ren_Numerical_2023} cases. The difference in the linear error evolution model suggests that the shock stability of the fifth-order scheme may be different from that of the first- and second-order schemes.
%（这里可以画出示意图）

Equation (\ref{eq linear error evolution model for conservative variables}) is applicable to all control volumes in the computational domain. As a result, the perturbation error evolution of all cells is as follows
\begin{equation}\label{eq linear error evolution model for all cells}
  \frac{\mathrm{d}}{\mathrm{d} t}\left(\begin{array}{c}
		\delta \mathbf{U}_{1,1} \\
		\vdots \\
		\delta \mathbf{U}_{imax,jmax}
		\end{array}\right)=\mathbf{S} \cdot\left(\begin{array}{c}
		\delta \mathbf{U}_{1,1} \\
		\vdots \\
		\delta \mathbf{U}_{imax,jmax}
		\end{array}\right),
\end{equation}
where \textbf{S} denotes the stability matrix in the present study. When considering only the evolution of initial perturbations, the solution of the linear time-invariant system (\ref{eq linear error evolution model for all cells}) is
\begin{equation}\label{eq solution of linear error evolution model}
	\left(\begin{array}{c}
		\delta \mathbf{U}_{1,1} \\
		\vdots \\
		\delta \mathbf{U}_{imax,jmax}
		\end{array}\right)(t)=\mathrm{e}^{\mathbf{S} t} \cdot\left(\begin{array}{c}
		\delta \mathbf{U}_{1,1} \\
		\vdots \\
		\delta \mathbf{U}_{imax,jmax}
	\end{array}\right)_{t=0}.
\end{equation}
All the perturbations will remain bounded if the maximum of the real part of the eigenvalues of \textbf{S} is negative. So the stability criterion is
\begin{equation}\label{eq stability criterion}
	\max (\operatorname{Re}(\lambda(\mathbf{S}))) \leq 0.
\end{equation}

%（这里表述不好）
Note that (\ref{eq simplify second-order flux linearized}) is accurate only when the flux function is differentiable at the mean value, which is not always holding and will seriously affect the analysis. As mentioned in \cite{Dumbser_Matrix_2004}, when the shock is not exactly between two cells ($ \varepsilon \neq 0 \enspace \text{and} \enspace \varepsilon \neq 1$), the non-differentiable solvers, such as Roe, will be differentiable. To Furtherly enhance the reliability of matrix stability analysis method, a smoothing procedure proposed in \cite{Shen_Stability_2014} is also employed
\begin{equation}   
	| \lambda | =\left\{\begin{array}{ll}|\lambda|,&\text{if~}|\lambda|\geq\delta_0,\\(\lambda^2+\delta_0^2)/(2\delta_0),&\text{if~}|\lambda|<\delta_0,\end{array}\right.
	\label{eq smoothing procedure}
\end{equation}
where $ \delta_0=10^{-4} $. $ \lambda $ denotes the eigenvalues of the Jacobian matrix and can be expressed as
\begin{equation}   
	\lambda_1 = \hat{q}-\hat{c},\quad
	\lambda_2 = \hat{q},\quad
	\lambda_2 = \hat{q},\quad
	\lambda_4 = \hat{q}+\hat{c}.
\end{equation} 
where $ \hat{(\cdot)} $ is the Roe averaged value.

The preceding paragraphs in this section describe the matrix stability analysis method for the fifth-order scheme, which is based on conservative variables. It is noted that a similar matrix stability analysis can also be performed for the primitive variables as well as the characteristic variables, which can be found in \ref{Appendix A} and \ref{Appendix B}, respectively.

\subsection{The planar steady shock instability for high-order schemes}
\begin{figure}[htbp]
	\centering
	\includegraphics[width=0.4\textwidth]{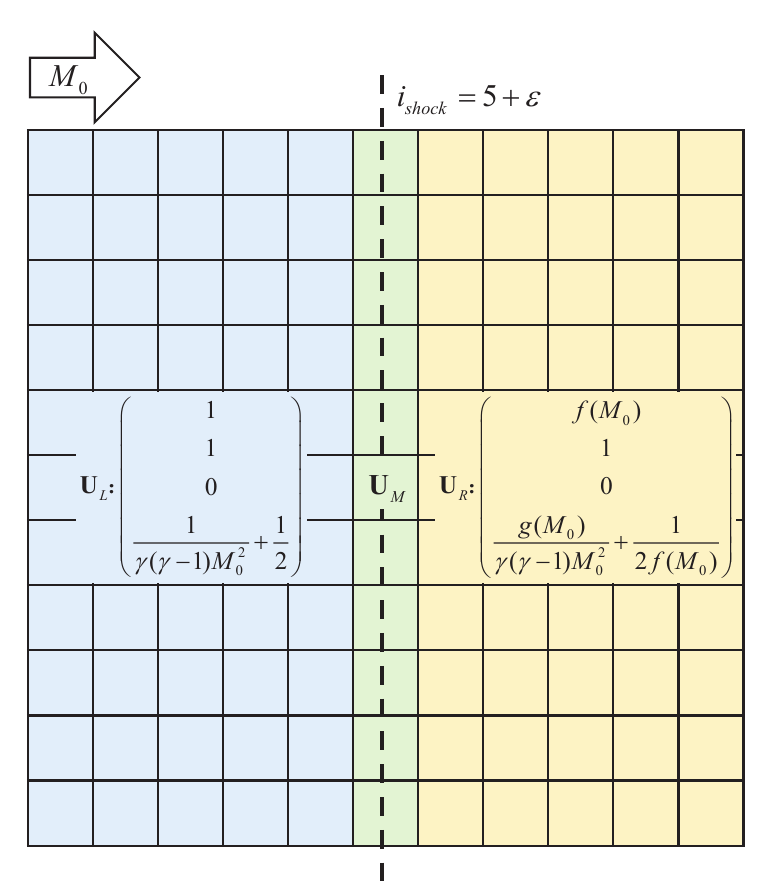}
	\caption{Set of the 2D normal shock problem.}
	\label{fig 2Dshock}
\end{figure}

The matrix stability analysis method proposed in \ref{subsection 3.2} can be employed to analyze the shock instability problem of all 2D flows with strong shocks. The carbuncle phenomenon is conventionally referred to as the distorted shock ahead of the blunt-body in supersonic or hypersonic flow \cite{Quirk_Contribution_1994}. However, in order to analyze the mechanism of shock instabilities, it is necessary to choose an appropriate test case. In the current work, the 2D steady normal shock problem is employed, which is simple enough and shares the fundamental characteristics of the carbuncle phenomenon. It has been well demonstrated that if a scheme yields unstable solutions for the steady normal shock problem, it will also suffer from the blunt-body carbuncle \cite{Dumbser_Matrix_2004, Ismail_Reliable_2006, Kitamura_Evaluation_2009}. The set of the 2D normal shock problem is shown in Fig.\ref{fig 2Dshock}, where $ M_0 $ is the upstream Mach number and 
\begin{equation}  
	f\left(M_0\right)=\left(\frac{2}{(\gamma+1) M_0^2}+\frac{\gamma-1}{\gamma+1}\right)^{-1}, \quad g\left(M_0\right)=\frac{2 \gamma M_0^2}{\gamma+1}-\frac{\gamma-1}{\gamma+1}.
\end{equation}
$ \mathbf{U}_M $ is the intermediate state within the shock ($ M: i=6 $) and set according to the Hugoniot curve
\begin{equation}   \label{eq numerical shock structure}
	\mathbf{U}_M=\left(\begin{array}{c}
		         \rho_M \\
		         u_M \\
		         0 \\
		         p_M
		        \end{array}\right)\quad \text{with} \quad 
	\left\{\begin{gathered}
		\rho_M=\left(1-\alpha_\rho\right) \rho_L+\alpha_\rho \rho_R \\
		u_M=\left(1-\alpha_u\right) u_L+\alpha_u u_R \\
		p_M=\left(1-\alpha_p\right) p_L+\alpha_p p_R
	\end{gathered}\right.,
\end{equation}
where
\begin{equation}   \label{eq epsilon}
	\begin{aligned}
		& \alpha_\rho=\varepsilon \\
		& \alpha_u=1-(1-\varepsilon)\left(1+\varepsilon \frac{M_0^2-1}{1+(\gamma-1) M_0^2 / 2}\right)^{-1 / 2}\left(1+\varepsilon \frac{M_0^2-1}{1-2 \gamma M_0^2 /(\gamma-1)}\right)^{-1 / 2}\\
		& \alpha_p=\varepsilon\left[1+(1-\varepsilon) \frac{\gamma+1}{\gamma-1} \frac{M_0^2-1}{M_0^2}\right]^{-1 / 2}
	\end{aligned}.
\end{equation}
$ \varepsilon \in  [0,1]$ is a weighting average that describes the initial state of the internal cell and is called shock position here. Detailed descriptions of boundary conditions and carbuncle phenomenon in the 2D normal shock problem can be found in \cite{Ren_Numerical_2023}, which are not introduced in this paper for the sake of simplicity.

As shown in (\ref{eq variables decomposition}), the matrix stability analysis method is based on the convergent and stable solution. In the current study, the 2D flow field for stability analysis is initialized by projecting the steady flow field from 1D computation onto the 2D domain, which is also employed in \cite{Dumbser_Matrix_2004, Ren_Numerical_2023, Sanders_Multidimensional_1998}. The detailed implementation of this initialization method can be found in our previous work in \cite{Ren_Numerical_2023}.

\subsection{Quantitative validation of the stability theory}
\begin{figure}[htbp]
	\centering
	\includegraphics[width=0.6\textwidth]{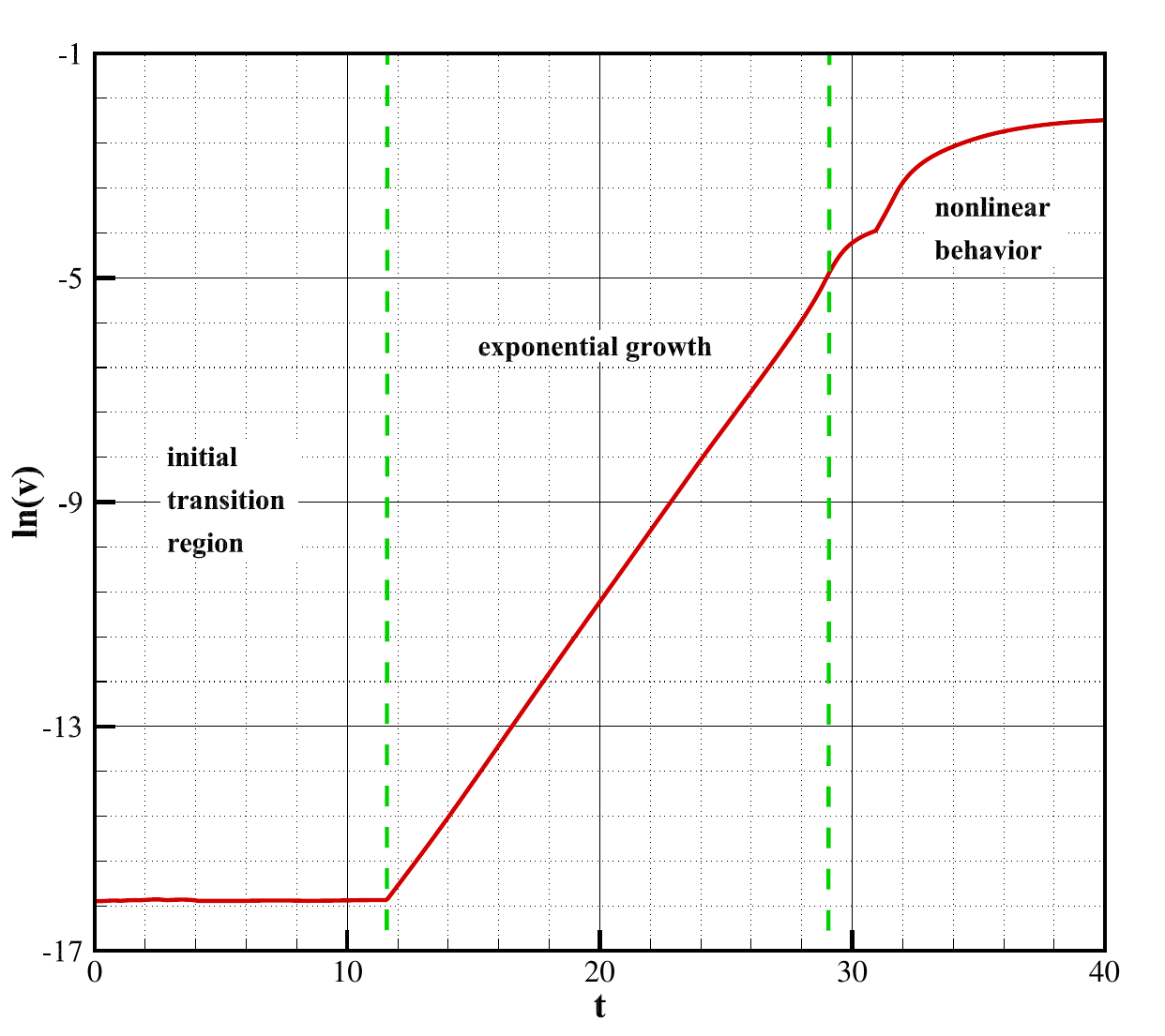}
	\caption{Evolution of the perturbation error in numerical experiment.(Grid with 11$ \times $11 cells, fifth-order scheme with Roe solver, $ M_0=20 $ and $ \varepsilon=0.1 $.)}
	\label{fig dv}
\end{figure}

Fig.\ref{fig dv} shows the evolution of the perturbation error in the numerical experiment of the 2D steady normal shock problem. In this paper, we take the norm $\|v\|_{\infty}(t)$ of the transverse velocity as the indicator of perturbation error since it is known that the exact solution for mean flows is $v = 0$ \cite{Dumbser_Matrix_2004, Henderson_Grid_2007}. The computation is performed by the fifth-order scheme with the Roe solver. An initial perturbation of size $10^{-7}$ is added to each cell at the beginning of the computation. Note that for the 2D steady normal shock problem employed throughout this paper, the computational grid is an $11\times11$ Cartesian grid, and the conditions are $M_0 = 20$ and $\varepsilon = 0.1$. The CFL number is 0.1 if not mentioned specifically.

\begin{figure}[htbp]
	\centering
	\includegraphics[width=0.6\textwidth]{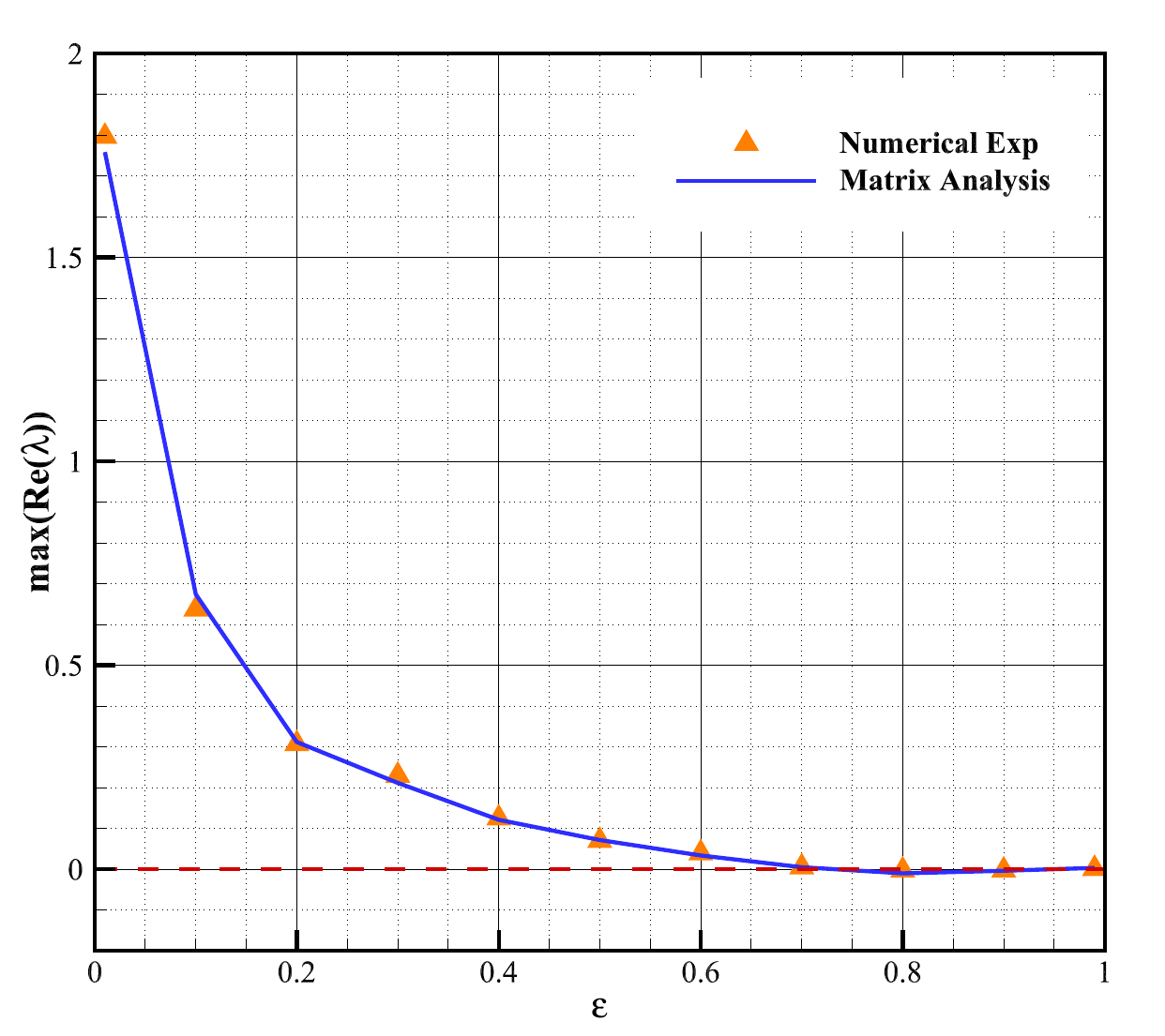}
	\caption{Quantitative validation of the matrix stability analysis method.(Grid with 11$ \times $11 cells, fifth-order scheme with Roe solver, $ M_0=20 $ and $ \varepsilon=\left\{0.01,0.1,0.2,\dots,0.8,0.9,0.99\right\} $.)}
	\label{fig validation}
\end{figure}

As shown in Fig.\ref{fig dv}, the evolution of the perturbation can be divided into three stages: the initial transition region, exponential growth stage, and nonlinear behavior stage. A detailed description of these stages can be found in \cite{Dumbser_Matrix_2004, Ren_Numerical_2023}. Among them, we are particularly interested in the exponential growth stage. This stage plays a vital role in the evolution of the perturbation error as it determines whether the error will increase or decrease, further stable or unstable, and how quickly it will develop towards instability. The evolution of the perturbation error in the exponential growth stage can be described by the equation:
\begin{equation}\label{eq exponential law}
  \|v\|_{\infty}(t)=v_{0} \mathrm{e}^{\lambda_{num}\left(t-t_{0}\right)},
\end{equation}
where $ \lambda_{num} $ is the temporal error growth rate in the present work. The values of $ v_{0} $, $ \lambda_{num} $, and $ t_{0} $ can be easily obtained from Fig.\ref{fig dv}. By analyzing both (\ref{eq solution of linear error evolution model}) and (\ref{eq exponential law}), it can be observed that the maximal real part of the eigenvalues can describe the temporal error growth rate of the perturbation error. This provides a reliable method to quantitatively validate the accuracy of the matrix stability analysis method proposed in section \ref{subsection 3.2}. The comparison is shown in Fig.\ref{fig validation}, where good agreements can be observed between them in both the unstable and stable regions. Additionally, similar good agreements can be obtained using other fifth-order schemes with different solvers and computational conditions, which further confirms the reliability of the matrix stability analysis method developed in the current study.

Furthermore, based on Fig.\ref{fig dv} and equation (\ref{eq solution of linear error evolution model}), two important points need to be emphasized. Firstly, $\max(\operatorname{Re}(\lambda(\mathbf{S})))$ can be used to predict the exponential growth of the perturbation error, rather than the magnitude of the perturbation error in the finial flow field. Thus, this method is supposed to predict the trend (increase or decrease) and speed of the perturbation error evolution, providing insights into whether the calculation will be stable or unstable and how quickly it will develop. Secondly, the matrix stability analysis is based on linear analysis and the initial perturbation. It can be deduced that such a linear analysis may not be valid in the nonlinear stage where the perturbations become too large to be analyzed linearly. Nonetheless, as mentioned in \cite{Shen_Stability_2014}, the matrix stability analysis method is a potent tool that offers quantitative insights into the behavior of numerical schemes and helps illustrate the mechanisms underlying shock instability problems. If the analysis suggests that the computation is unstable, it is always possible to detect the instability, although such unstable phenomena only manifest in specific numerical cases. It is important to note that when referring to instability in the remainder of this paper, we specifically signify the occurrence of positive exponential growth in the perturbation error during the exponential growth stage.

\section{Stability analysis of high-order finite-volume schemes}\label{section 4}

The shock instability problem is conventionally associated with the dissipation characteristics of shock-capturing schemes. When using dissipative solvers such as HLL and van Leer, the scheme can stably capture strong shocks in most cases. Whereas, for low-dissipative solvers, such as Roe and HLLC, the computation is more prone to instability. This has been well demonstrated in previous investigations for first- and second-order schemes \cite{Ren_Numerical_2023, Dumbser_Matrix_2004, Quirk_Contribution_1994, Kitamura_Evaluation_2009, Kitamura_Evaluation_2010, Pandolfi_Numerical_2001}. The characteristics and mechanisms of shock instability in low-order schemes have been extensively studied. However, when the spatial accuracy is increased to the fifth-order, it raises questions about whether the inherent dissipation of dissipative solvers is still sufficient to stably capture strong shocks. Moreover, it is well known that high-order schemes are more vulnerable to shock instabilities. The characteristics and mechanisms of instability problem for strong shocks need to be investigated in depth.

\subsection{Stability results of fifth-order schemes with different solvers}\label{subsection 4.1}

% \begin{figure}[htbp]
% 	\centering

% 	\subfigure[first-order scheme]{
% 	\begin{minipage}[t]{0.46\linewidth}
% 	\centering
% 	\includegraphics[width=0.9\textwidth]{5_1_HLLC_1nd_matrix.pdf}
% 	\end{minipage}
% 	}
% 	\subfigure[second-order scheme]{
% 	\begin{minipage}[t]{0.46\linewidth}
% 	\centering
% 	\includegraphics[width=0.9\textwidth]{5_1_HLLC_2nd_matrix.pdf}
% 	\end{minipage}
% 	}

% 	\subfigure[fifth-order scheme]{
% 	\begin{minipage}[t]{0.46\linewidth}
% 	\centering
% 	\includegraphics[width=0.9\textwidth]{5_1_HLLC_5nd_matrix.pdf}
% 	\end{minipage}
% 	}
% 	\subfigure[density contour of fifth-order scheme]{
% 	\begin{minipage}[t]{0.46\linewidth}
% 	\centering
% 	\includegraphics[width=0.9\textwidth]{5_1_HLLC_5nd_contour.pdf}
% 	\end{minipage}
% 	}

% 	\centering
% 	\caption{HLLC solver with the wave speeds estimate proposed in \cite{Einfeldt_Godunovtype_1991}.(Grid with 11$ \times $11 cells, $ M_0=20 $ and $ \varepsilon =0.1 $. (a) - (c) are the distributions of the eigenvalues in the complex plane; (d) is the density contour of the fifth-order scheme.)}\label{fig HLLC matrix}
% \end{figure}

\begin{figure}[htbp]
	\centering

	\subfigure[first-order scheme]{
	\begin{minipage}[t]{0.46\linewidth}
	\centering
	\includegraphics[width=0.9\textwidth]{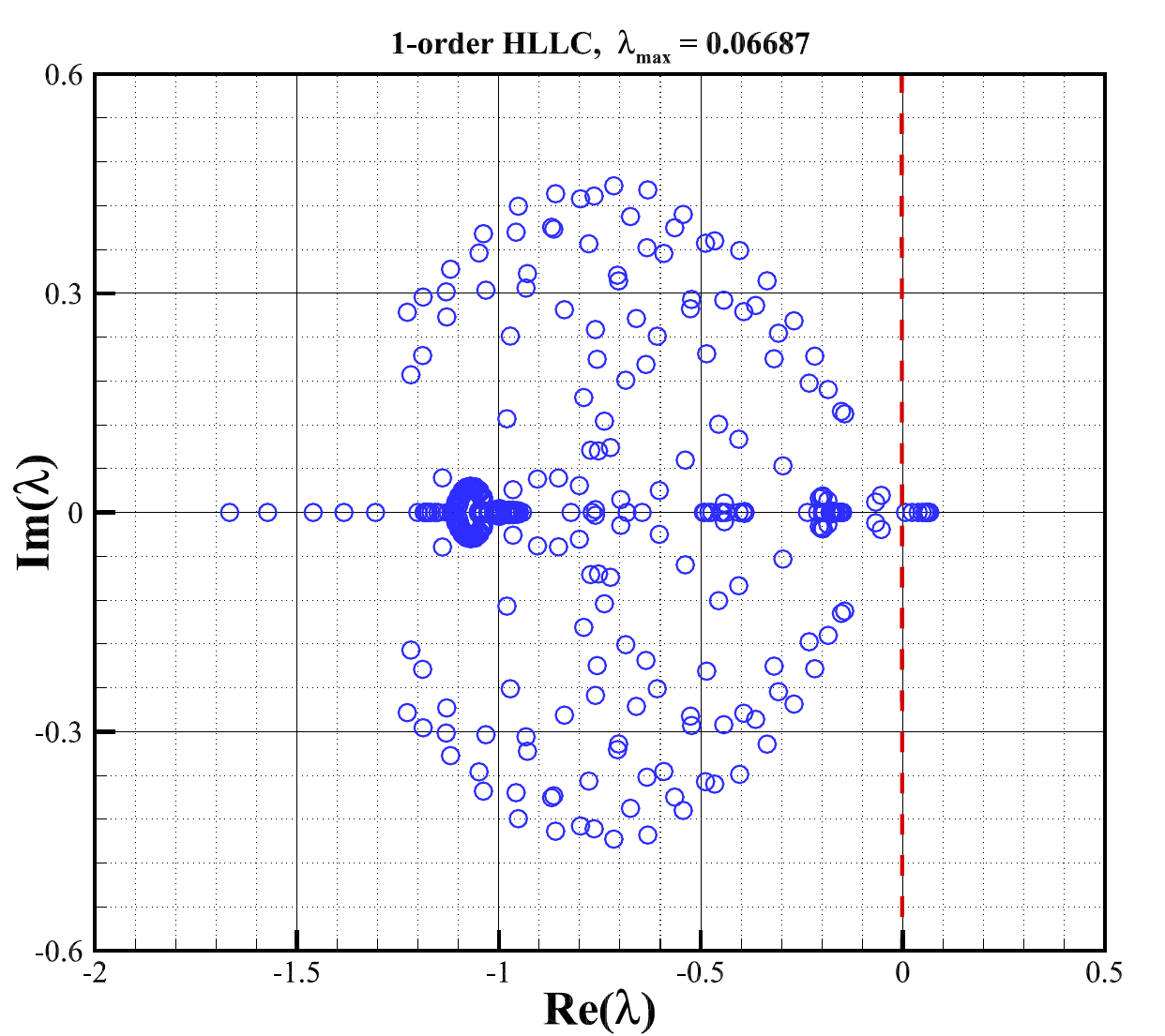}
	\end{minipage}
	}
	\subfigure[second-order scheme]{
	\begin{minipage}[t]{0.46\linewidth}
	\centering
	\includegraphics[width=0.9\textwidth]{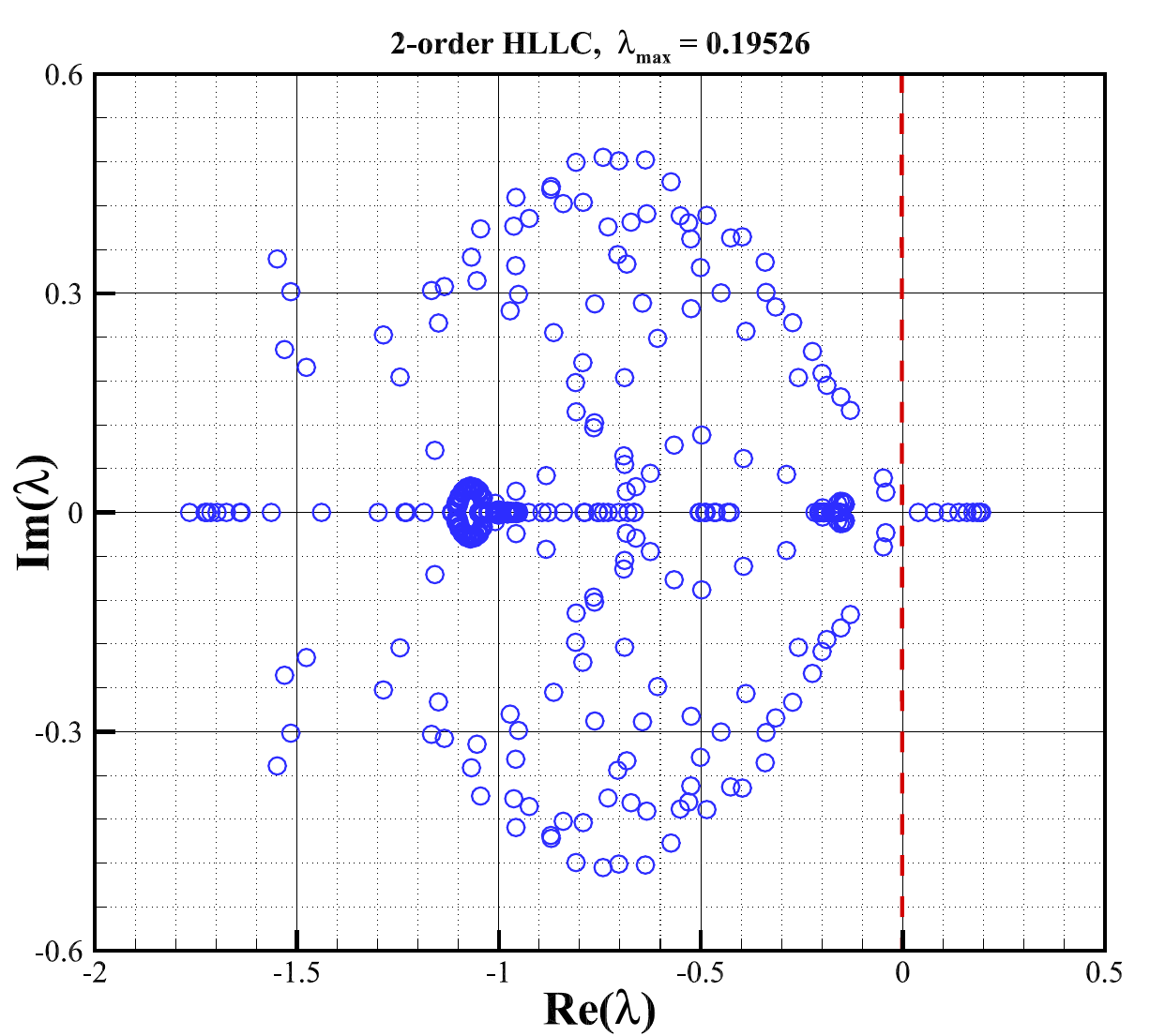}
	\end{minipage}
	}

	\subfigure[fifth-order scheme]{
	\begin{minipage}[t]{0.46\linewidth}
	\centering
	\includegraphics[width=0.9\textwidth]{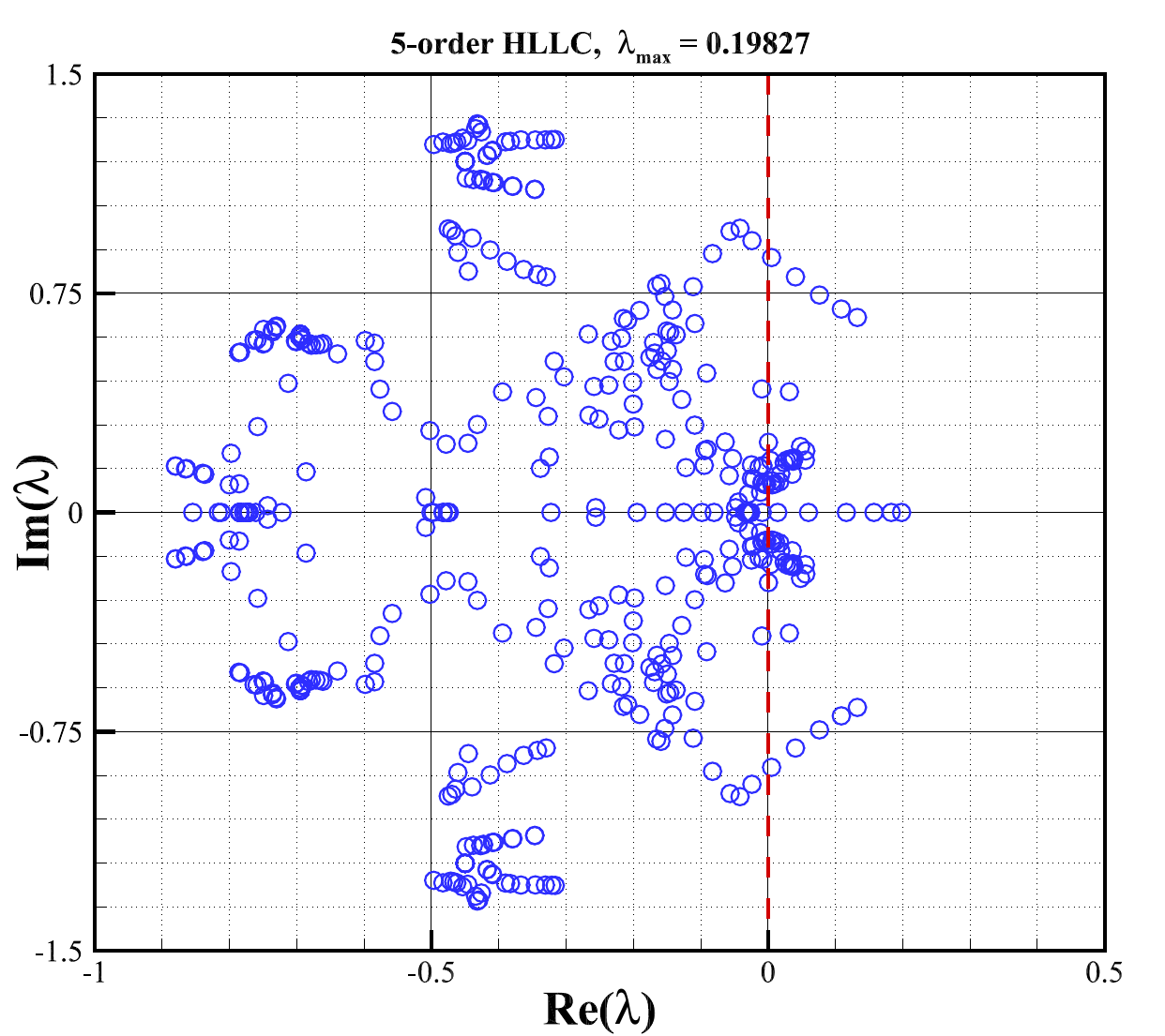}
	\end{minipage}
	}
	\subfigure[density contour of fifth-order scheme]{
	\begin{minipage}[t]{0.46\linewidth}
	\centering
	\includegraphics[width=0.9\textwidth]{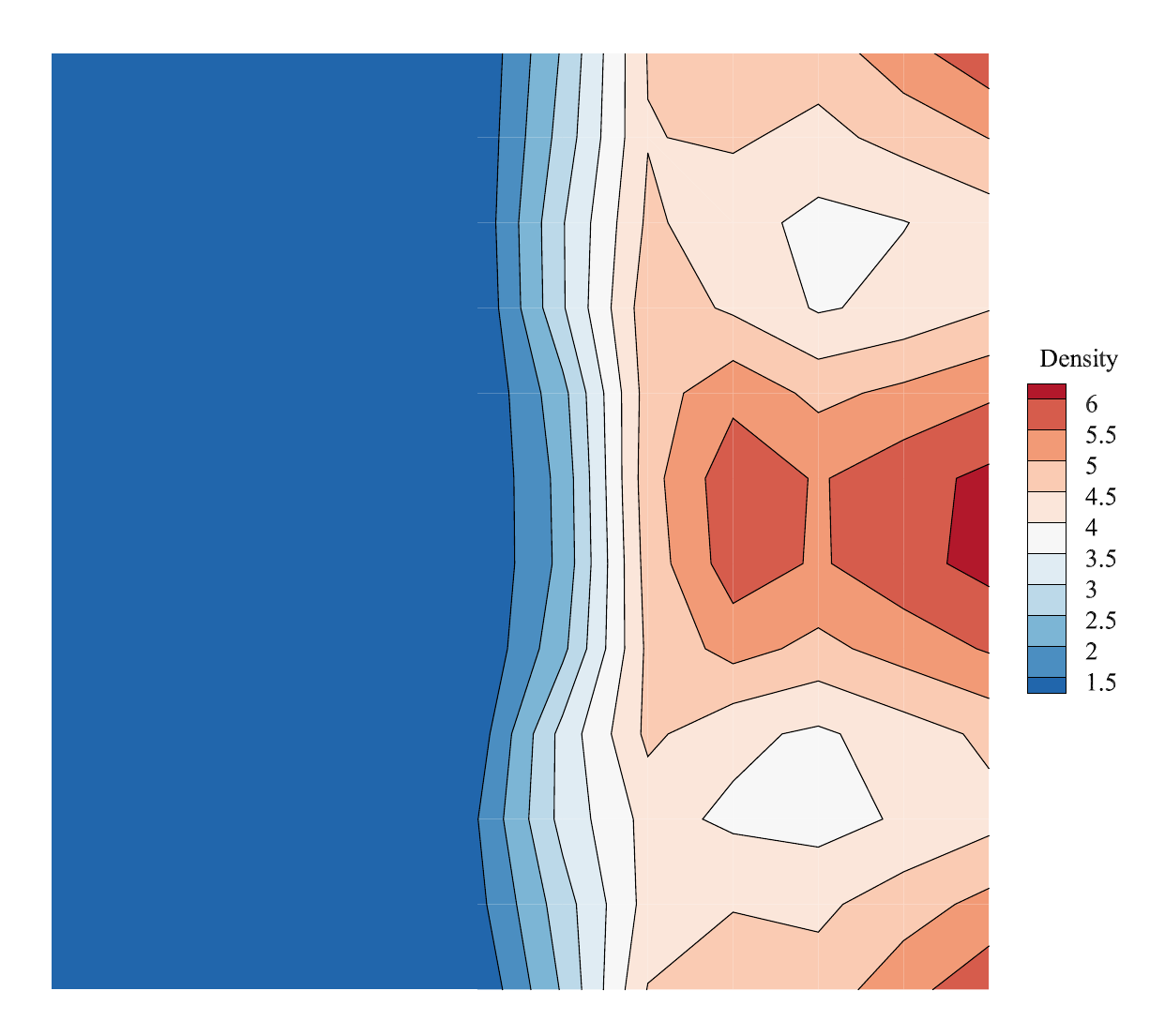}
	\end{minipage}
	}

	\centering
	\caption{HLLC solver.(Grid with 11$ \times $11 cells, $ M_0=20 $ and $ \varepsilon =0.1 $. (a) - (c) are the distributions of the eigenvalues in the complex plane; (d) is the density contour of the fifth-order scheme.)}\label{fig dissipative HLLC matrix}
\end{figure}

\begin{figure}[htbp]
	\centering

	\subfigure[first-order scheme]{
	\begin{minipage}[t]{0.46\linewidth}
	\centering
	\includegraphics[width=0.9\textwidth]{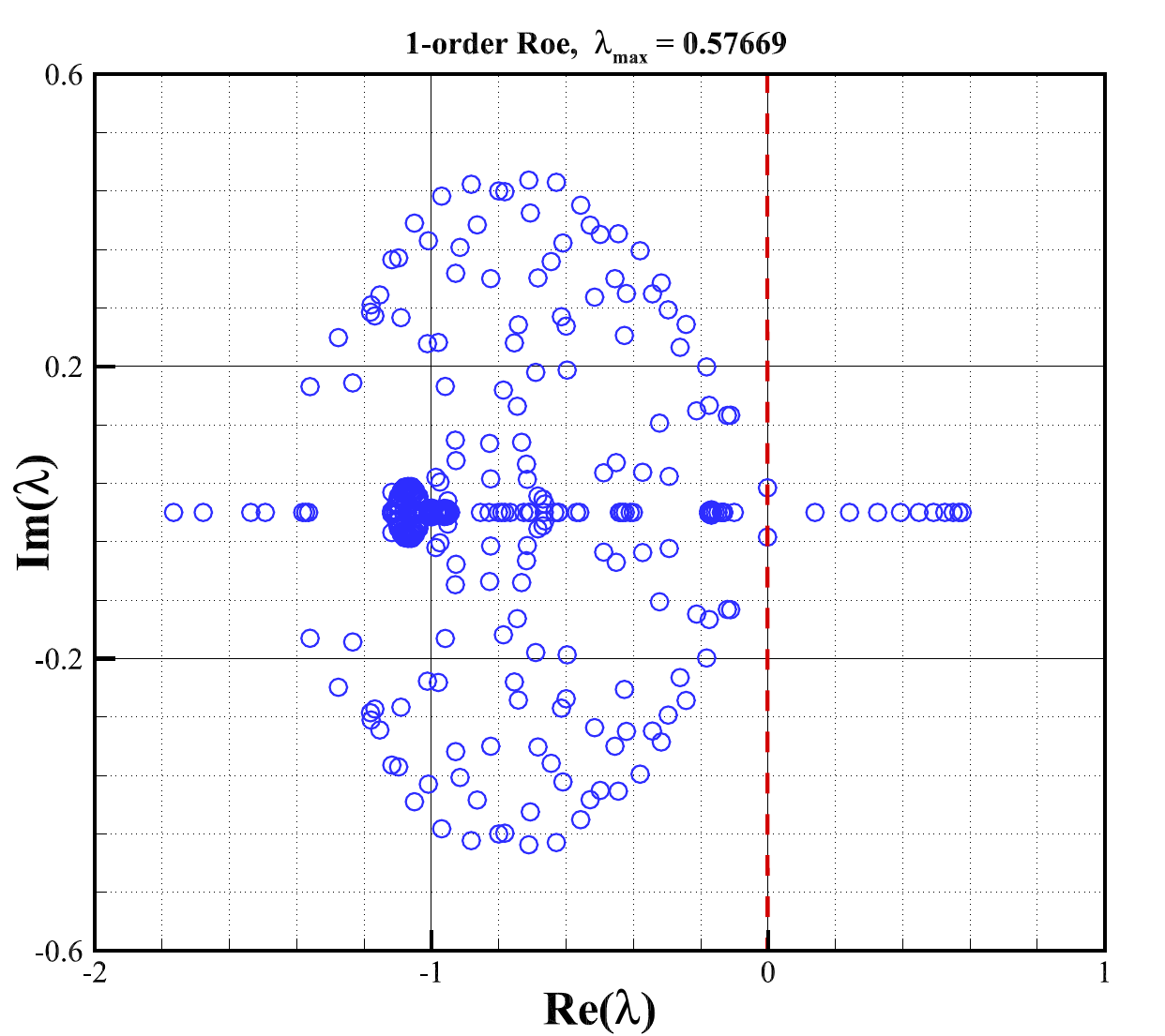}
	\end{minipage}
	}
	\subfigure[second-order scheme]{
	\begin{minipage}[t]{0.46\linewidth}
	\centering
	\includegraphics[width=0.9\textwidth]{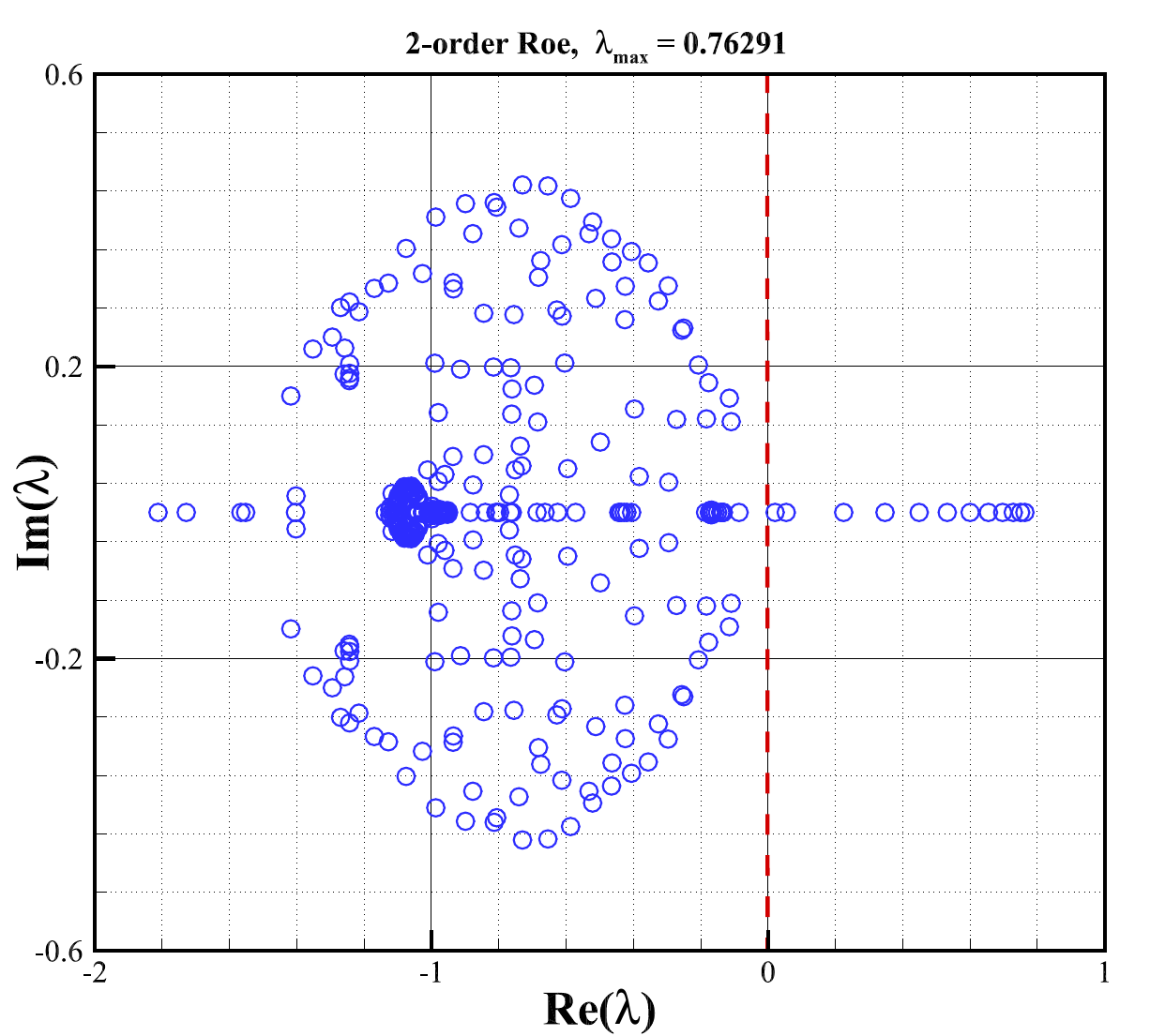}
	\end{minipage}
	}

	\subfigure[fifth-order scheme]{
	\begin{minipage}[t]{0.46\linewidth}
	\centering
	\includegraphics[width=0.9\textwidth]{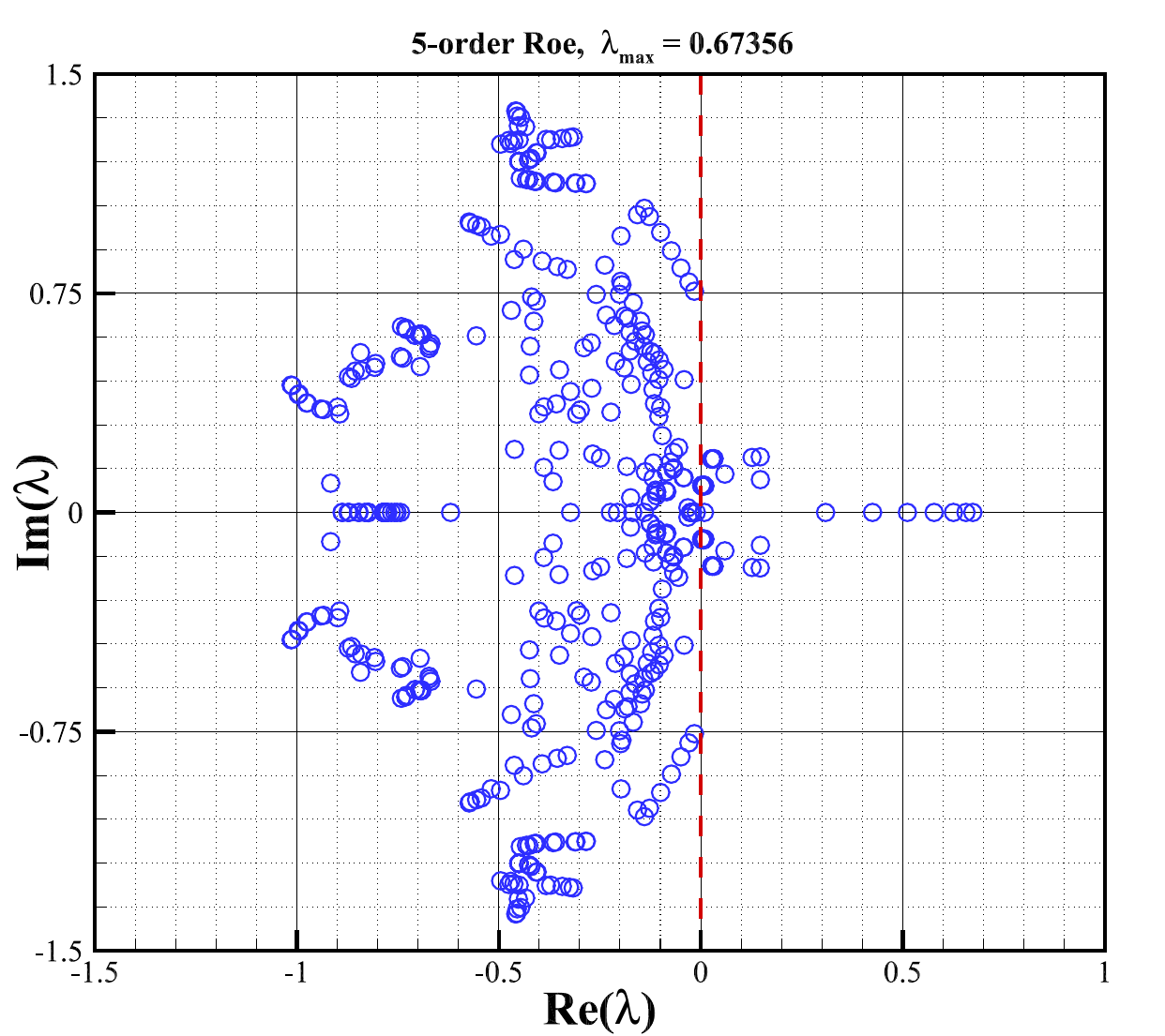}
	\end{minipage}
	}
	\subfigure[density contour of fifth-order scheme]{
	\begin{minipage}[t]{0.46\linewidth}
	\centering
	\includegraphics[width=0.9\textwidth]{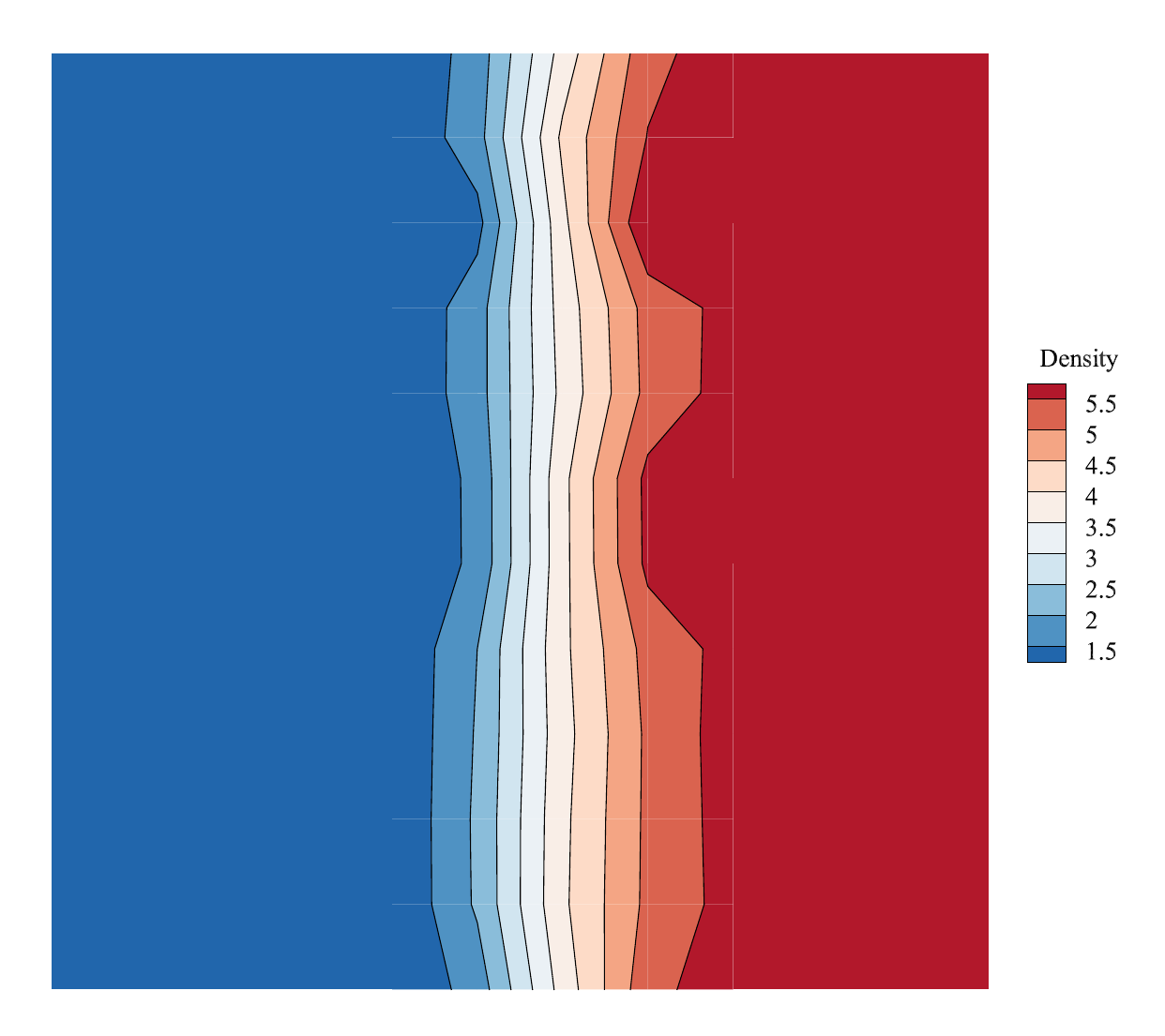}
	\end{minipage}
	}

	\centering
	\caption{Roe solver.(Grid with 11$ \times $11 cells, $ M_0=20 $ and $ \varepsilon =0.1 $. (a) - (c) are the distributions of the eigenvalues in the complex plane; (d) is the density contour of the fifth-order scheme.)}\label{fig Roe matrix}
\end{figure}

\begin{figure}[htbp]
	\centering

	\subfigure[first-order scheme]{
	\begin{minipage}[t]{0.46\linewidth}
	\centering
	\includegraphics[width=0.9\textwidth]{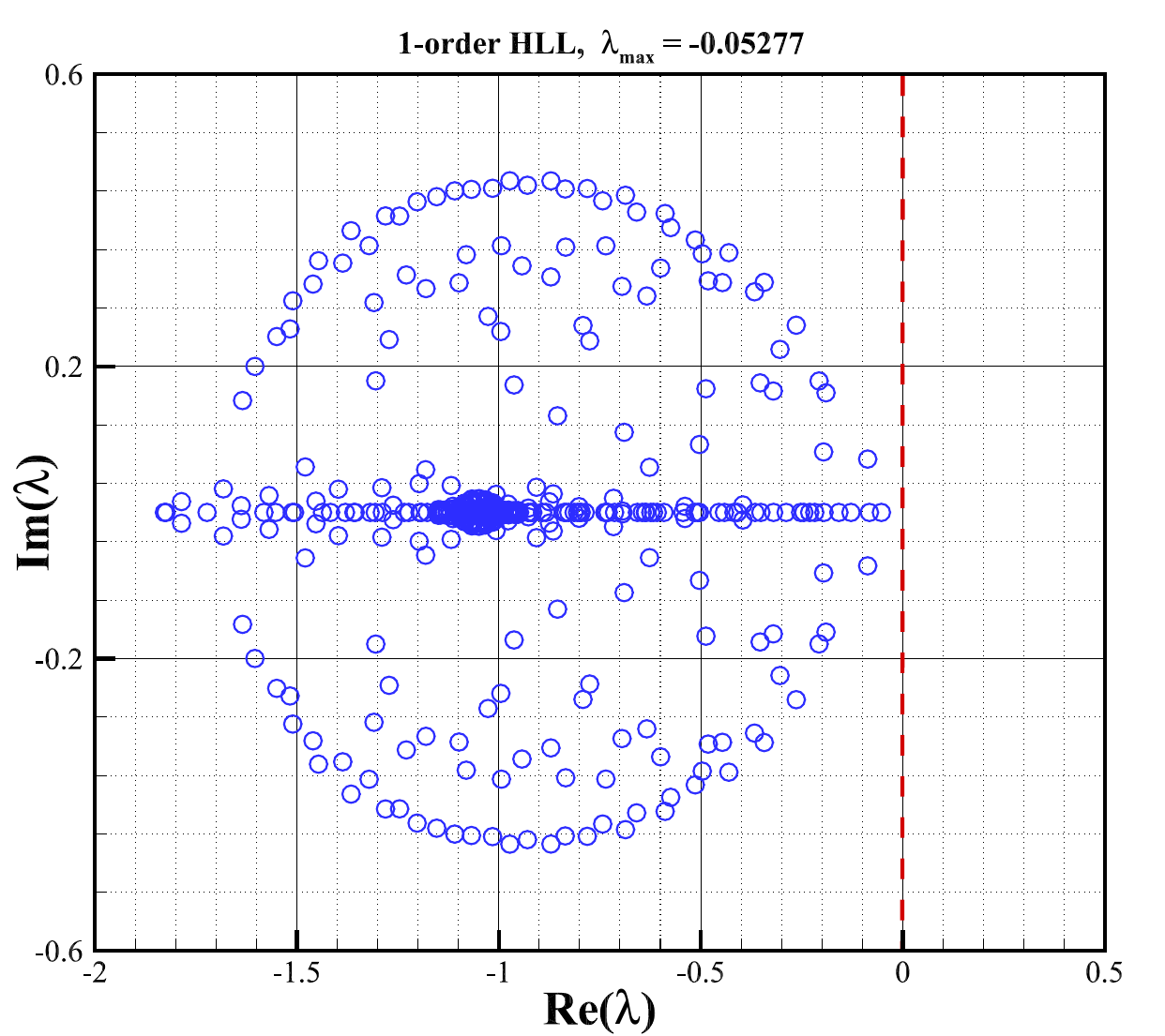}
	\end{minipage}
	}
	\subfigure[second-order scheme]{
	\begin{minipage}[t]{0.46\linewidth}
	\centering
	\includegraphics[width=0.9\textwidth]{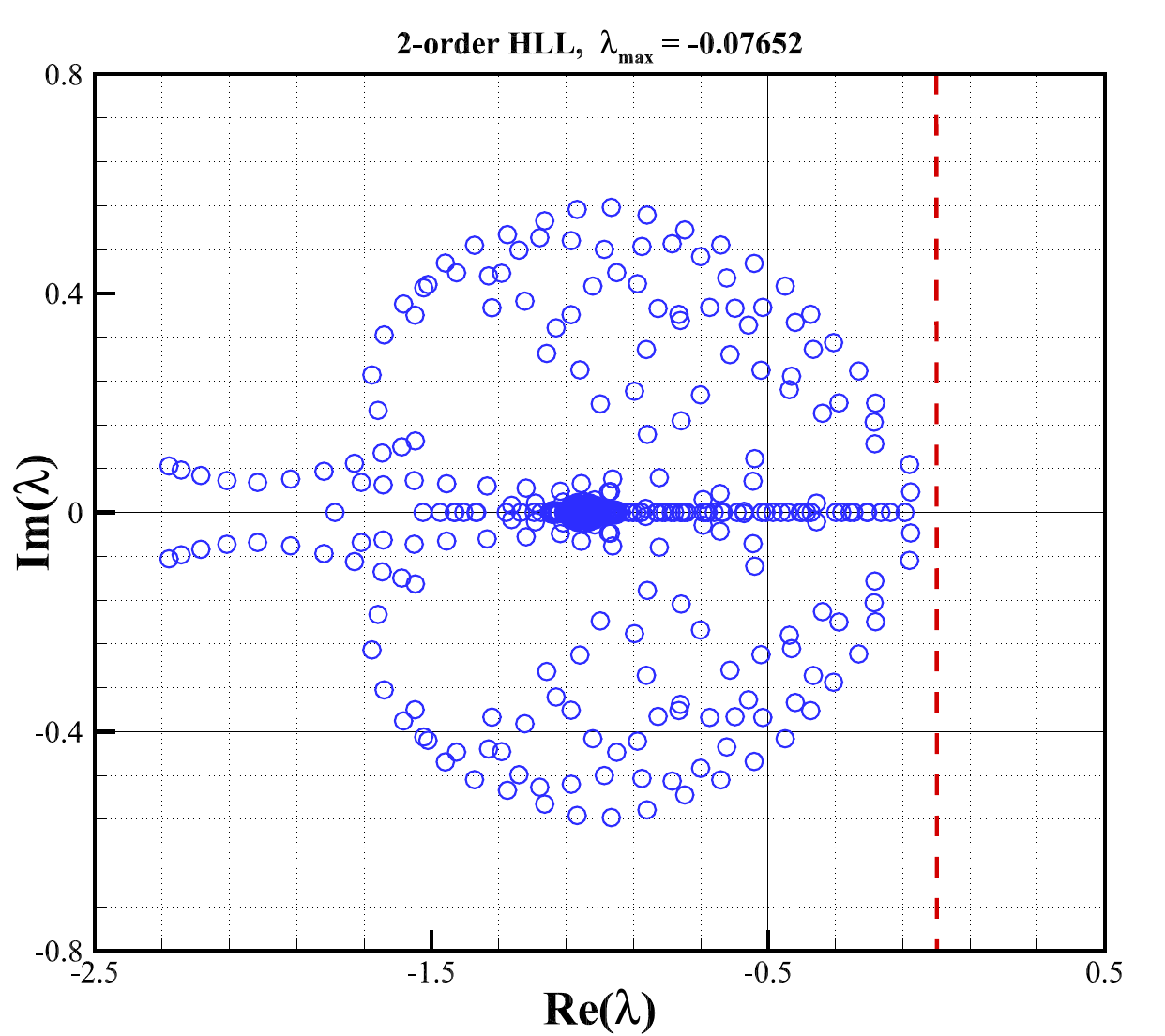}
	\end{minipage}
	}

	\subfigure[fifth-order scheme]{
	\begin{minipage}[t]{0.46\linewidth}
	\centering
	\includegraphics[width=0.9\textwidth]{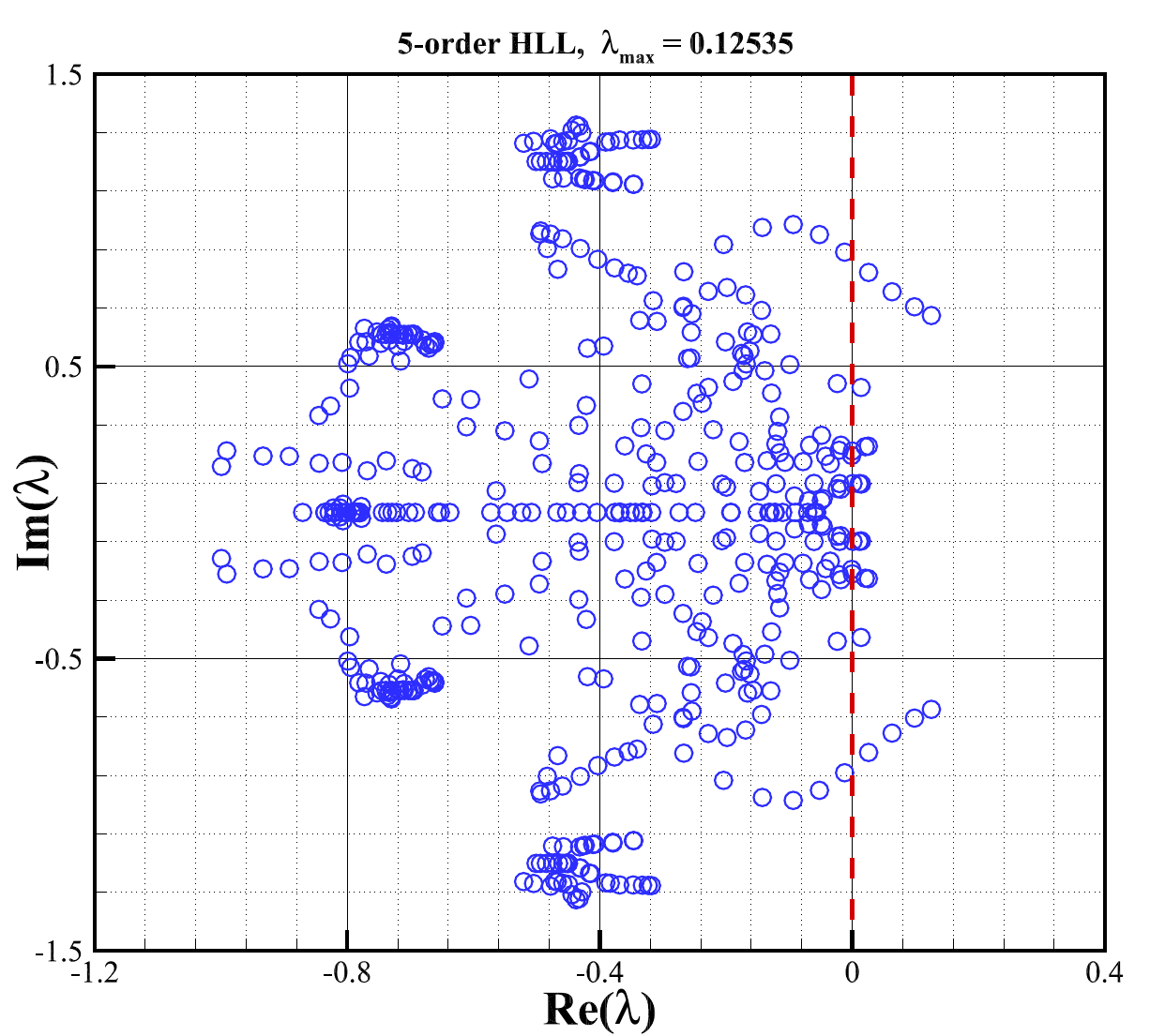}
	\end{minipage}
	}
	\subfigure[density contour of fifth-order scheme]{
	\begin{minipage}[t]{0.46\linewidth}
	\centering
	\includegraphics[width=0.9\textwidth]{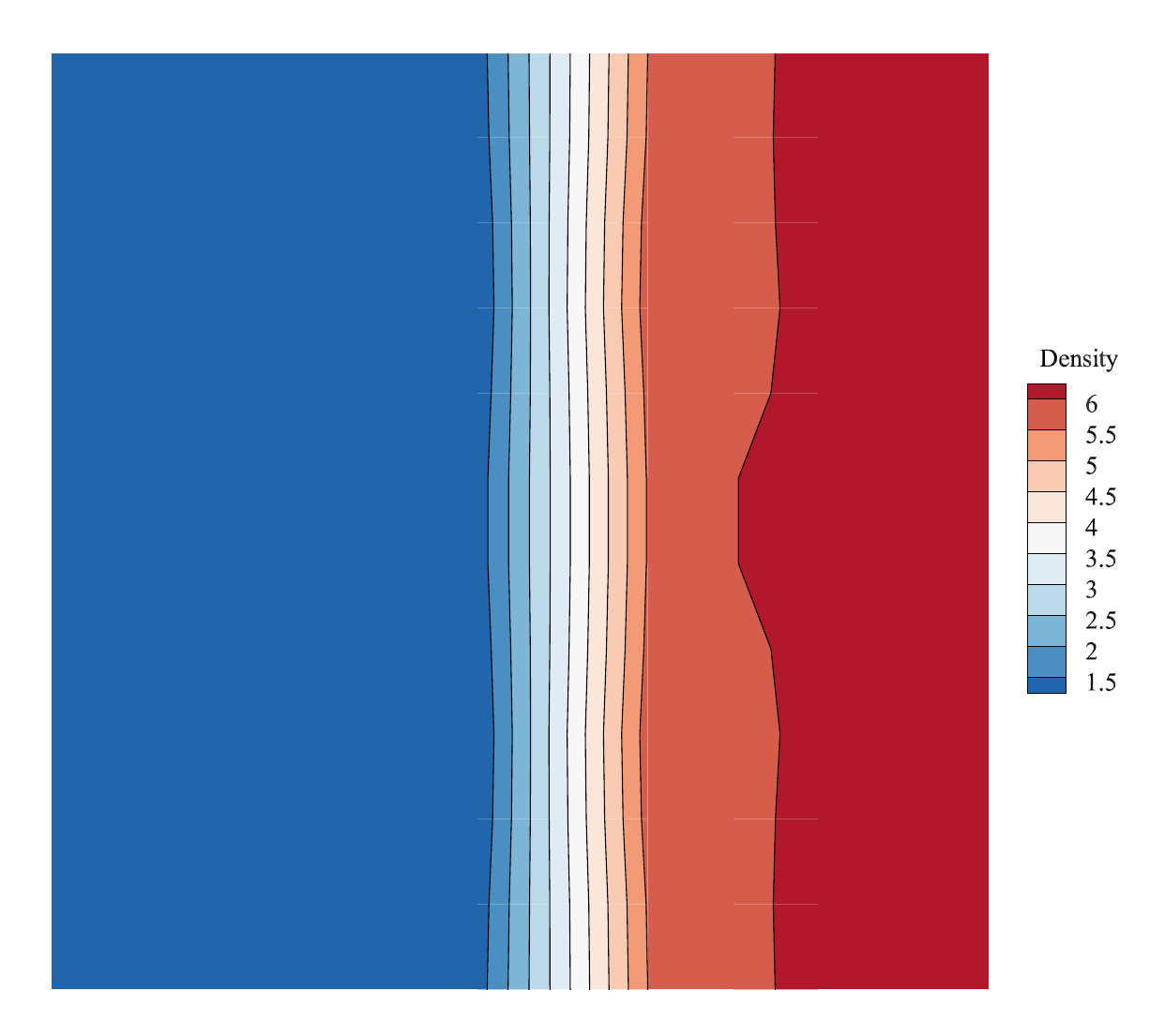}
	\end{minipage}
	}

	\centering
	\caption{HLL solver.(Grid with 11$ \times $11 cells, $ M_0=20 $ and $ \varepsilon =0.1 $. (a) - (c) are the distributions of the eigenvalues in the complex plane; (d) is the density contour of the fifth-order scheme.)}\label{fig HLL matrix}
\end{figure}

\begin{figure}[htbp]
	\centering

	\subfigure[first-order scheme]{
	\begin{minipage}[t]{0.46\linewidth}
	\centering
	\includegraphics[width=0.9\textwidth]{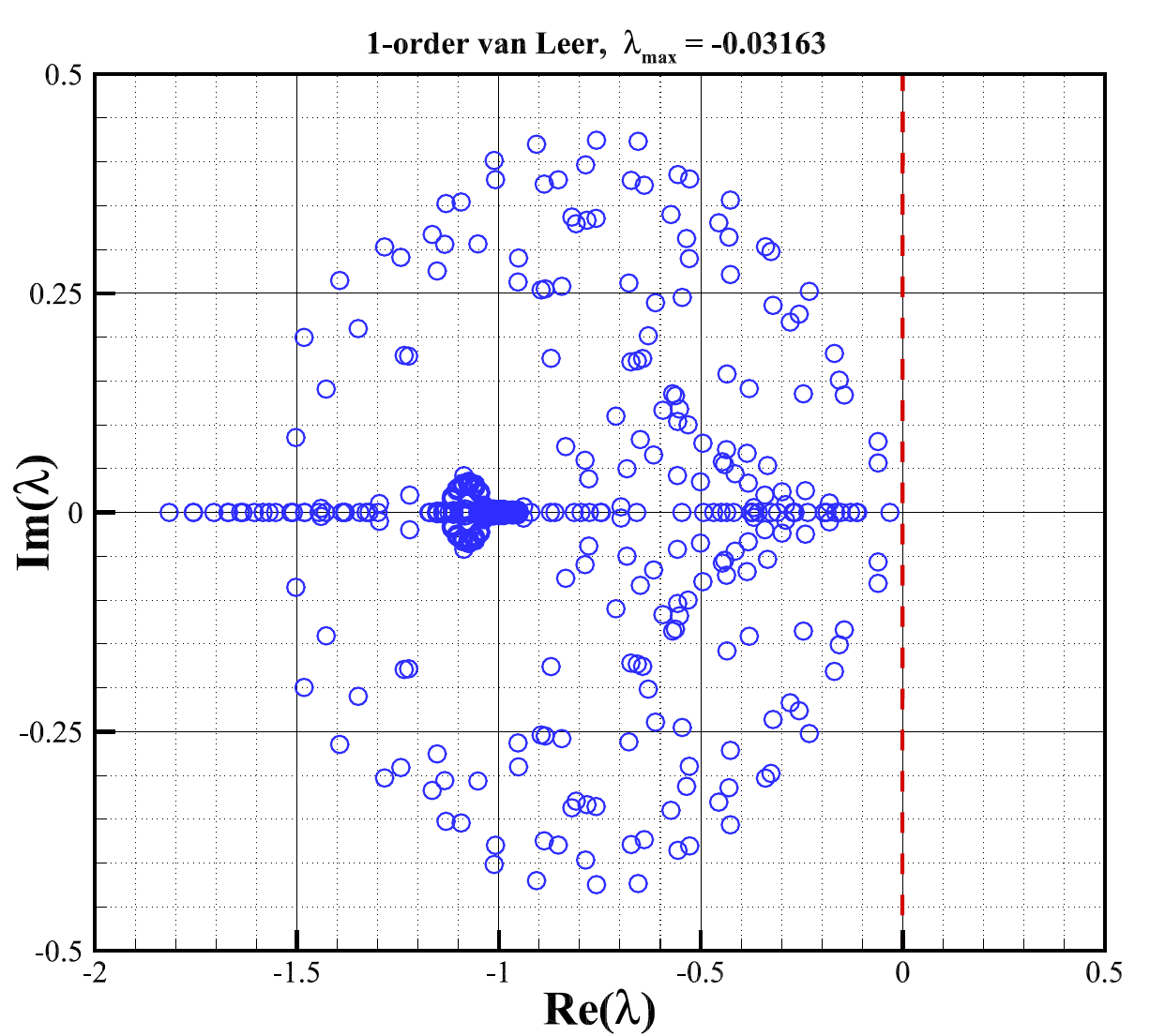}
	\end{minipage}
	}
	\subfigure[second-order scheme]{
	\begin{minipage}[t]{0.46\linewidth}
	\centering
	\includegraphics[width=0.9\textwidth]{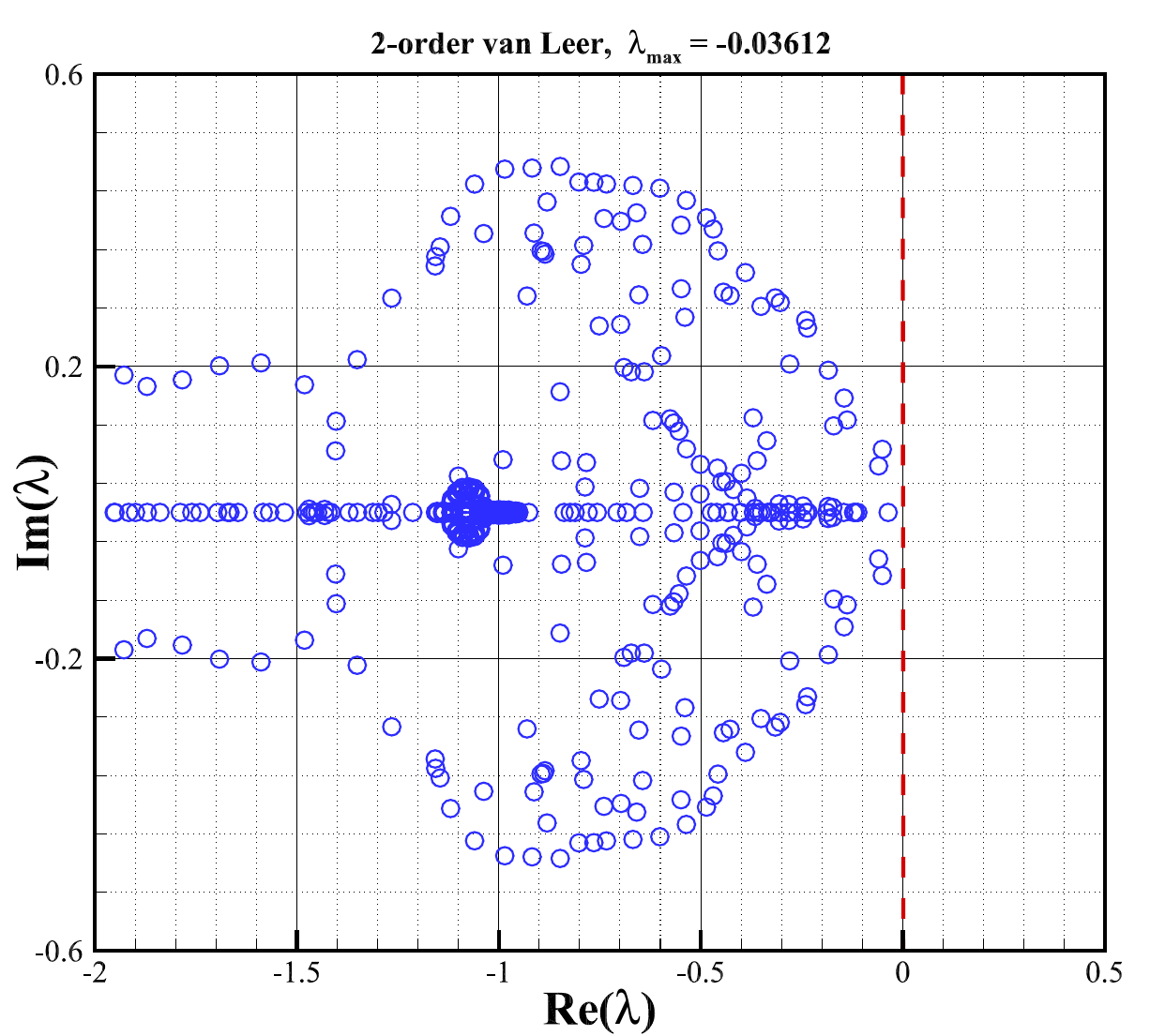}
	\end{minipage}
	}

	\subfigure[fifth-order scheme]{
	\begin{minipage}[t]{0.46\linewidth}
	\centering
	\includegraphics[width=0.9\textwidth]{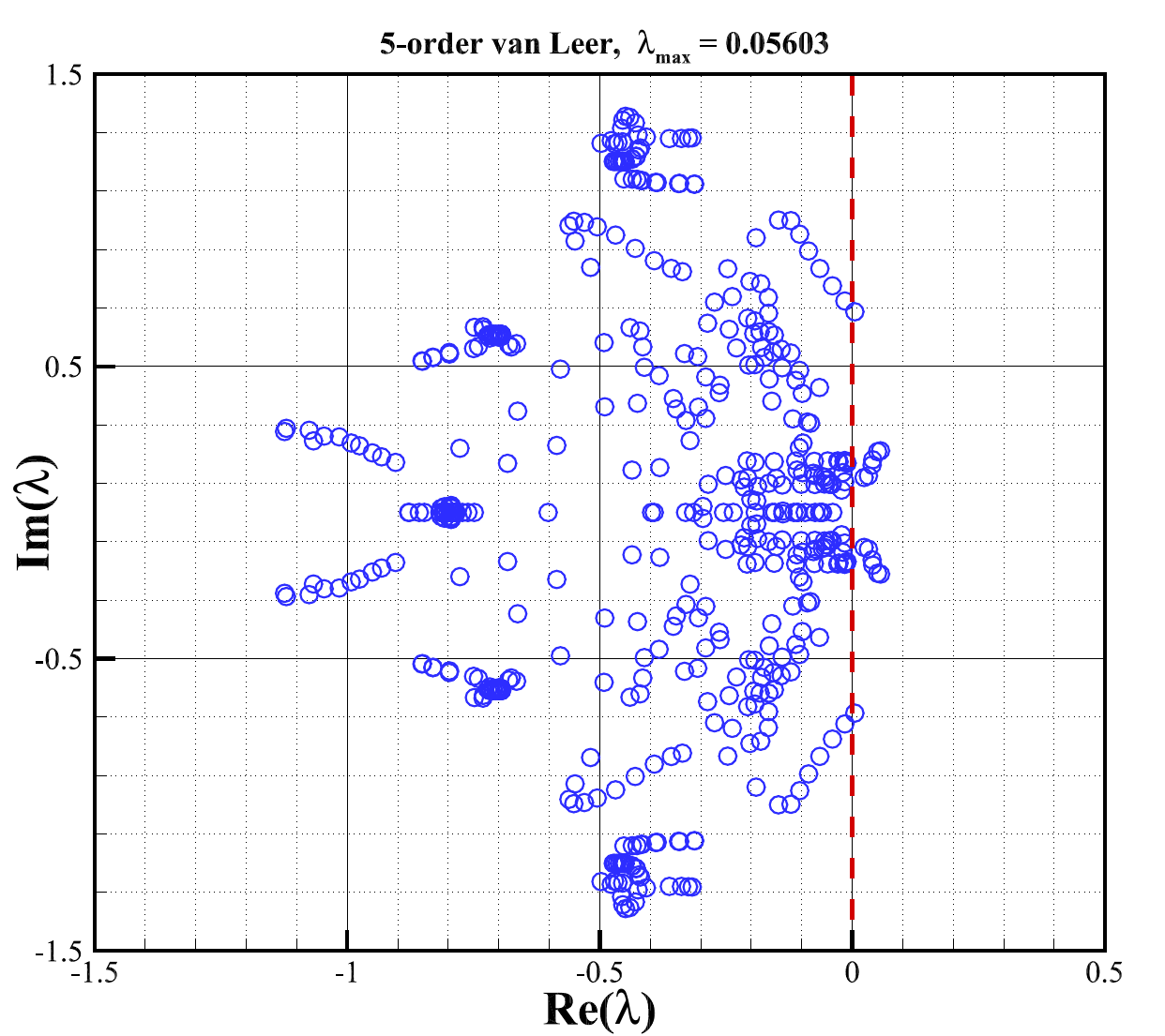}
	\end{minipage}
	}
	\subfigure[density contour of fifth-order scheme]{
	\begin{minipage}[t]{0.46\linewidth}
	\centering
	\includegraphics[width=0.9\textwidth]{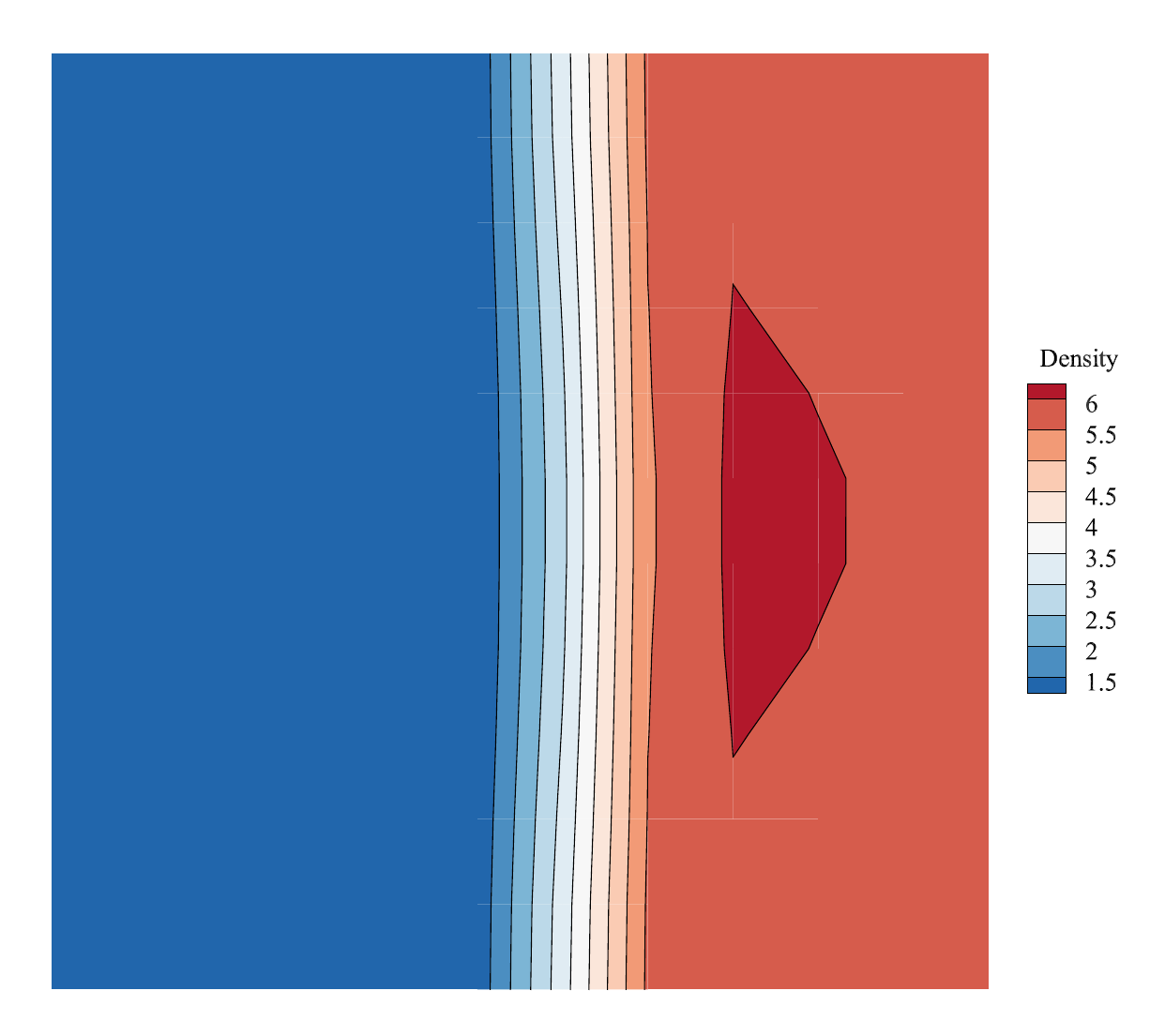}
	\end{minipage}
	}

	\centering
	\caption{van Leer solver.(Grid with 11$ \times $11 cells, $ M_0=20 $ and $ \varepsilon =0.1 $. (a) - (c) are the distributions of the eigenvalues in the complex plane; (d) is the density contour of the fifth-order scheme.)}\label{fig van Leer matrix}
\end{figure}

The stability of fifth-order schemes is investigated and compared with the first- and second-order schemes in this section. To this end, the solvers of HLL \cite{Harten_Upstream_1983}, HLLC \cite{Toro_Restoration_1994}, Roe \cite{Roe_Approximate_1981}, and van Leer \cite{VanLeer_Fluxvector_1982} are employed. Note that in this paper, the HLL and HLLC solvers employ the wave speeds estimate proposed by Davis \cite{Davis_Simplified_1988}. Conventionally, it has been demonstrated that the first- and second-order schemes equipped with the HLL and van Leer solvers can capture the shocks stably, while will suffer from the shock instability problem if the Roe and HLLC solvers are employed. The matrix stability analysis of these solvers and their corresponding density contours of fifth-order schemes are shown in Fig.\ref{fig dissipative HLLC matrix} - Fig.\ref{fig van Leer matrix}. The matrix stability analysis method for the first- and second-order schemes can be found in \cite{Dumbser_Matrix_2004,Ren_Numerical_2023}. As shown in Fig.\ref{fig dissipative HLLC matrix} - Fig.\ref{fig van Leer matrix}, the results of the first- and second-order MUSCL schemes (van Albada limiter is used) are consistent with previous understanding. However, for the fifth-order schemes, the situations are different.

\subsubsection{WENO schemes with low-dissipative solvers}
Firstly, it can be observed from Fig.\ref{fig dissipative HLLC matrix} (c) and \ref{fig Roe matrix} (c) that there are eigenvalues with positive real parts for the fifth-order schemes equipped with the HLLC and Roe solvers, indicating there instability. So do the first- and second-order cases. The corresponding density contours validate the analysis results. It is worth noting that for both solvers, the value of $ \max(\text{Re}(\lambda)) $ for fifth-order schemes is smaller than that of second-order schemes but greater than that of the first-order cases. This suggests that when using fifth-order schemes, the computation will develop towards instability more quickly than with first-order schemes but more slowly than with second-order schemes. This trend is further confirmed by the numerical experiment shown in Fig.\ref{fig evolution of first-, second- and fifth-order schemes}. This result is not completely consistent with our previous understanding. At first sight, it is believed that higher accuracy leads to less dissipation, consequently leading to faster instability growth. We speculate that this can attribute to the weighted combination of perturbations with the WENO scheme for reconstruction, which results in the ``smearing" of large perturbations, thus slowing down the evolution of the perturbation error. This is supported by the results in \cite{Zangeneh_Stability_2019}, which find that a larger stencil size is helpful to stabilize the instability. However, the ``smearing'' process can only slow down the evolution of perturbation errors, but can not determine that the computation is stable or unstable. In fact, as shown in Fig.\ref{fig evolution of first-, second- and fifth-order schemes}, although the perturbation error develop more slowly, when it is large (in the nonlinear stage), the fifth-order scheme gets NAN and collapses at $ t=47.8 $. This indicates that fifth-order schemes face more severe stability problems when the perturbation error reaches a certain magnitude.

\begin{figure}[htbp]
	\centering
	\includegraphics[width=0.6\textwidth]{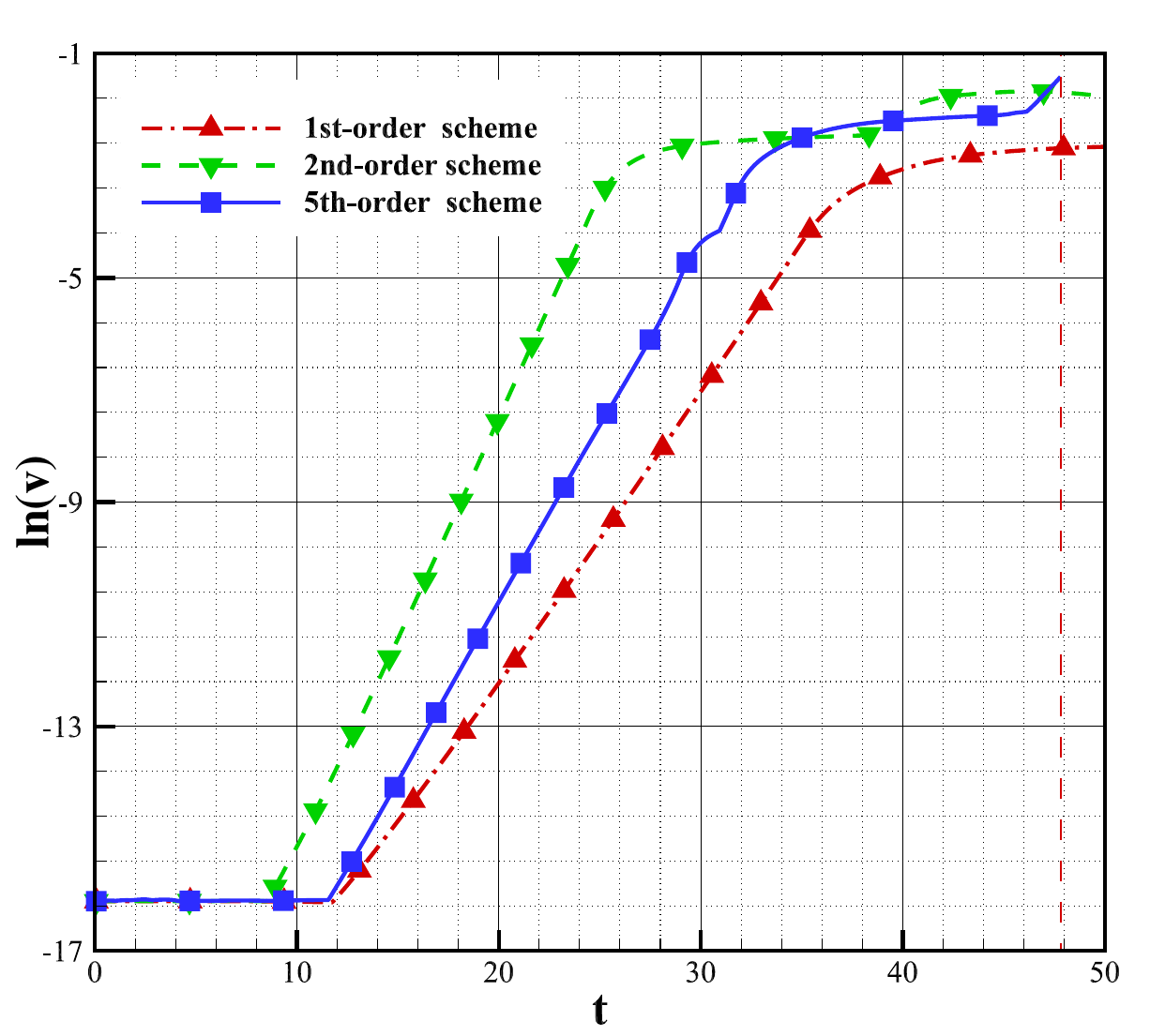}
	\caption{Evolution of the perturbation error calculated by first-, second-, and fifth-order schemes.(Grid with 11$ \times $11 cells, $ M_0=20 $ and $ \varepsilon=0.1 $.)}
	\label{fig evolution of first-, second- and fifth-order schemes}
\end{figure}

\subsubsection{WENO schemes with dissipative solvers}
\begin{figure}[htbp]
	\centering
	\includegraphics[width=0.6\textwidth]{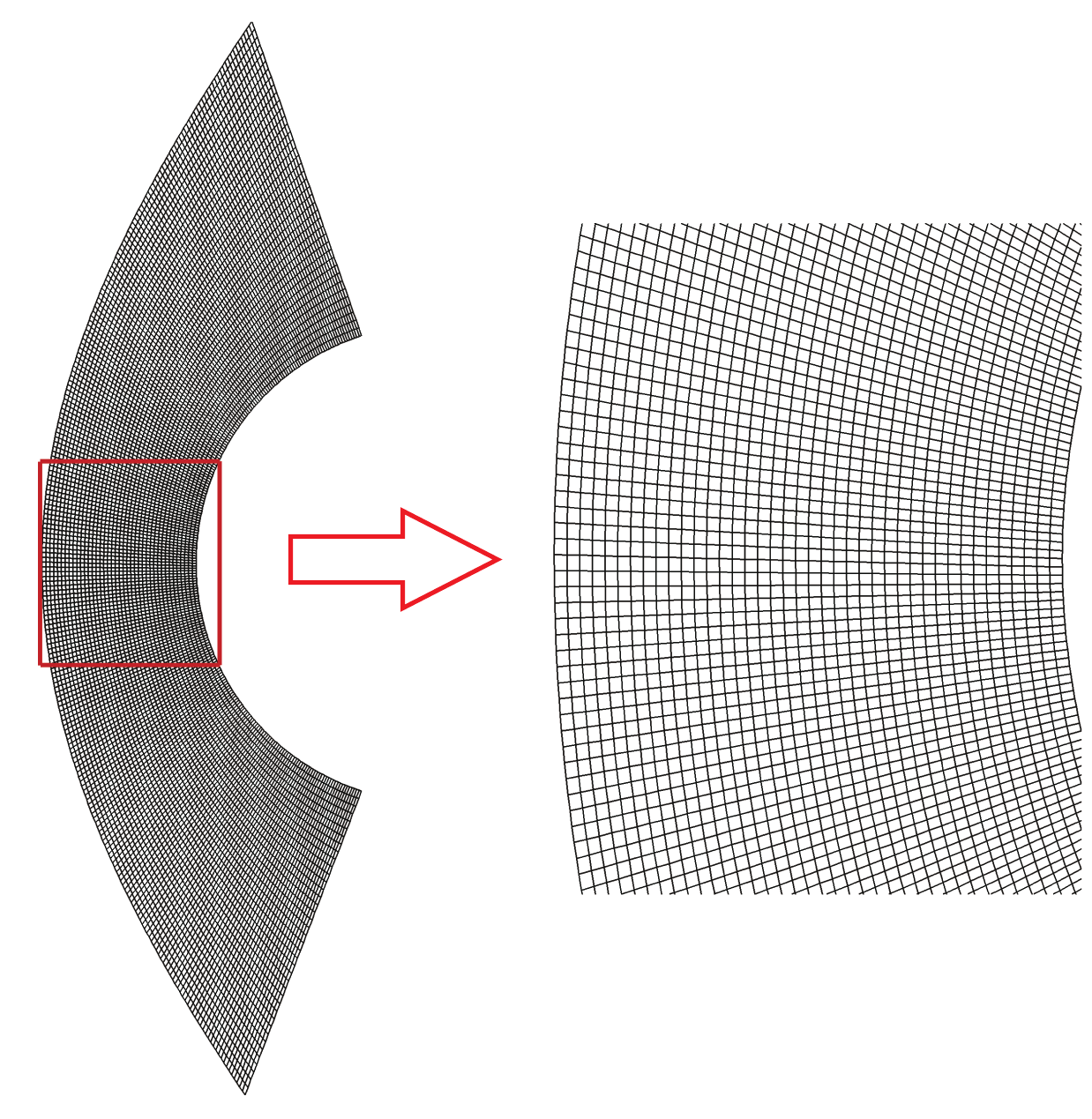}
	\caption{The computational grid of the blunt-body.}
	\label{fig cylinder grid}
\end{figure}

\begin{figure}[htbp]
	\centering

	\subfigure[HLL]{
	\begin{minipage}[t]{0.3\linewidth}
	\centering
	\includegraphics[width=0.9\textwidth]{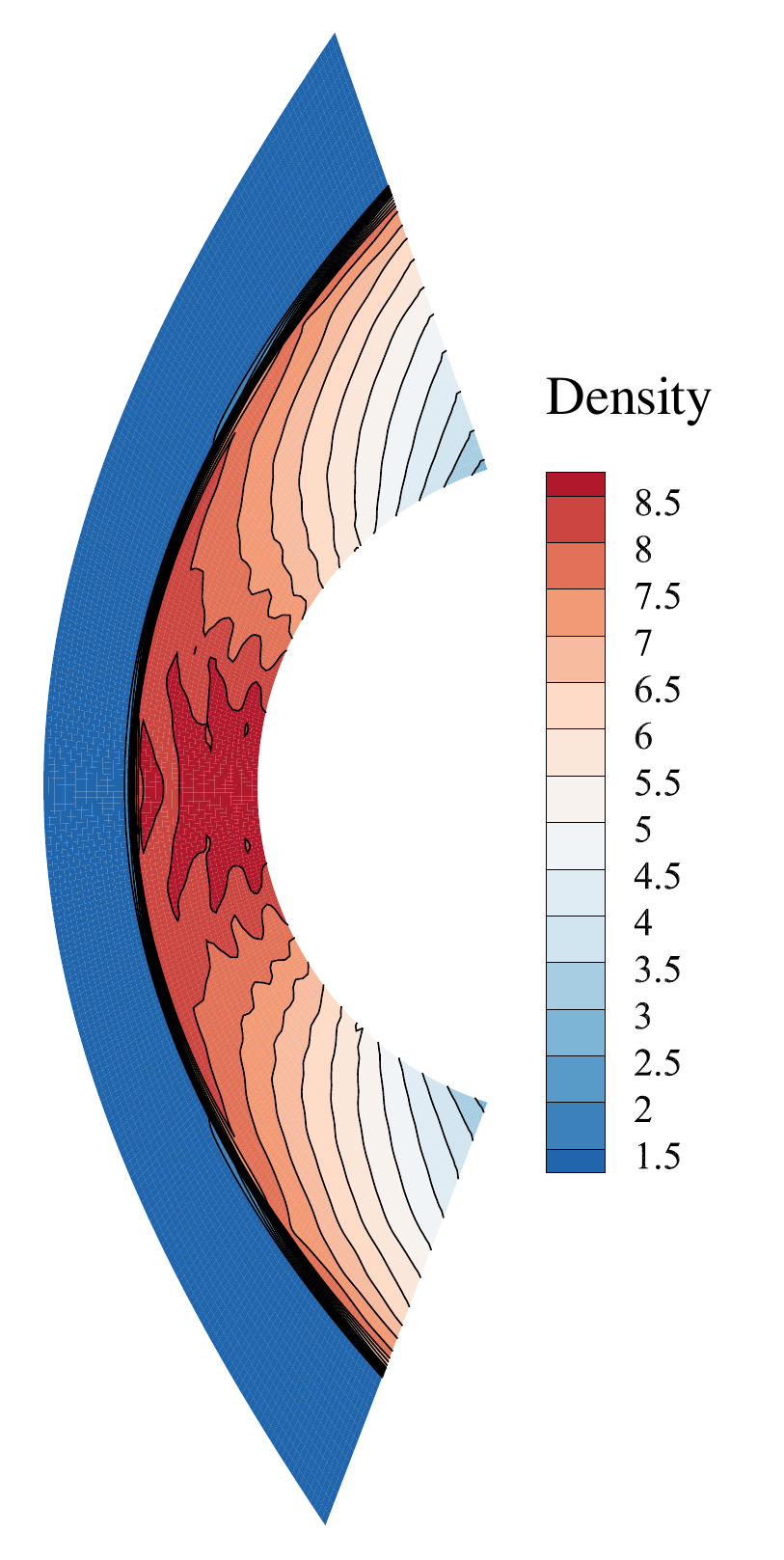}
	\end{minipage}
	}
    \subfigure[van Leer]{
	\begin{minipage}[t]{0.3\linewidth}
	\centering
	\includegraphics[width=0.9\textwidth]{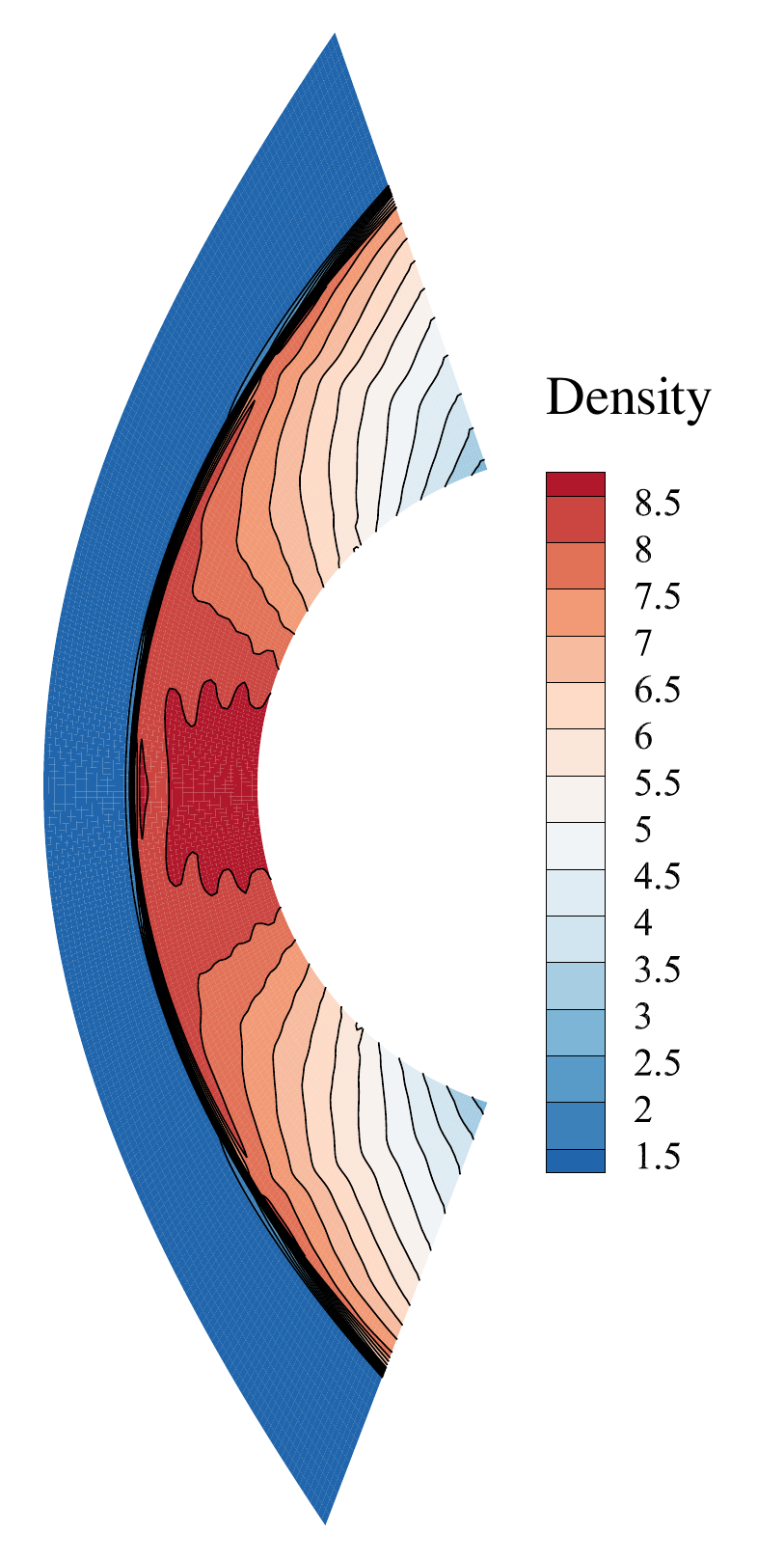}
	\end{minipage}
	}

	\centering
	\caption{The density contours of hypersonic flow over a cylinder computed by fifth-order schemes with different solvers.}\label{fig cylinder}
\end{figure}

It can be observed from Fig.\ref{fig HLL matrix} (c) and \ref{fig van Leer matrix} (c) that, the fifth-order schemes with HLL and van Leer solvers still have eigenvalues with positive real parts indicating their instability. Note that such a result is inconsistent with our previous understanding that when these dissipative solvers are employed, the numerical schemes can capture strong shocks stably. The conclusion is further supported by the corresponding density contours shown in Fig.\ref{fig HLL matrix} (d) and Fig.\ref{fig van Leer matrix} (d), where the shock profiles begin to distort. In order to further validate the analysis results, another typical problem, the hypersonic flow over a blunt-body, is employed. The computational conditions for this problem are as follows
\begin{equation}
	M_0=20.0 \quad \rho=1.4  \quad p=1.0.
\end{equation}
As shown in Fig.\ref{fig cylinder grid}, a structured grid with 240 cells in the circumferential direction and 40 cells in the radial direction covers the computational domain. Fig.\ref{fig cylinder} shows the density contours. As shown, the fifth-order schemes with the two dissipative solvers both produce unstable solutions, confirming the results of the stability analysis.

\begin{figure}[htbp]
	\centering

	\subfigure[distribution of the eigenvalues in the complex plane]{
	\begin{minipage}[t]{0.46\linewidth}
	\centering
	\includegraphics[width=0.85\textwidth]{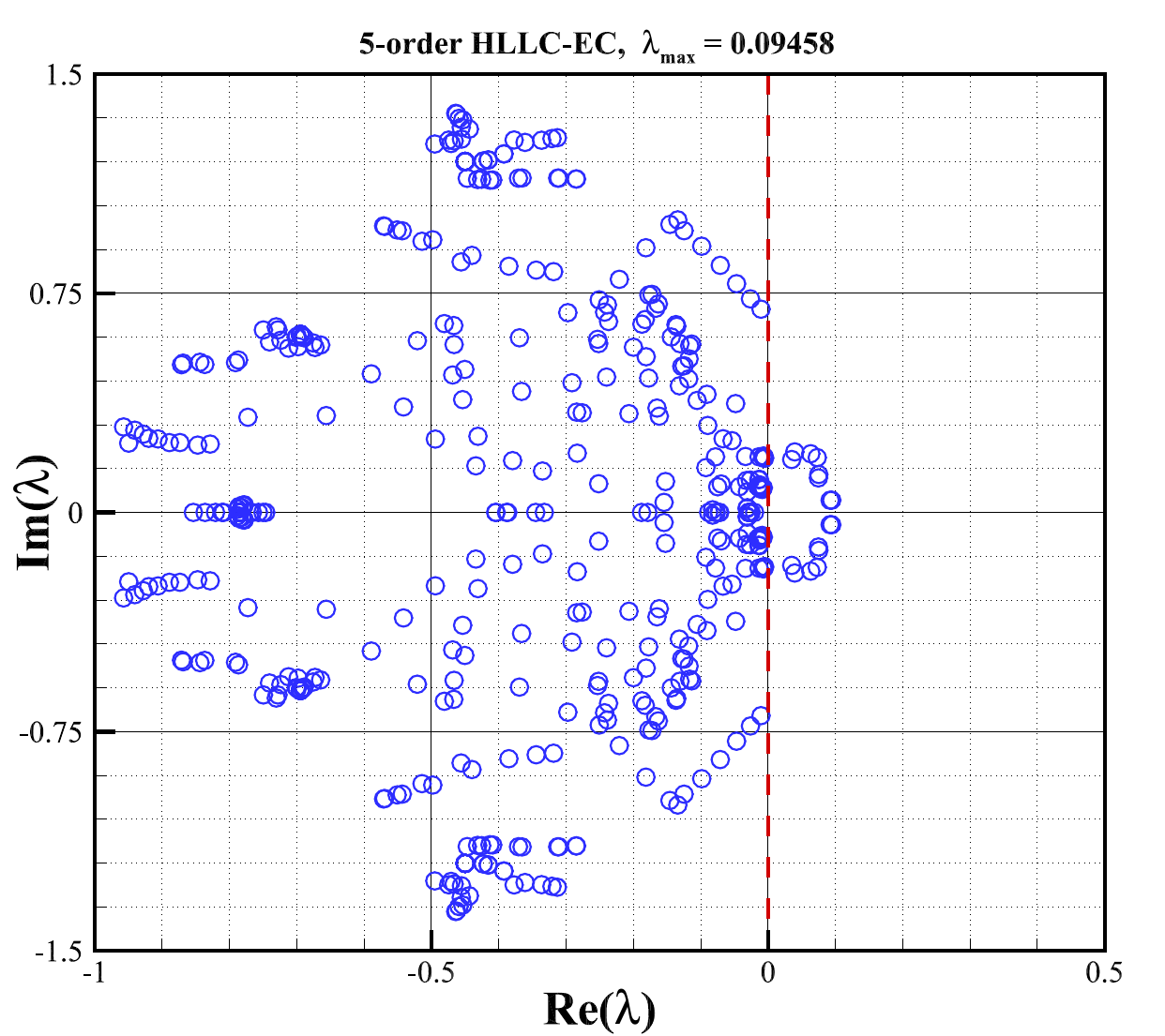}
	\end{minipage}
	}
	\subfigure[density contour]{
	\begin{minipage}[t]{0.46\linewidth}
	\centering
	\includegraphics[width=0.85\textwidth]{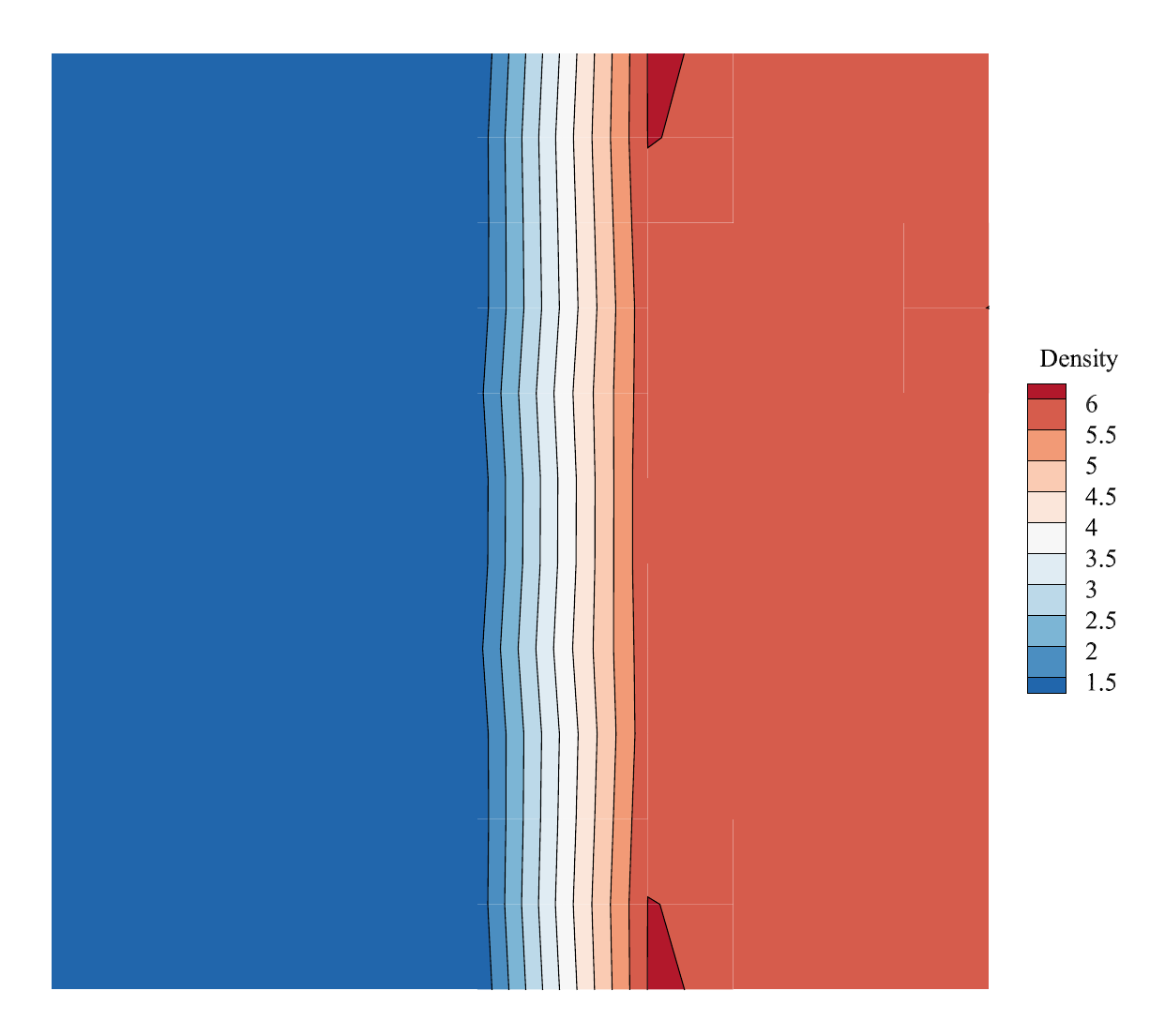}
	\end{minipage}
	}

	\centering
	\caption{Analysis of the fifth-order scheme with HLLC-EC solver.(Grid with 11$ \times $11 cells, $ M_0=20 $ and $ \varepsilon =0.1 $.)}\label{fig HLLC-EC}
\end{figure}

\begin{figure}[htbp]
	\centering

	\subfigure[distribution of the eigenvalues in the complex plane]{
	\begin{minipage}[t]{0.46\linewidth}
	\centering
	\includegraphics[width=0.85\textwidth]{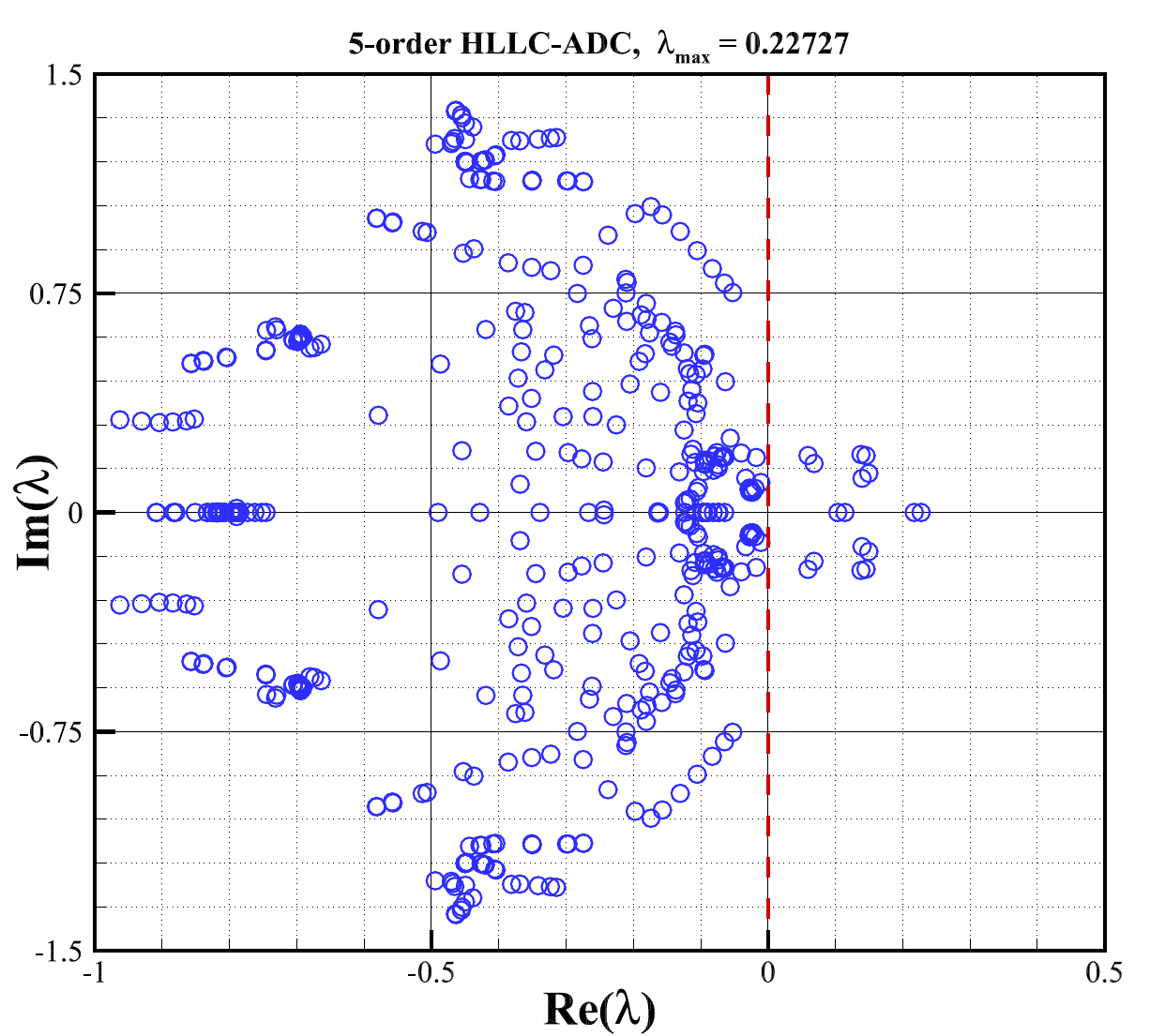}
	\end{minipage}
	}
	\subfigure[density contour]{
	\begin{minipage}[t]{0.46\linewidth}                  
	\centering
	\includegraphics[width=0.85\textwidth]{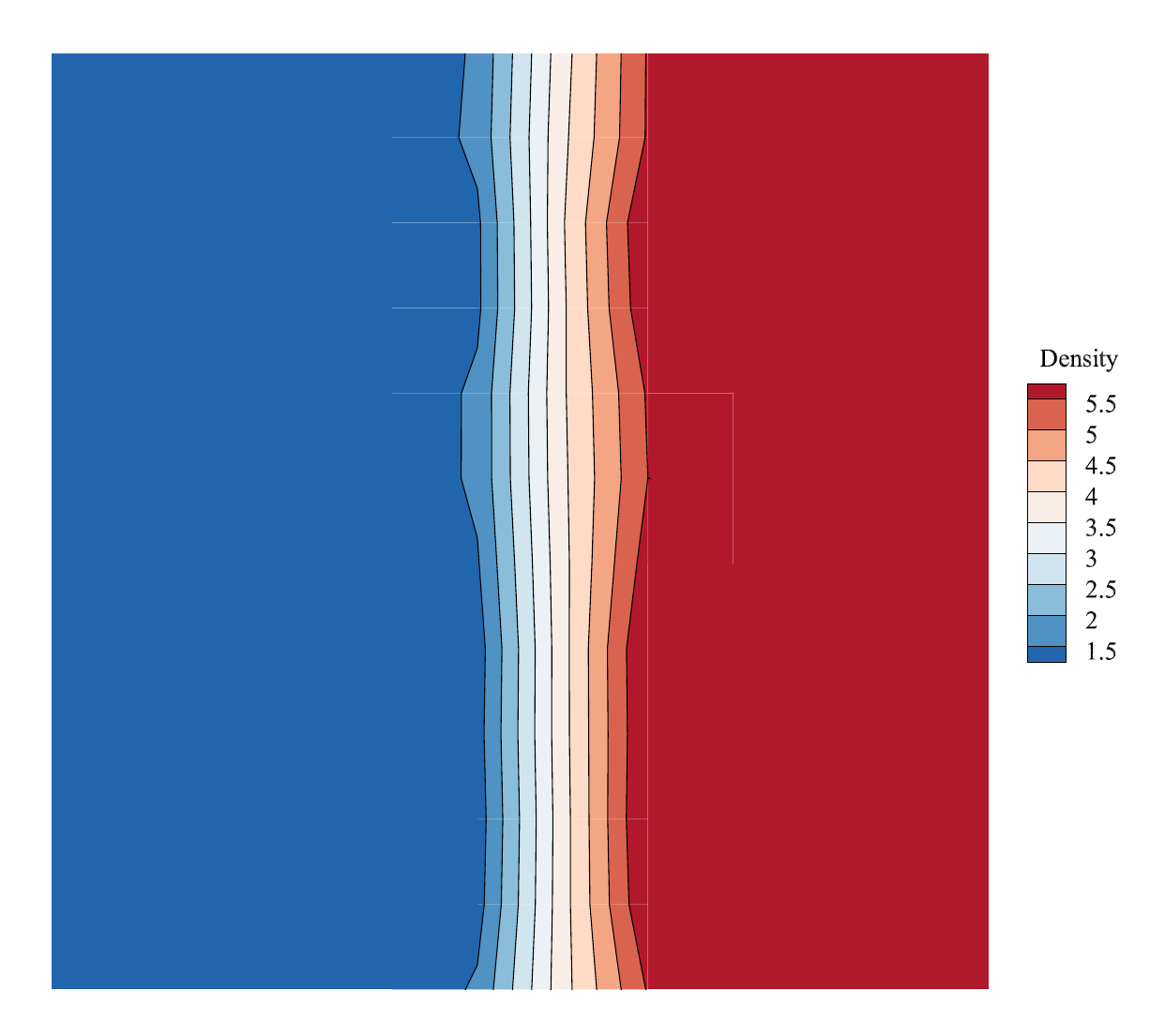}
	\end{minipage}
	}

	\centering
	\caption{Analysis of the fifth-order scheme with HLLC-ADC solver.(Grid with 11$ \times $11 cells, $ M_0=20 $ and $ \varepsilon =0.1 $.)}\label{fig HLLC-ADC}
\end{figure}

\begin{figure}[htbp]
	\centering

	\subfigure[distribution of the eigenvalues in the complex plane]{
	\begin{minipage}[t]{0.46\linewidth}
	\centering
	\includegraphics[width=0.85\textwidth]{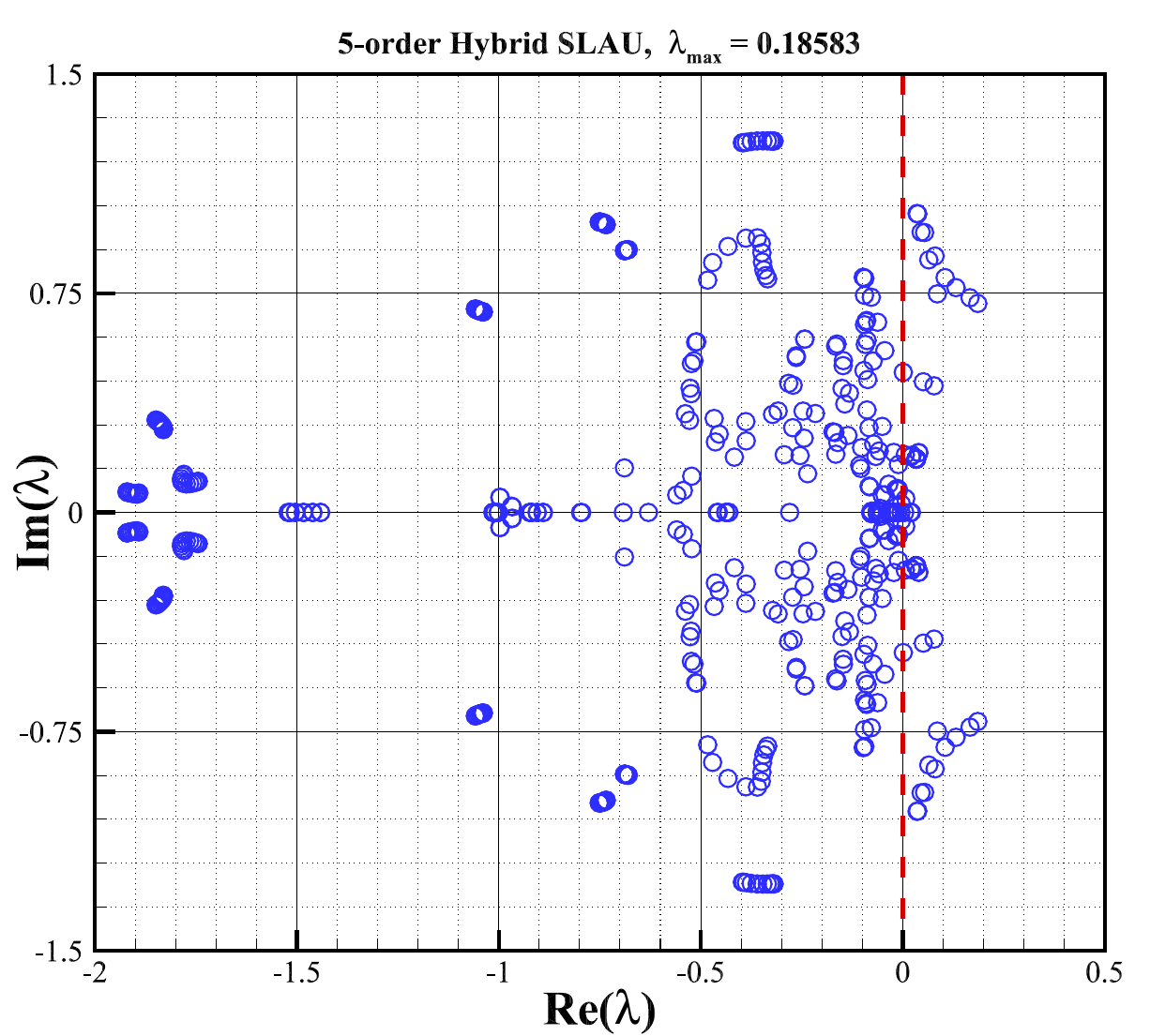}
	\end{minipage}
	}
	\subfigure[density contour]{
	\begin{minipage}[t]{0.46\linewidth}
	\centering
	\includegraphics[width=0.85\textwidth]{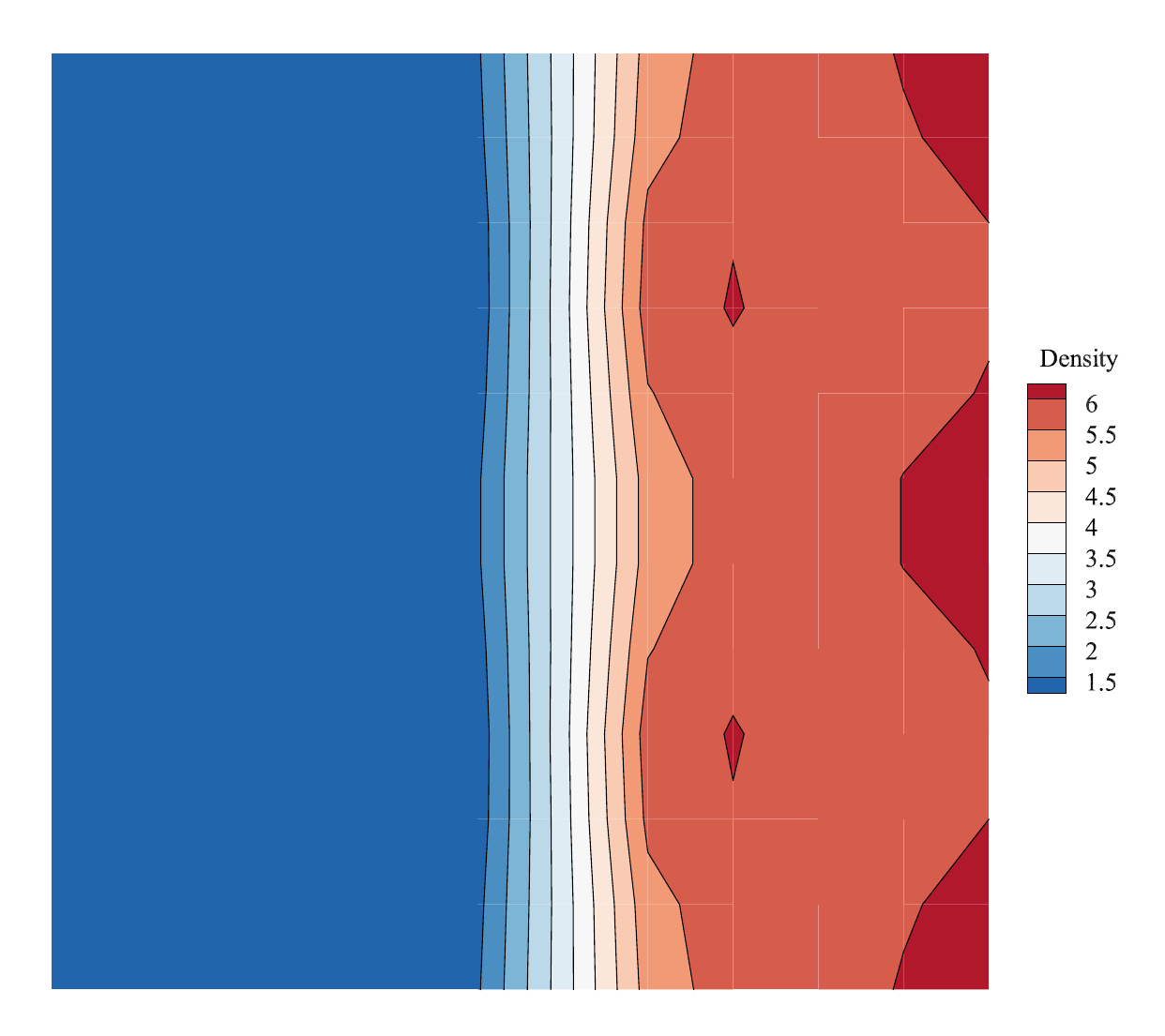}
	\end{minipage}
	}

	\centering
	\caption{Analysis of the fifth-order scheme with the hybrid SLAU solver.(Grid with 11$ \times $11 cells, $ M_0=20 $ and $ \varepsilon =0.1 $.)}\label{fig hybrid SLAU}

\end{figure}

\begin{figure}[htbp]
	\centering

	\subfigure[HLLC-ADC]{
	\begin{minipage}[t]{0.3\linewidth}
	\centering
	\includegraphics[width=0.9\textwidth]{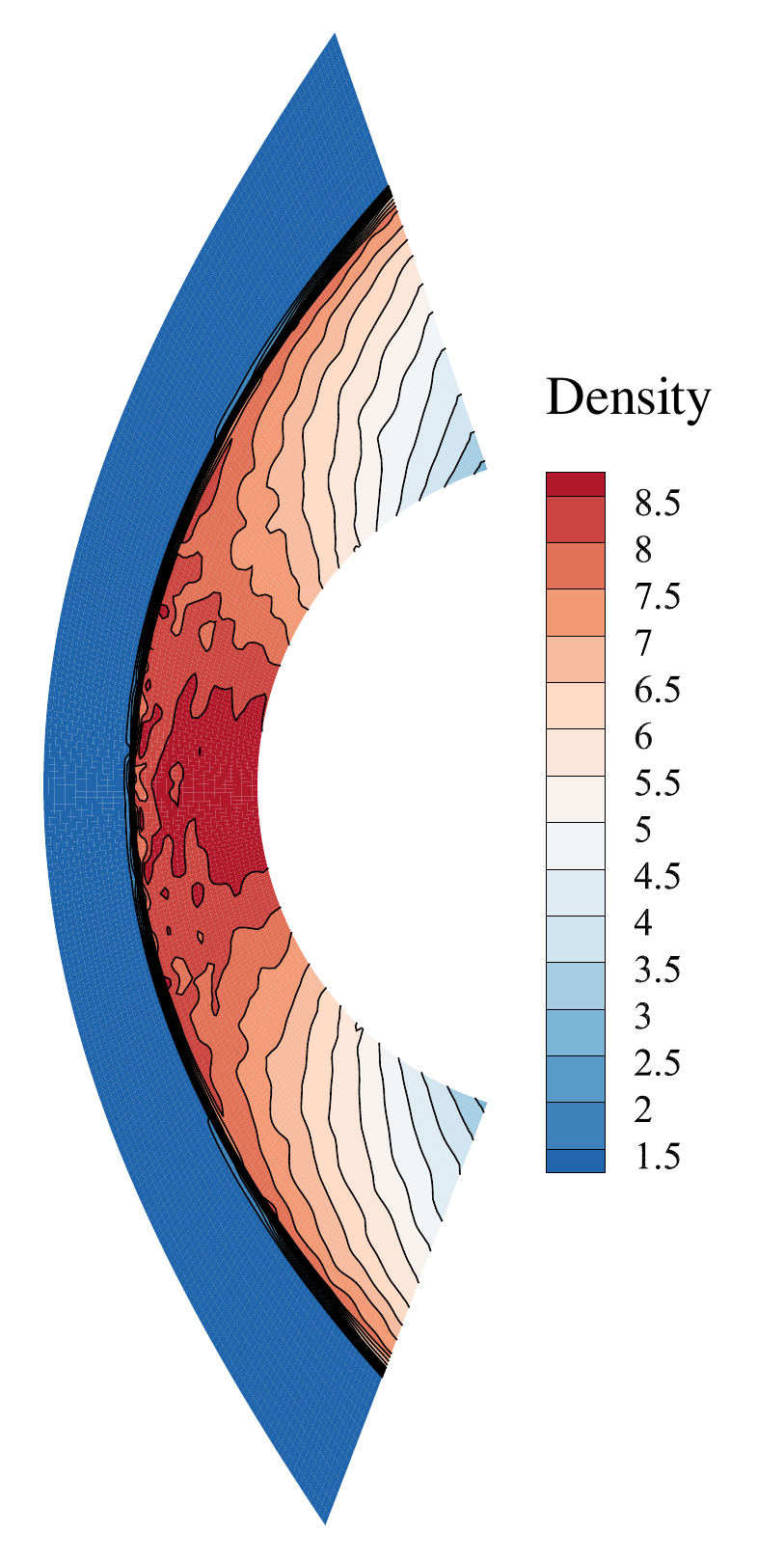}
	\end{minipage}
	}
	\subfigure[HLLC-EC]{
	\begin{minipage}[t]{0.3\linewidth}
	\centering
	\includegraphics[width=0.9\textwidth]{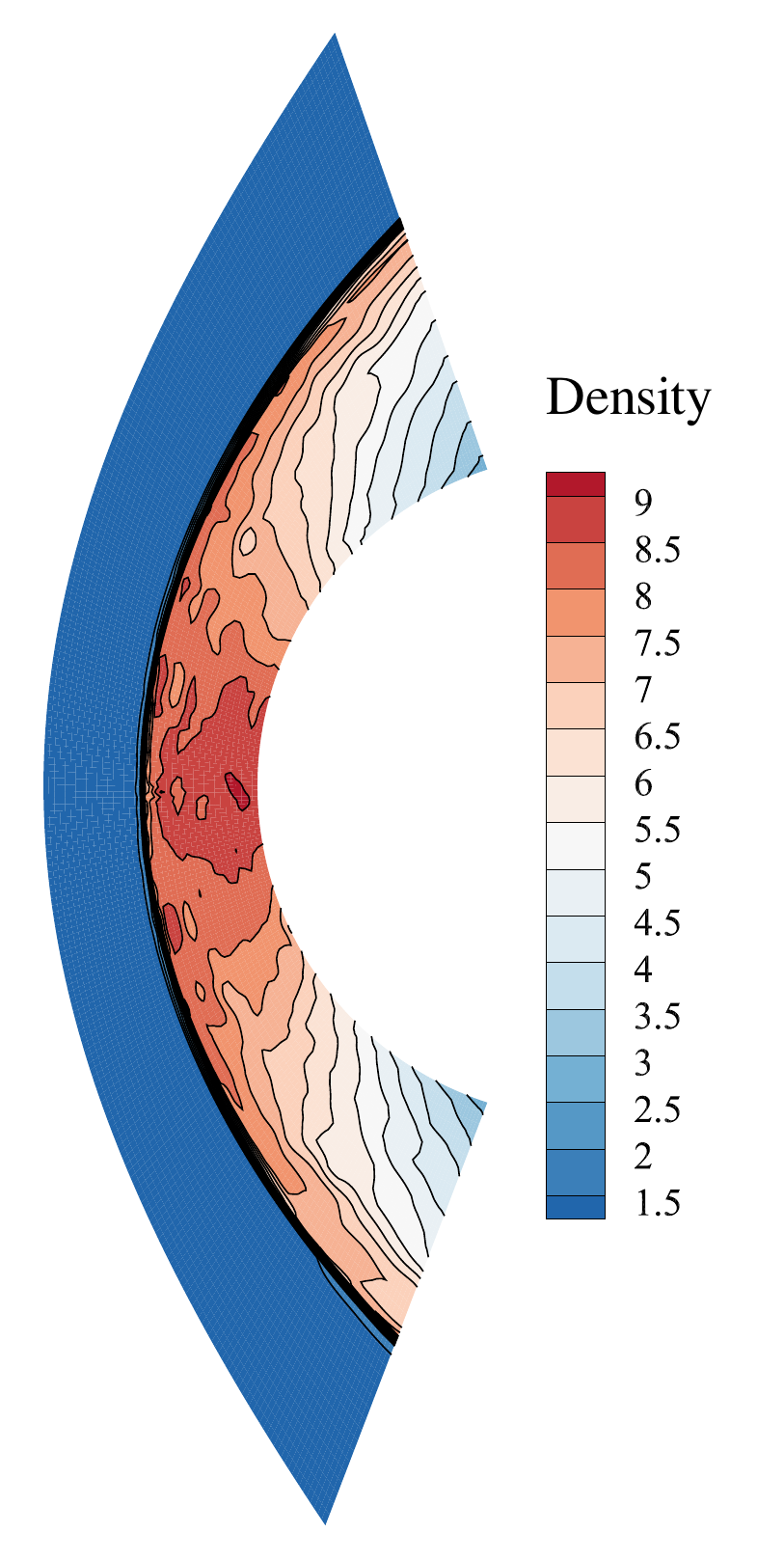}
	\end{minipage}
	}
  \subfigure[Hybrid SLAU]{
	\begin{minipage}[t]{0.3\linewidth}
	\centering
	\includegraphics[width=0.9\textwidth]{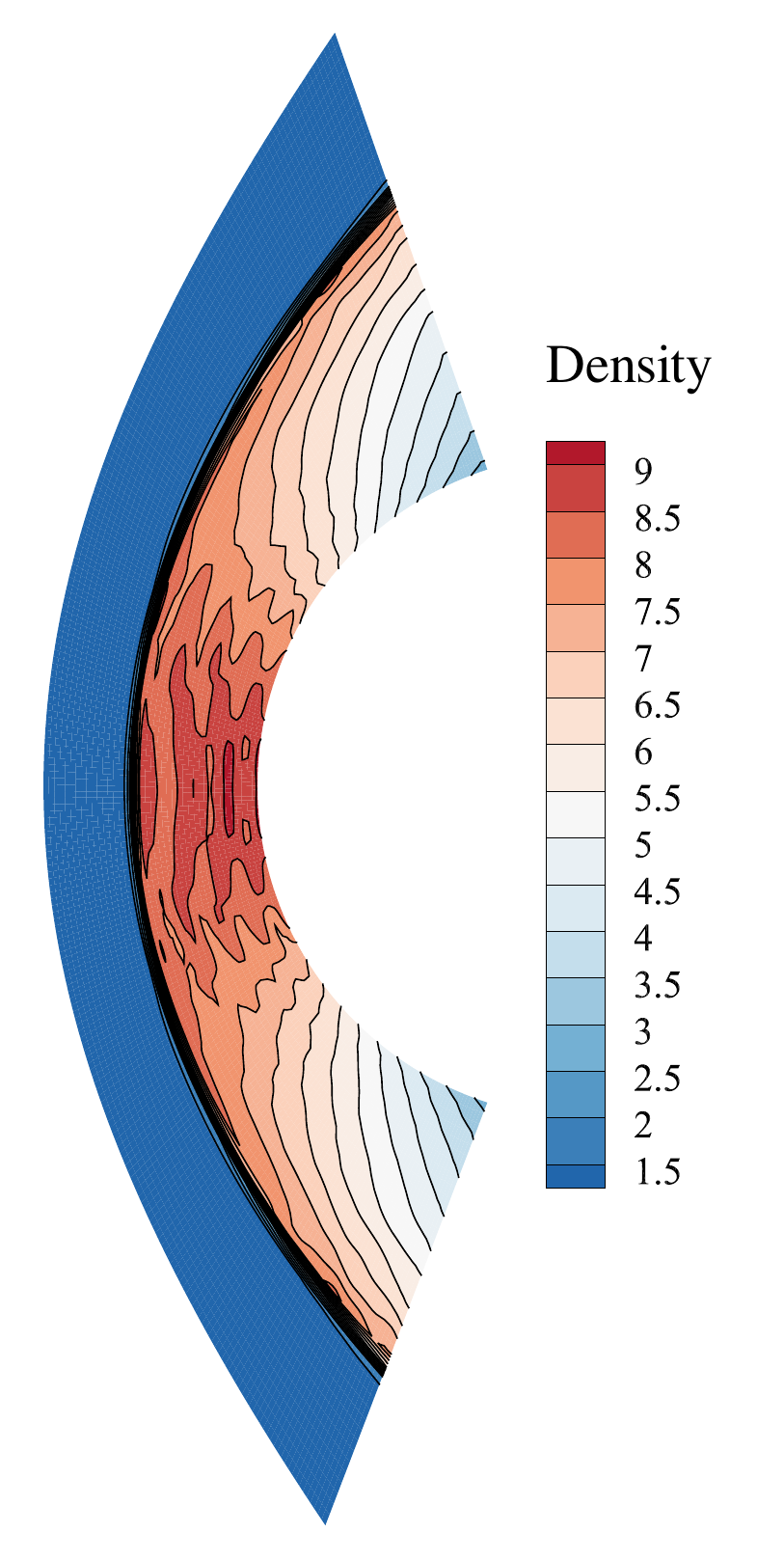}
	\end{minipage}
	}

	\centering
	\caption{The density contours of hypersonic flow over a cylinder computed by fifth-order schemes with different solvers.}\label{fig cylinder 2}
\end{figure}

Since the fifth-order scheme will suffer the shock instability problem even with dissipative solvers, it can be deduced that many methods aiming to capture shock stably will also have the risk of shock instabilities when the spatial accuracy is increased to the fifth-order. Examples of such methods include the HLLC-ADC \cite{Simon_Cure_2018}, hybrid SLAU \cite{Zhang_Robust_2017}, and the HLLC-EC \cite{Xie_Further_2021} and so on. These methods have been shown to achieve good results in first and second-order accuracy calculations. However, the matrix analysis presented in Fig.\ref{fig HLLC-EC} - Fig.\ref{fig hybrid SLAU} reveals that when these methods are embedded into the fifth-order scheme, the eigenvalues still have positive real parts, indicating their failure in maintaining stability at fifth-order accuracy. This conclusion is further supported by the corresponding density contours shown in Fig.\ref{fig HLLC-EC} - Fig.\ref{fig cylinder 2}, where unstable results are clearly observed.

\subsubsection{A further discussion on WENO reconstruction and shock instability}

\begin{figure}[htbp]
	\centering
	\includegraphics[width=0.6\textwidth]{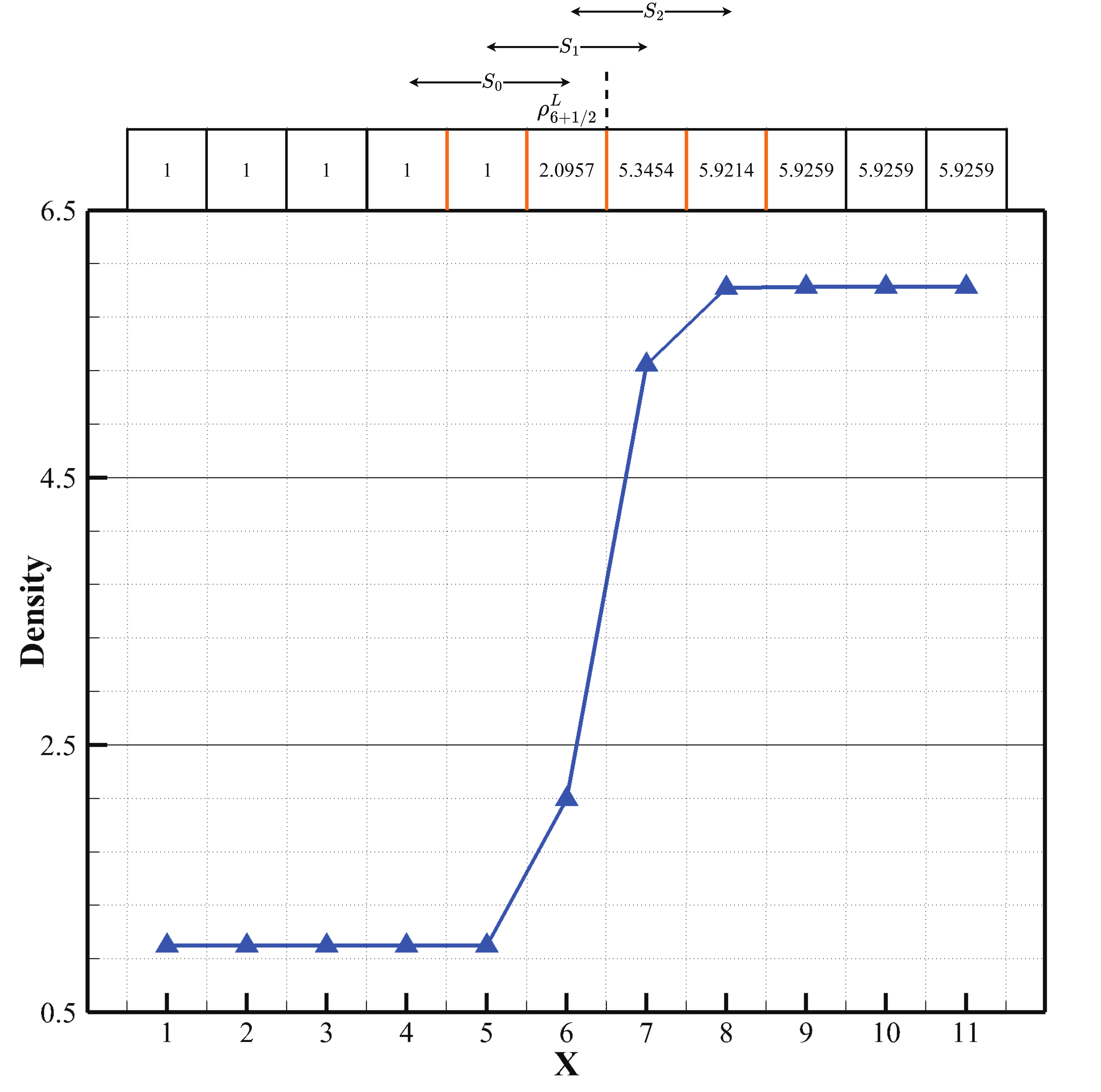}
	\caption{The initial data of the fifth-order scheme with HLL solver and stencils of the reconstruction for $\rho_{6+1/2}^L$.(The values in the table is the density in different cells.)}
	\label{fig InitialData HLL}
\end{figure}

\begin{figure}[htbp]
	\centering

	\subfigure[distribution of the eigenvalues in the complex plane]{
	\begin{minipage}[t]{0.46\linewidth}
	\centering
	\includegraphics[width=0.9\textwidth]{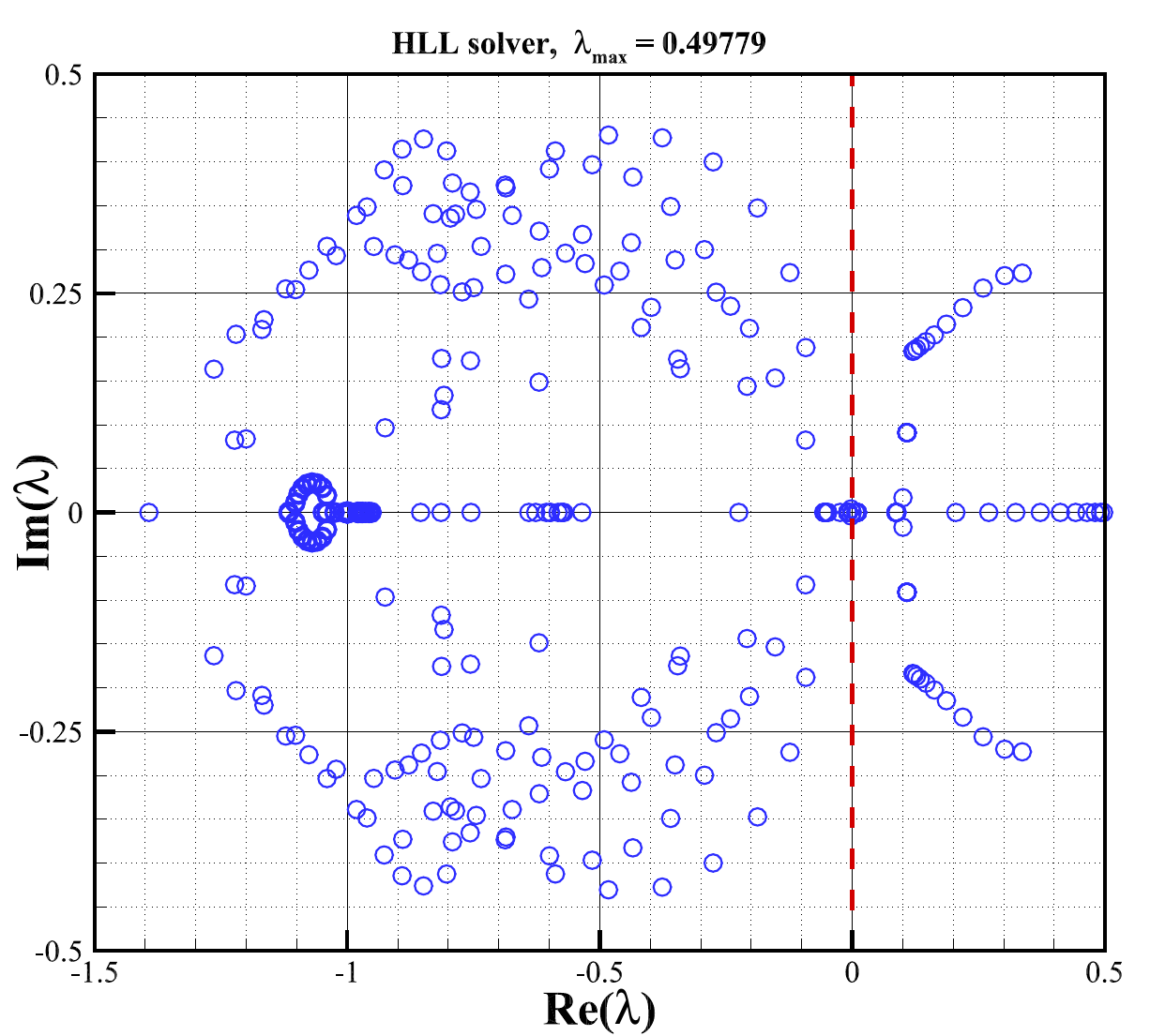}
	\end{minipage}
	}
	\subfigure[density contour]{
	\begin{minipage}[t]{0.46\linewidth}
	\centering
	\includegraphics[width=0.9\textwidth]{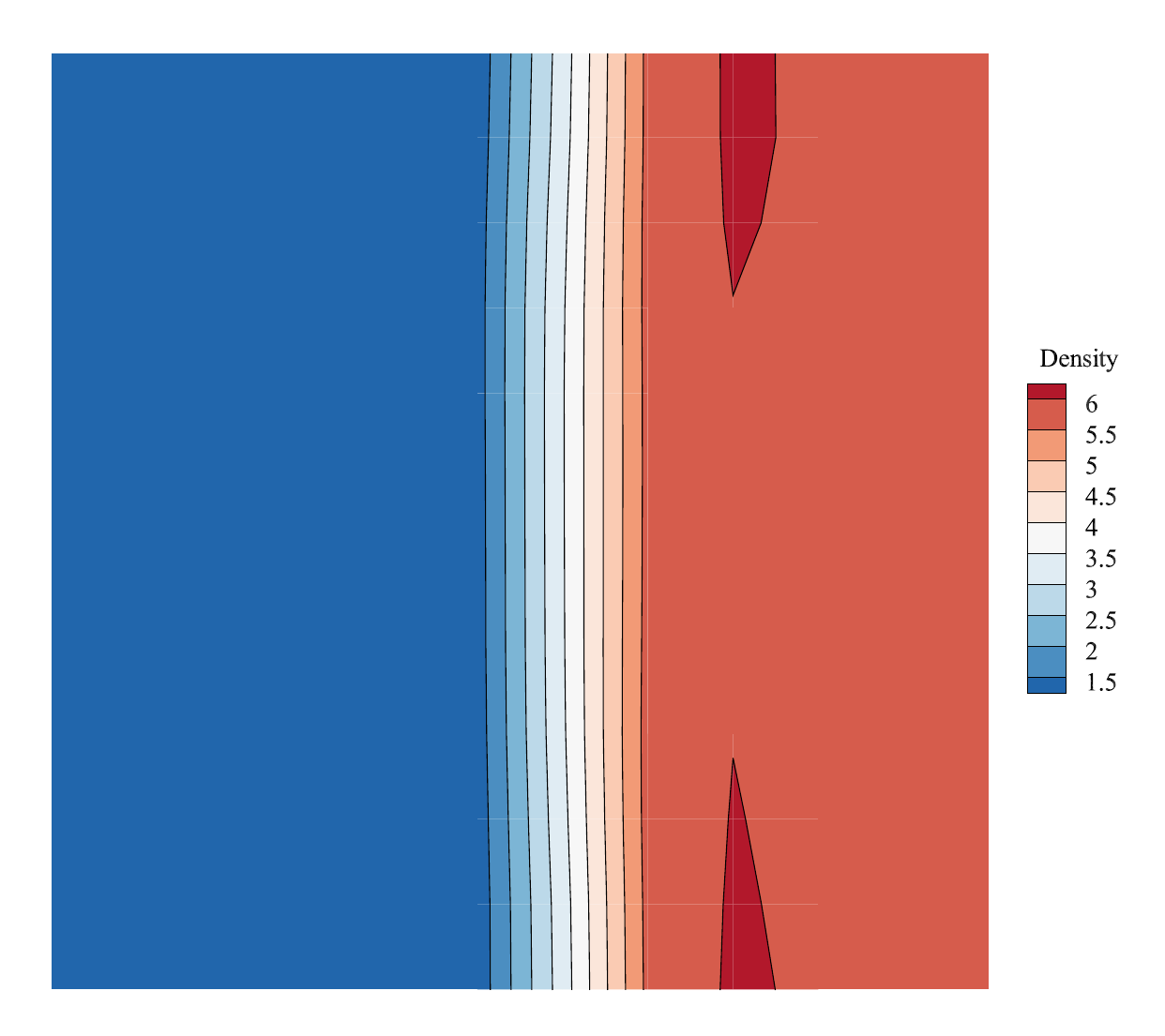}
	\end{minipage}
	}

	\centering
	\caption{Resluts of the HLL solver when the faces near the numerical shock structure employ third-order scheme.(Grid with 11$ \times $11 cells, $ M_0=20 $ and $ \varepsilon =0.1 $.)}\label{fig HLL_ENO}
\end{figure}

\begin{figure}[htbp]
	\centering

	\subfigure[distribution of the eigenvalues in the complex plane]{
	\begin{minipage}[t]{0.46\linewidth}
	\centering
	\includegraphics[width=0.9\textwidth]{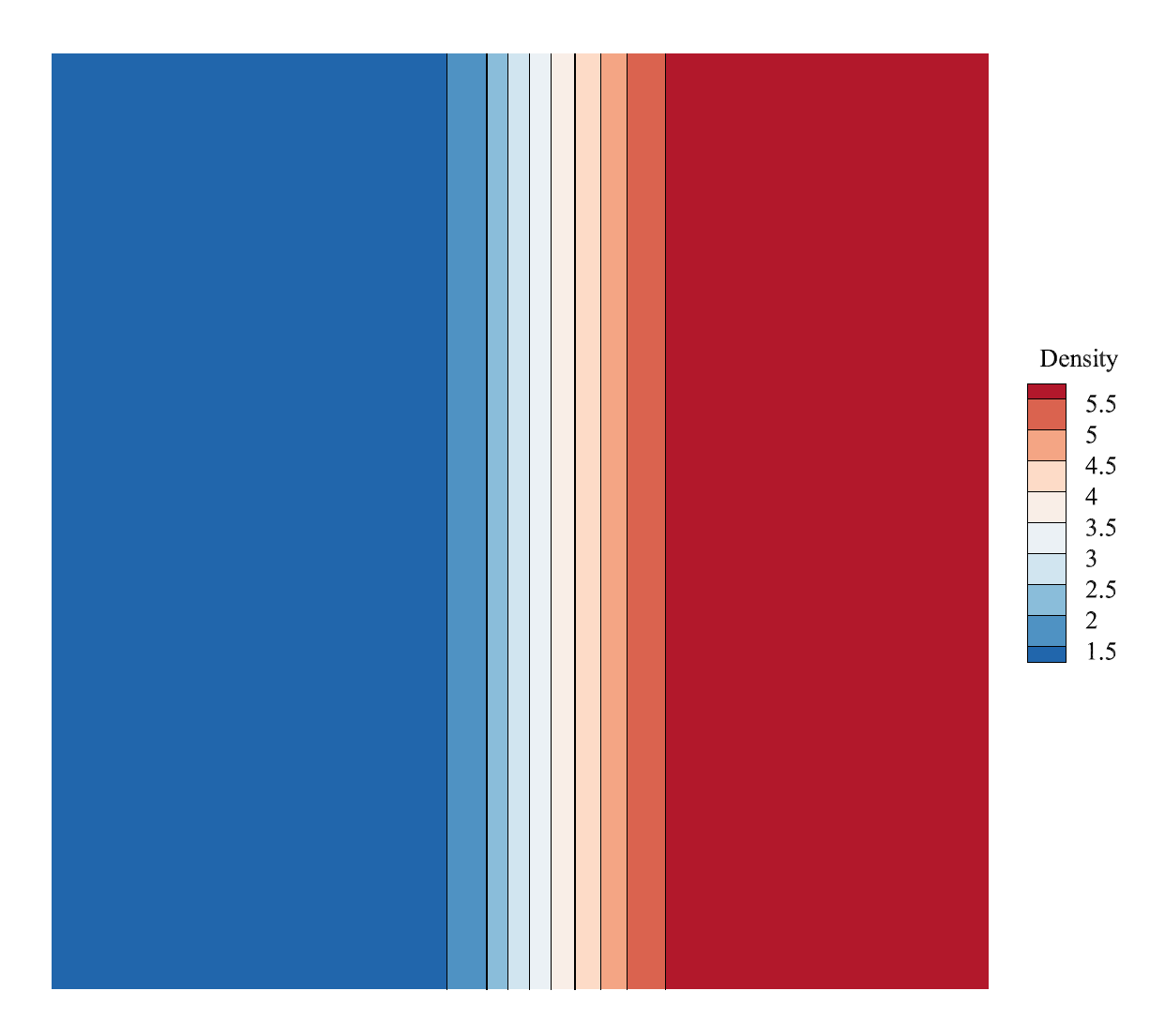}
	\end{minipage}
	}
	\subfigure[density contour]{
	\begin{minipage}[t]{0.46\linewidth}
	\centering
	\includegraphics[width=0.9\textwidth]{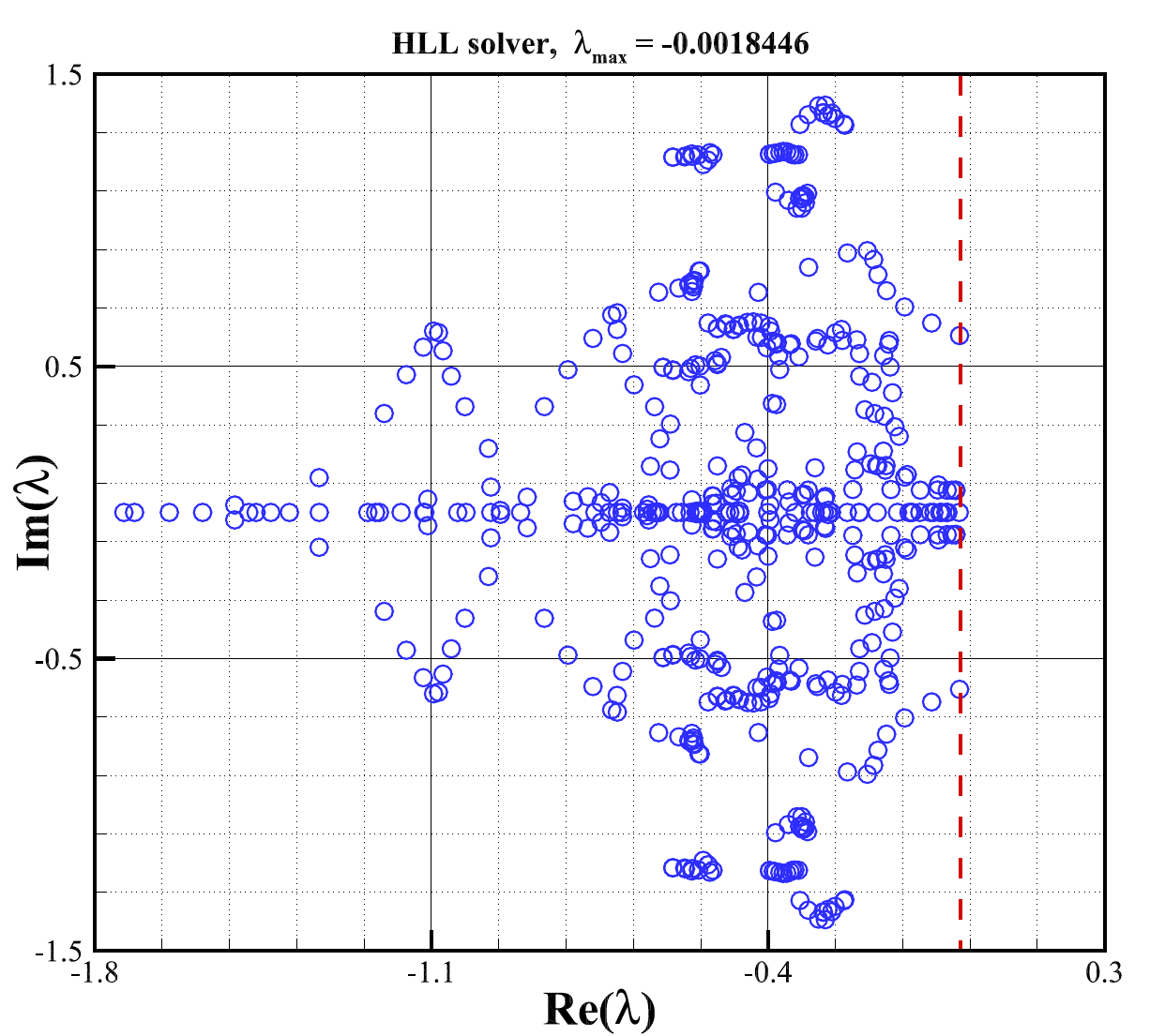}
	\end{minipage}
	}

	\centering
	\caption{Resluts of the HLL solver when the faces near the numerical shock structure employ first-order scheme.(Grid with 11$ \times $11 cells, $ M_0=20 $ and $ \varepsilon =0.1 $.)}\label{fig hybrid HLL}
\end{figure}

The matrix stability analysis and numerical experiments both indicate that in the case of fifth-order schemes, even dissipative schemes will be susceptible to shock instabilities. In the following analysis, we will further investigate the reasons behind this phenomenon using both matrix analysis and numerical experiments. Firstly, the core idea of the WENO method is to construct an interpolation function by combining candidate functions using nonlinear weights. These weights play a crucial role in the performance of the WENO scheme. They are specifically designed to achieve the highest-order linear scheme in smooth regions while still maintaining the ENO property in non-smooth regions. However, when reconstructing variables within the numerical shock structure, such as $ \rho_{6+i/2}^L $ in Fig.\ref{fig InitialData HLL}, all third-order stencils will include at least one discontinuity. Consequently, all nonlinear weights are nonzero (as shown in Table \ref{tab smooth factors and nonlinear weights of different stencils}), implying that all three third-order functions are utilized in the reconstruction of $ \rho_{6+i/2}^L $. As a result, the presence of the numerical shock structure prevents the WENO scheme from maintaining its ENO property near shocks. Secondly, Fig.\ref{fig HLL_ENO} demonstrates that even when using the smoothest third-order ENO stencils to reconstruct variables on both sides of the faces within the numerical shock structure (represented by the orange interfaces in Fig.\ref{fig InitialData HLL}), the computation still fails to stably capture the shocks. Conversely, the results in Fig.\ref{fig hybrid HLL} show that employing the first-order scheme on these faces yields stable results. Similarly, stable results can also be obtained when employing the second-order scheme on these faces. Therefore, it can be inferred that the third-order spatial accuracy near the numerical shock structure is too high to maintain stability.

\begin{table}[!ht]
    \centering
    \caption{The smooth factors and nonlinear weights of different stencils for the reconstruction of $\rho_{6+1/2}^L$.}
    \begin{tabular}{cccc}
		\toprule
        third-order stencils & $ S_0 $ & $ S_1 $ & $ S_2 $ \\ 
		\midrule
        smooth factors       & 4.00186 & 9.74719 & 28.78125  \\ 
        nonlinear weights    & 0.21135 & 0.62458 & 0.16407  \\ 
		\bottomrule
    \end{tabular}
    \label{tab smooth factors and nonlinear weights of different stencils}
\end{table}

Based on the analysis in this section, it can be concluded that the presence of the numerical shock structure causes the WENO scheme to lose its ENO property near shocks. Additionally, even the smoothest ENO stencil is used to reconstruct the variables within the shock structure, the spatial accuracy is still too high to stably capture the shock. These issues result in the inability of dissipative solvers in fifth-order WENO schemes to stably capture strong shocks. These findings provide valuable insights for the development of new high-order schemes. It is important to consider the influence of numerical shock structure. In particular, it is necessary to sacrifice accuracy near shocks and limit the order of accuracy to first- or second-order in the vicinity of shocks in order to avoid instabilities.

\subsection{Sensibility of cell face to shock instabilities}\label{subsection 4.2}

\begin{figure}[htbp]
	\centering
	\includegraphics[width=0.5\textwidth]{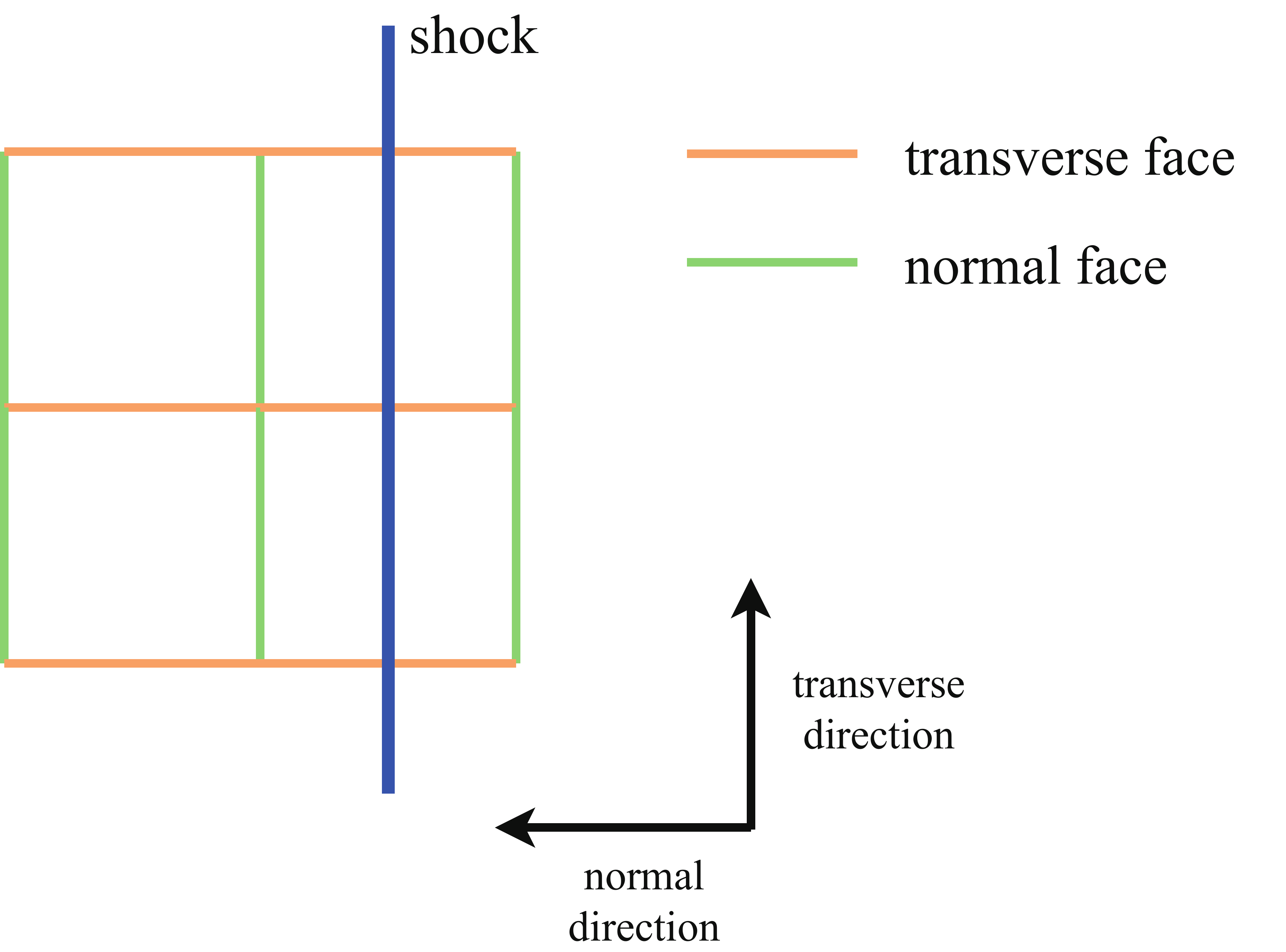}
	\caption{Diagram of transverse and normal faces near the shock.}
	\label{fig direction}
\end{figure}

It has been demonstrated that for low-dissipative solvers such as Roe and HLLC, shock instabilities originate from the transverse faces rather than the normal faces (identified as orange and green faces in Fig.\ref{fig direction}) near the shock \cite{Shen_Stability_2014, Simon_Cure_2018, Chen_Mechanism_2018, Xu_Does_1999, Xu_Does_1999}. This phenomenon is because there is less dissipation in the transverse faces compared to the normal faces \cite{Xu_Dissipative_2005, Xu_Does_1999}. Therefore, if there is sufficient dissipation in the transverse faces, the computation can remain stable. However, these studies primarily focus on capturing thin shocks without considering the numerical shock structure. In \cite{Ren_Numerical_2023}, it has been revealed that the shock instability problem is a multidimensional coupling problem. To stably capture strong shocks, there must be sufficient dissipation in transverse faces and more than one point within the numerical shock structure in the direction perpendicular to the shock. It is challenging to stably capture shocks completely if only one of these conditions is satisfied. Such a conclusion has been obtained for low-order schemes, and in this section, the matrix stability analysis method is employed to investigate the multidimensional coupling issue for fifth-order schemes.

\begin{figure}[htbp]
	\centering

	\subfigure[distribution of the eigenvalues in the complex plane]{
	\begin{minipage}[t]{0.46\linewidth}
	\centering
	\includegraphics[width=0.9\textwidth]{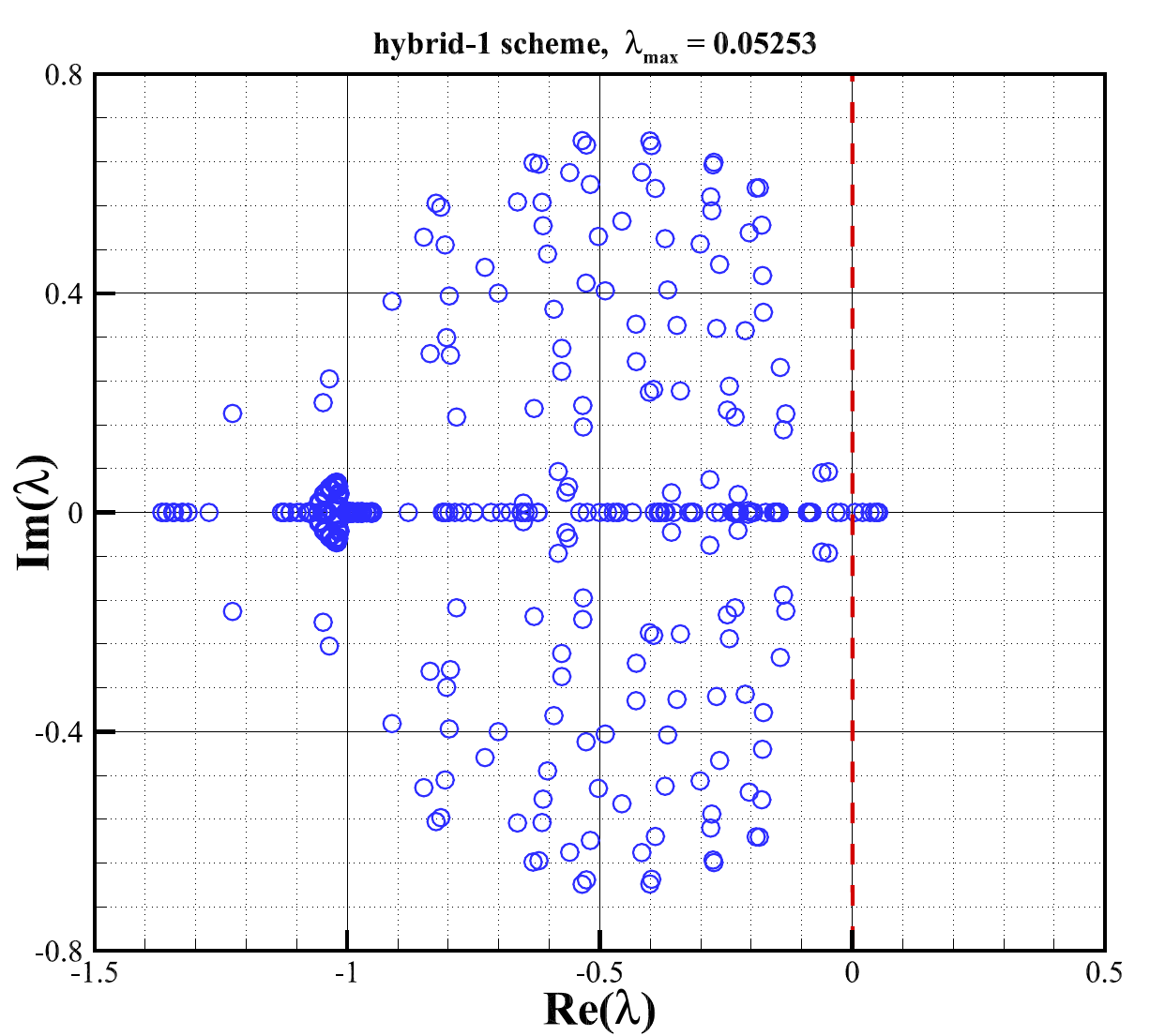}
	\end{minipage}
	}
	\subfigure[density contour]{
	\begin{minipage}[t]{0.46\linewidth}
	\centering
	\includegraphics[width=0.9\textwidth]{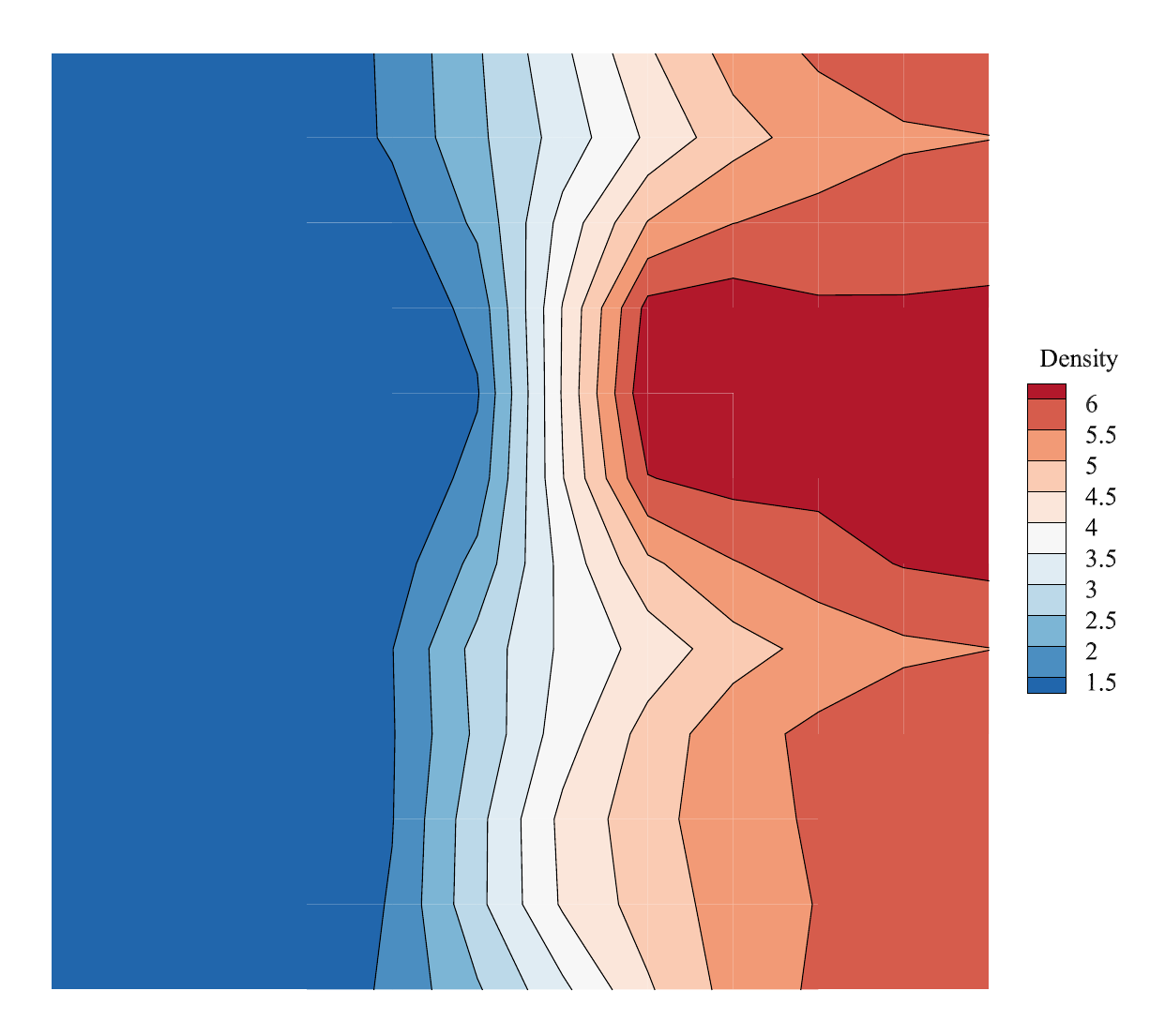}
	\end{minipage}
	}

	\centering
	\caption{Matrix stability analysis and the corresponding density contour of hybrid-1 scheme.(Grid with 11$ \times $11 cells, $ M_0=20 $ and $ \varepsilon =0.1 $.)}\label{fig hybrid-1}
\end{figure}

\begin{figure}[htbp]
	\centering

	\subfigure[distribution of the eigenvalues in the complex plane]{
	\begin{minipage}[t]{0.46\linewidth}
	\centering
	\includegraphics[width=0.9\textwidth]{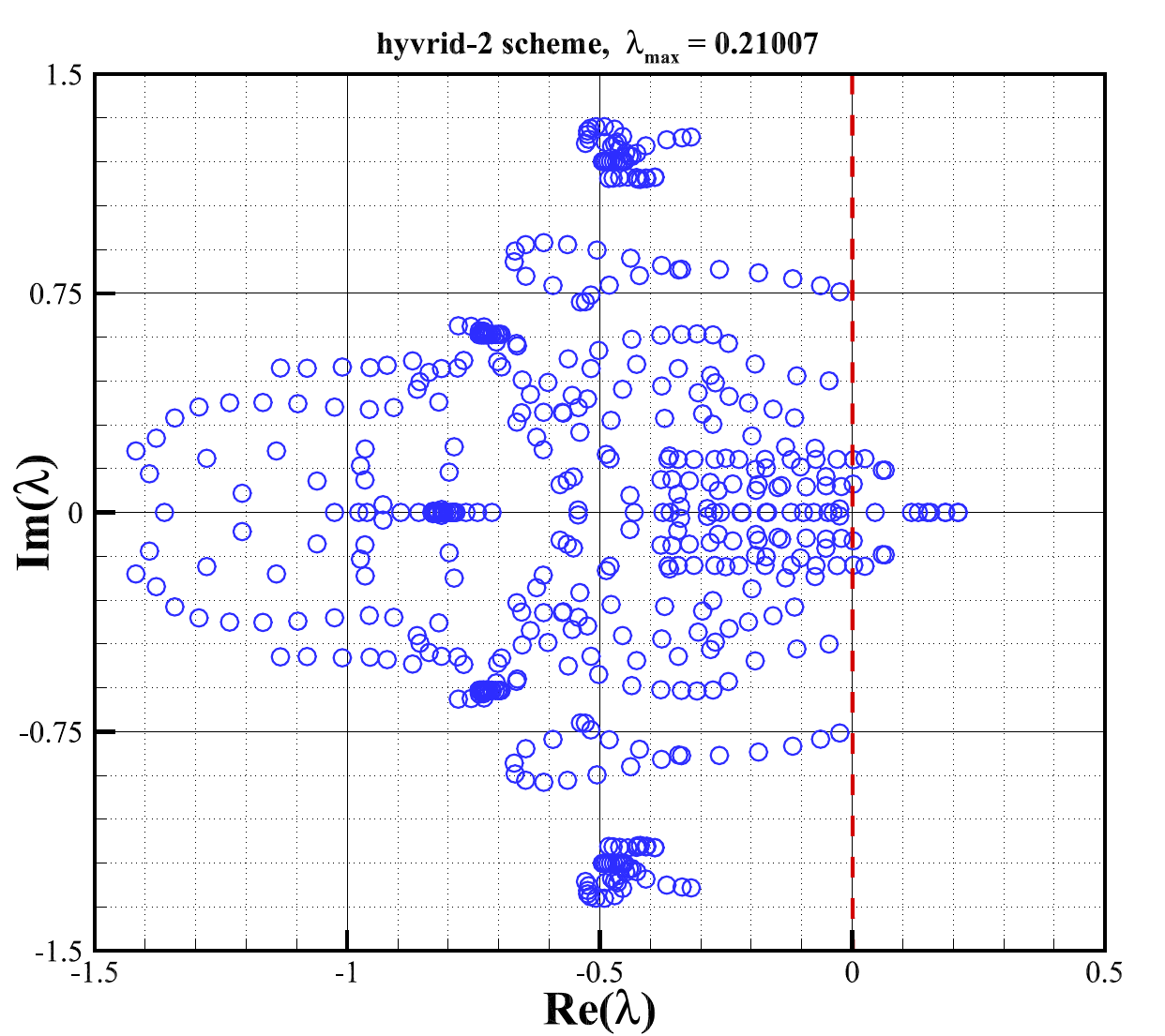}
	\end{minipage}
	}
	\subfigure[density contour]{
	\begin{minipage}[t]{0.46\linewidth}
	\centering
	\includegraphics[width=0.9\textwidth]{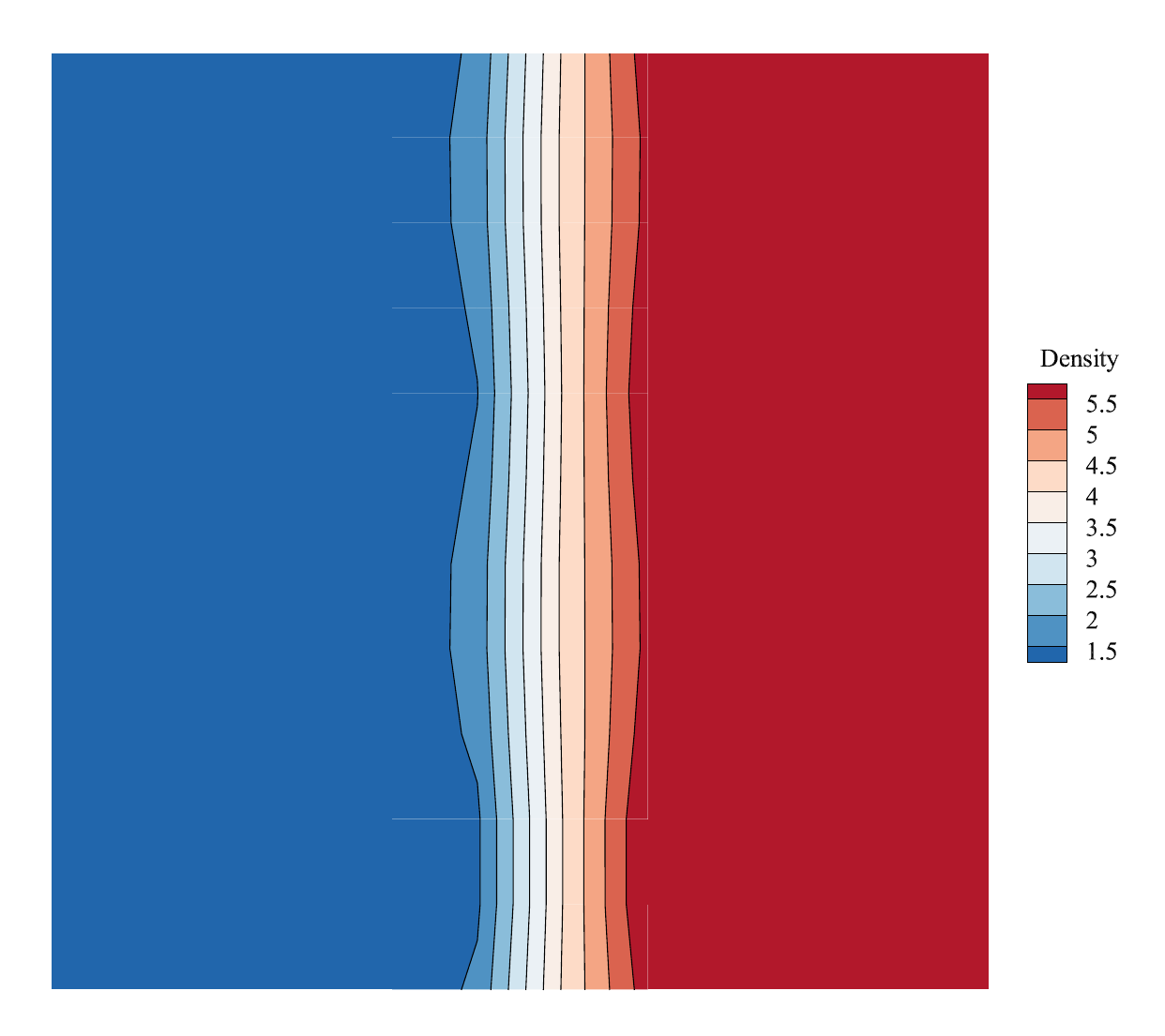}
	\end{minipage}
	}

	\centering
	\caption{Matrix stability analysis and the corresponding density contour of hybrid-2 scheme.(Grid with 11$ \times $11 cells, $ M_0=20 $ and $ \varepsilon =0.1 $.)}\label{fig hybrid-2}
\end{figure}

\begin{figure}[htbp]
	\centering

	\subfigure[hybrid-1 scheme]{
	\begin{minipage}[t]{0.3\linewidth}
	\centering
	\includegraphics[width=0.9\textwidth]{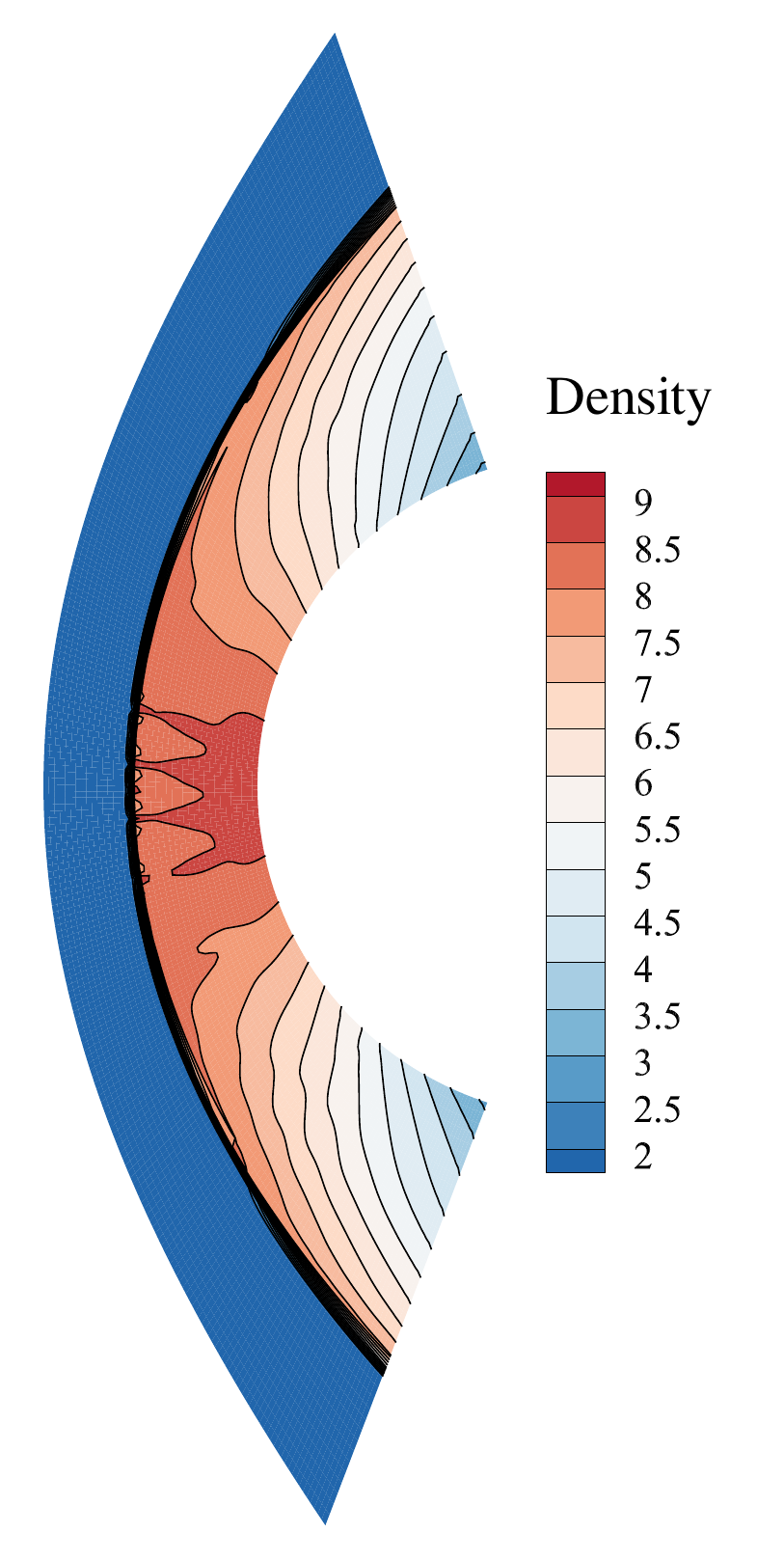}
	\end{minipage}
	}
    \subfigure[hybrid-2 scheme]{
	\begin{minipage}[t]{0.3\linewidth}
	\centering
	\includegraphics[width=0.9\textwidth]{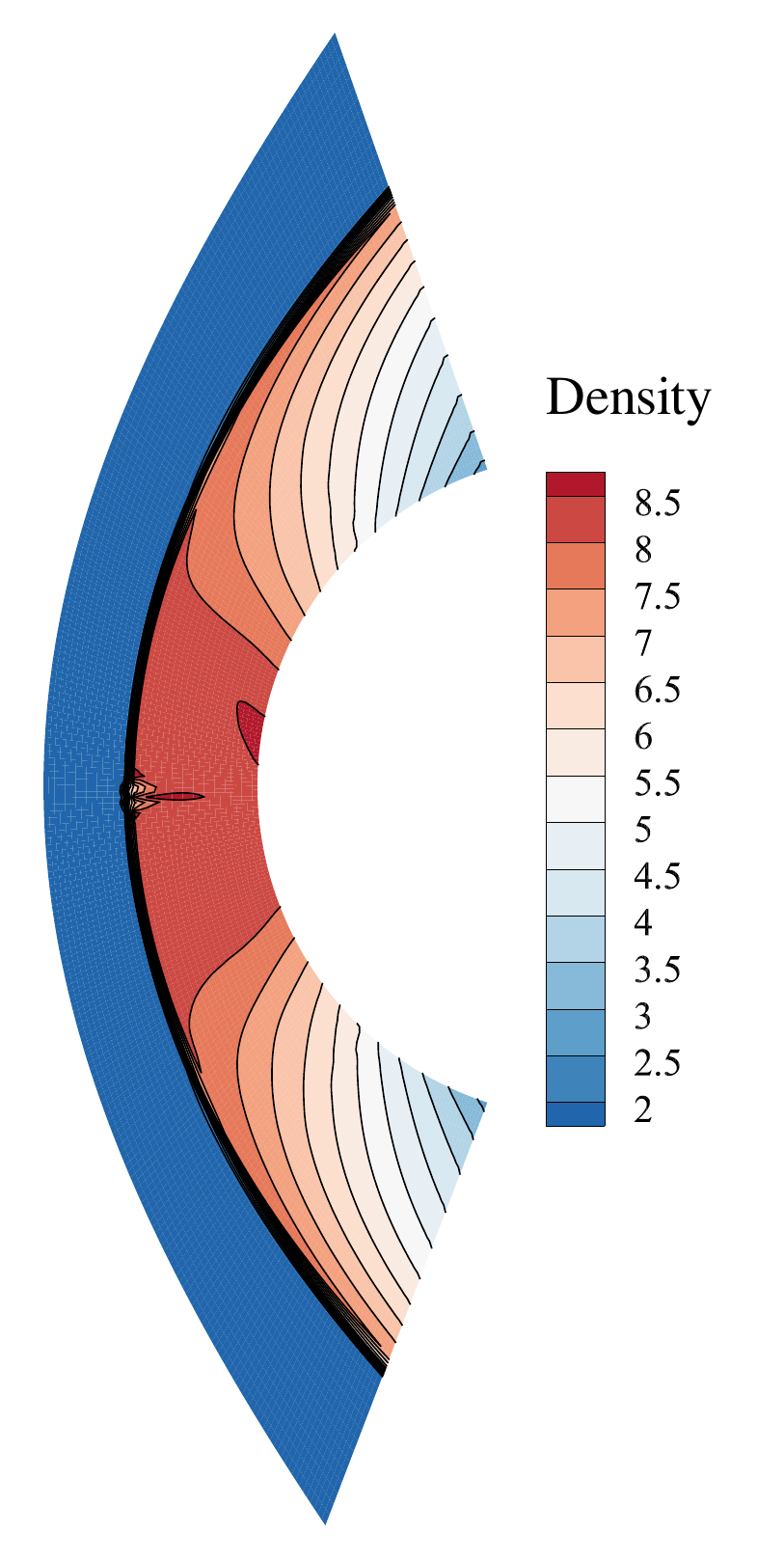}
	\end{minipage}
	}

	\centering
	\caption{The density contours of hypersonic flow over a cylinder computed the hybrid schemes.}\label{fig cylinder direction}
\end{figure}

To this end, we design two hybrid schemes, denoted as $\mathbf{F}_{\text{hybrid-1}}$ and $\mathbf{F}_{\text{hybrid-2}}$, which can be expressed as:
\begin{equation}
	\begin{aligned}
		& \mathbf{F}_{\text{hybrid-1}}= \begin{cases}\mathbf{F}_{\text{5th-order Roe}} & , \text{for transverse faces} \\
		\mathbf{F}_{\text{1st-order van Leer}} & , \text{for normal faces}\end{cases}, \\
		& \mathbf{F}_{\text{hybrid-2}}= \begin{cases}\mathbf{F}_{\text{1st-order van Leer}} & , \text{for transverse faces} \\
		\mathbf{F}_{\text{5th-order Roe}} & , \text{for normal faces}\end{cases}.
	\end{aligned}
\end{equation}
Here, $\mathbf{F}_{\text{hybrid-1/2}}$ represent the two hybrid flux schemes. $\mathbf{F}_{\text{5th-order Roe}}$ and $\mathbf{F}_{\text{1st-order van Leer}}$ correspond to the fifth-order scheme with the Roe solver and the first-order scheme with the van Leer solver, respectively. It can be observed that $\mathbf{F}_{\text{hybrid-1}}$ can produce a solution where there are more than one point within the numerical shock structure in the direction perpendicular to the shock, but with little dissipation on transverse faces. In contrast, if $\mathbf{F}_{\text{hybrid-2}}$ is used, there is sufficient dissipation on transverse faces, but only one point within the shock structure. By analyzing the stability of these two hybrid schemes, the multidimensional coupling problem of fifth-order schemes can be decoupled and easily analyzed. Fig.\ref{fig hybrid-1} and Fig.\ref{fig hybrid-2} show the results of matrix analysis and the corresponding density contours. As shown, both $\mathbf{F}_{\text{hybrid-1}}$ and $\mathbf{F}_{\text{hybrid/2}}$ have positive maximal real parts of all eigenvalues, suggesting that the two schemes will suffer from shock instabilities. Consistent with the analysis results, the corresponding density contours also display obvious instability. To validate the results further, the simulation of hypersonic flows over a blunt-body is performmed. The details of the simulation can be found in section \ref{subsection 4.1}. As shown in Fig.\ref{fig cylinder direction}, the computations become unstable regardless of the hybrid scheme employed. Therefore, it can be concluded from this section that the shock instability problem of fifth-order schemes is still a multidimensional coupling problem, just like in low-order cases. This conclusion provides guidance for developing robust high-order schemes, which should have sufficient dissipation on transverse faces and at least two points within the numerical shock structure in the direction perpendicular to the shock.

\subsection{Spatial localization of the source of shock instability}\label{subsection 4.3}
Investigating the spatial localization of the source of shock instability is of great significance. It can help to shed new light on understanding the mechanism of the shock instability problem and conduct targeted treatments for the location prone to instability. There are several studies aiming to investigate the spatial localization of the source of shock instability. In \cite{Dumbser_Matrix_2004,Chen_Mechanism_2018}, authors suppose that the shock instability originates from the upstream region and is convected downstream. However, the studies in \cite{Dumbser_Matrix_2004,Chen_Mechanism_2018} are based on the thin shock without the numerical shock structure, which is not common in actual simulations. In \cite{Xie_Numerical_2017,Ren_Numerical_2023}, investigations on the spatial location of shock instability with the numerical shock structure are carried out. They argue that the shock instability originates from the numerical shock structure. However, such a conclusion is dedeuced from studies for first- or second-order schemes. When the spatial accuracy is increased to fifth-order, does the instability still originate from the numerical shock structure? In this section, this problem will be investigated by matrix stability analysis and numerical experiments.

\begin{figure}[htbp]
	\centering

	\subfigure[$ \rho $]{
	\begin{minipage}[t]{0.46\linewidth}
	\centering
	\includegraphics[width=0.9\textwidth]{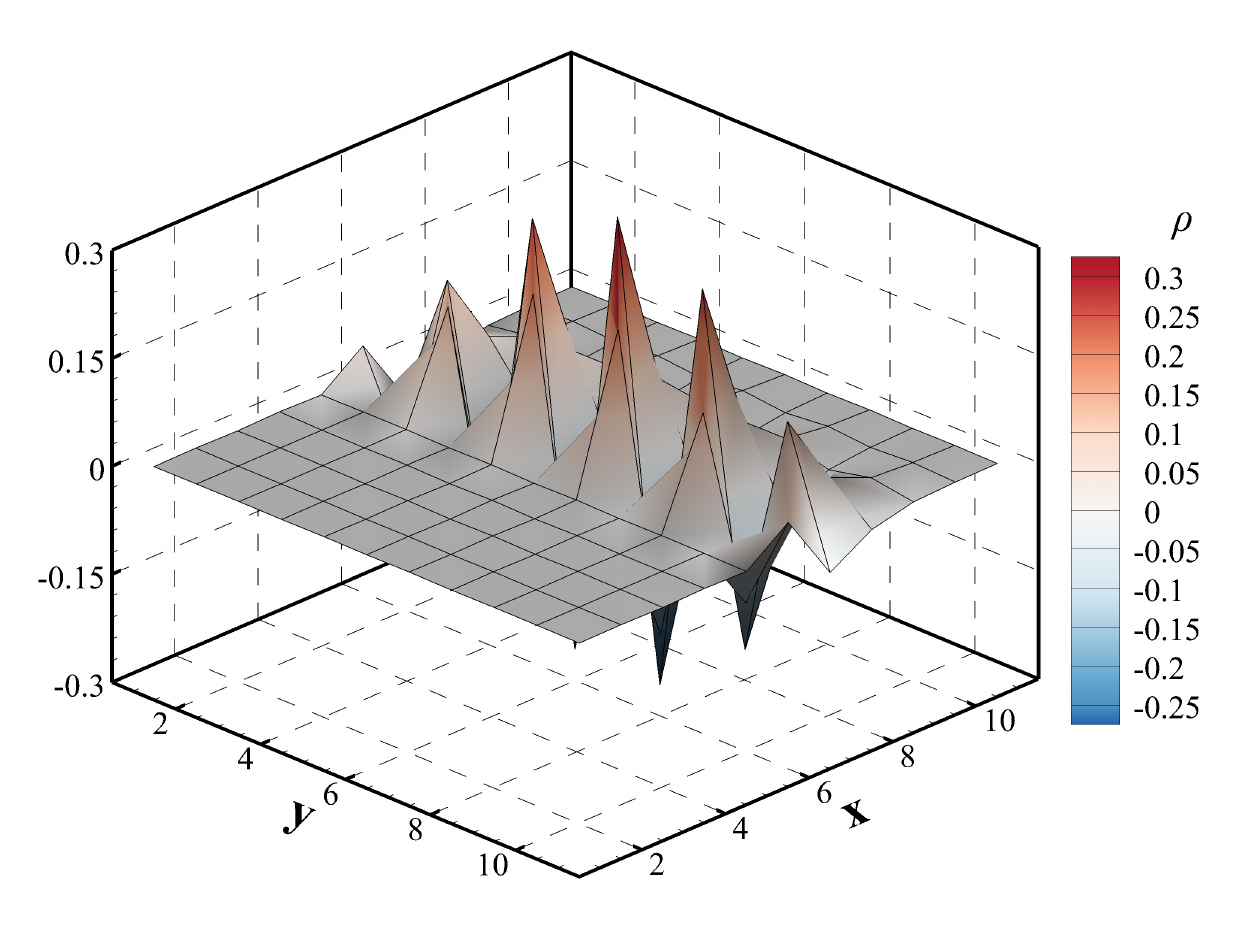}
	\end{minipage}
	}
	\subfigure[$ u $]{
	\begin{minipage}[t]{0.46\linewidth}
	\centering
	\includegraphics[width=0.9\textwidth]{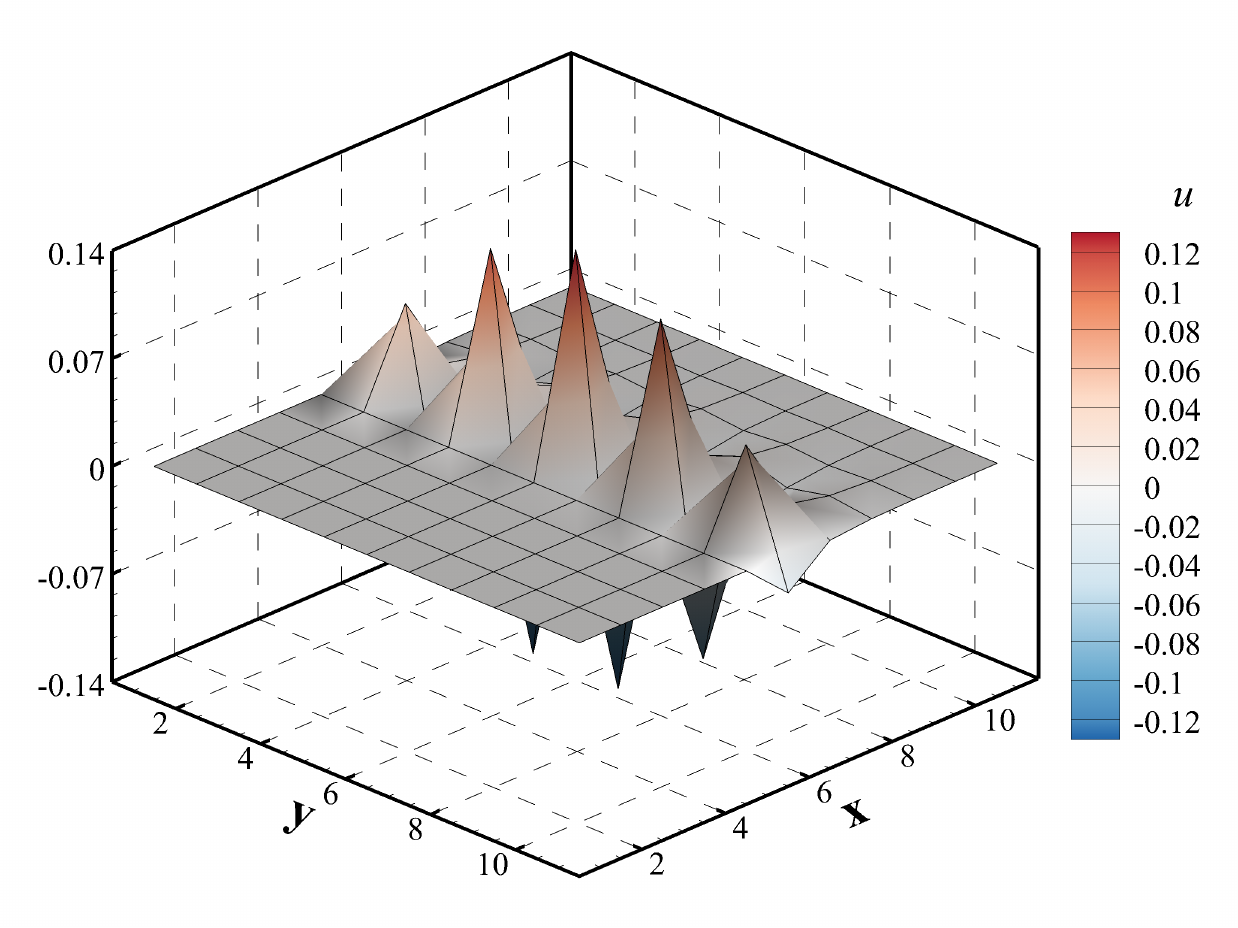}
	\end{minipage}
	}

	\subfigure[$ v $]{
	\begin{minipage}[t]{0.46\linewidth}
	\centering
	\includegraphics[width=0.9\textwidth]{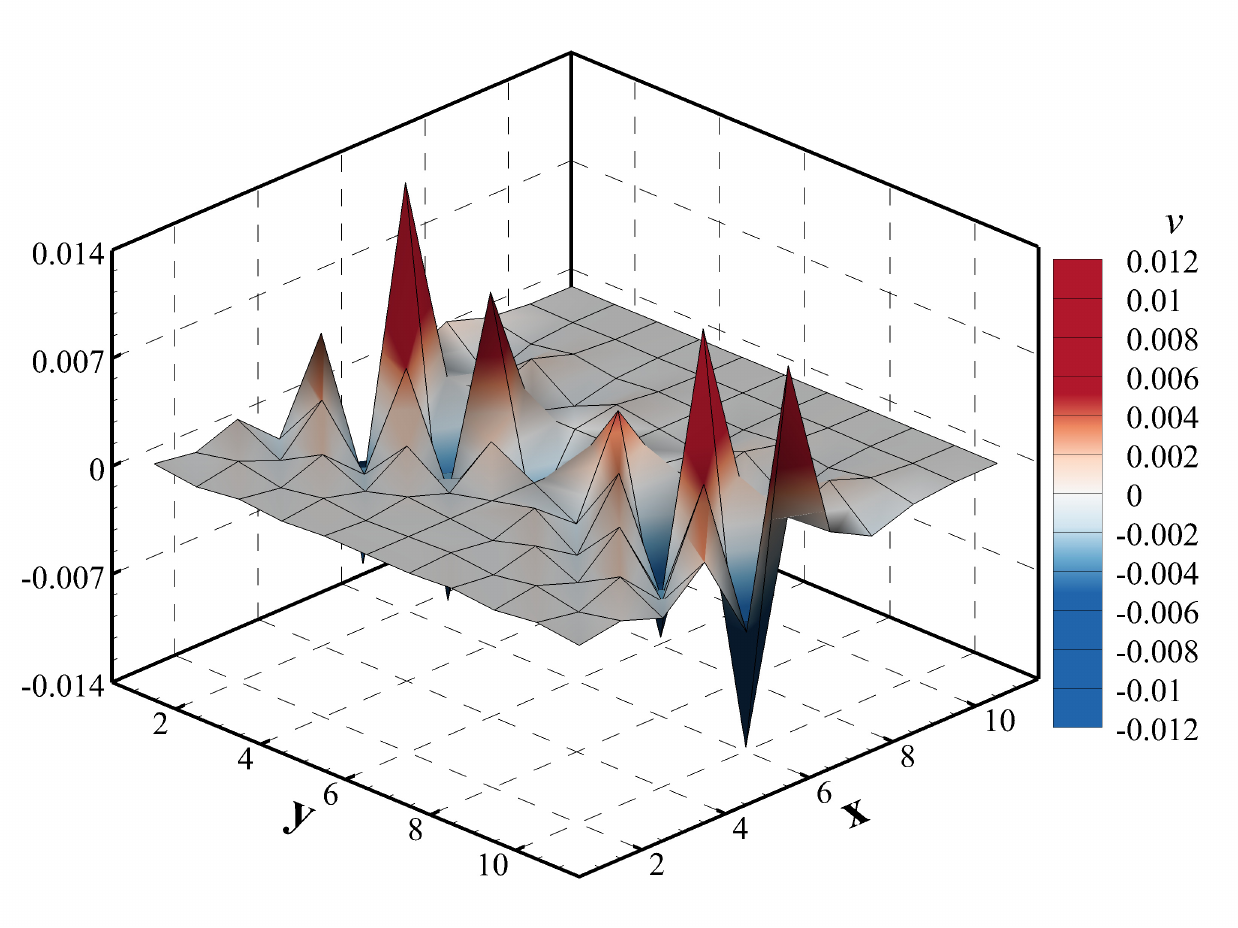}
	\end{minipage}
	}
	\subfigure[$ p $]{
	\begin{minipage}[t]{0.46\linewidth}
	\centering
	\includegraphics[width=0.9\textwidth]{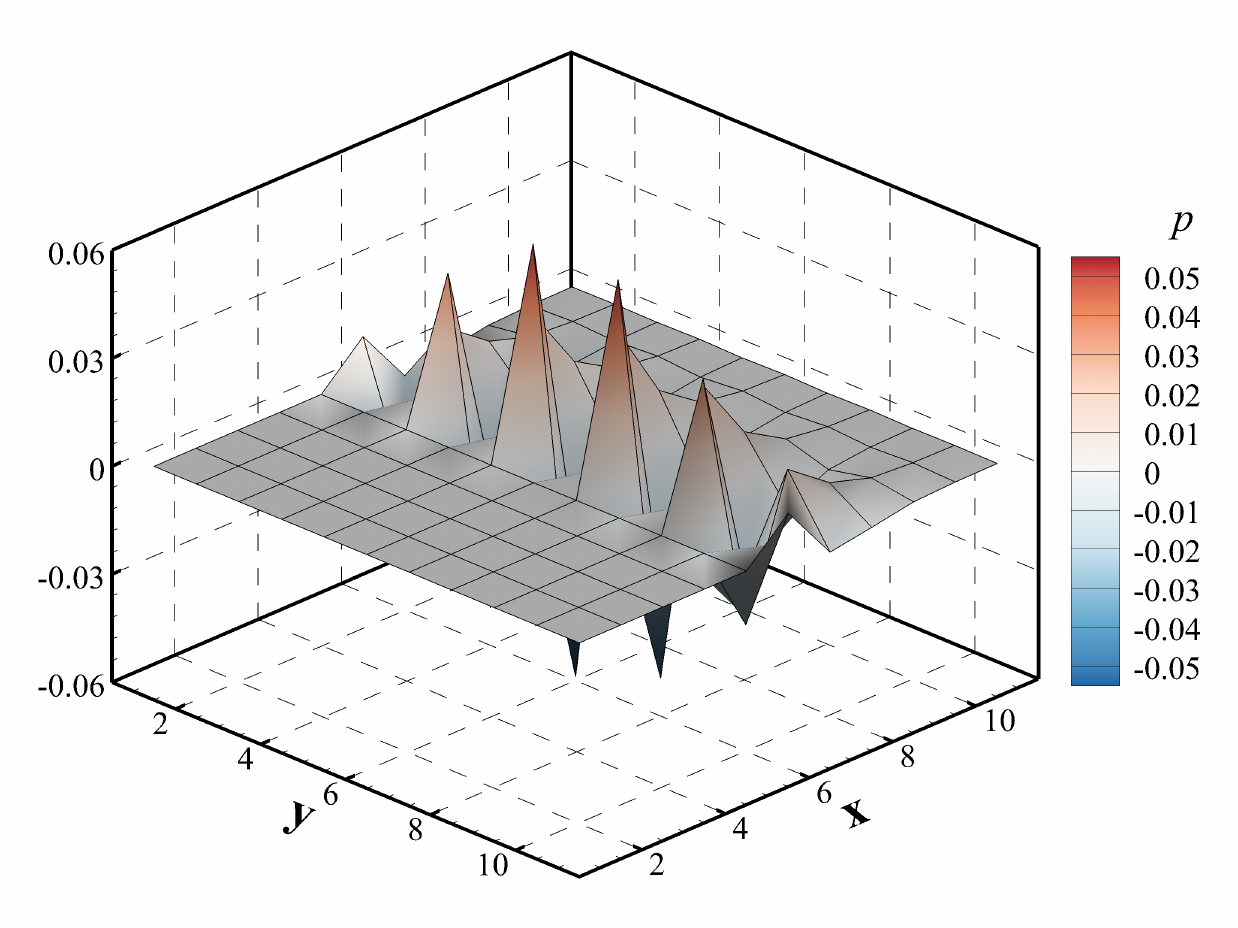}
	\end{minipage}
	}

	\centering
	\caption{Eigenvectors of the primitive variables of $ \lambda = 0.67356+0i $. (11$ \times $11 cells, fifth-order scheme with Roe solver, $ M_0=20 $ and $ \varepsilon=0.1 $.)}\label{fig eigenvector}
\end{figure}

\begin{figure}[htbp]
	\centering

	\subfigure[the profile of $ i=6 $(numerical shock structure)]{
	\begin{minipage}[t]{0.46\linewidth}
	\centering
	\includegraphics[width=0.9\textwidth]{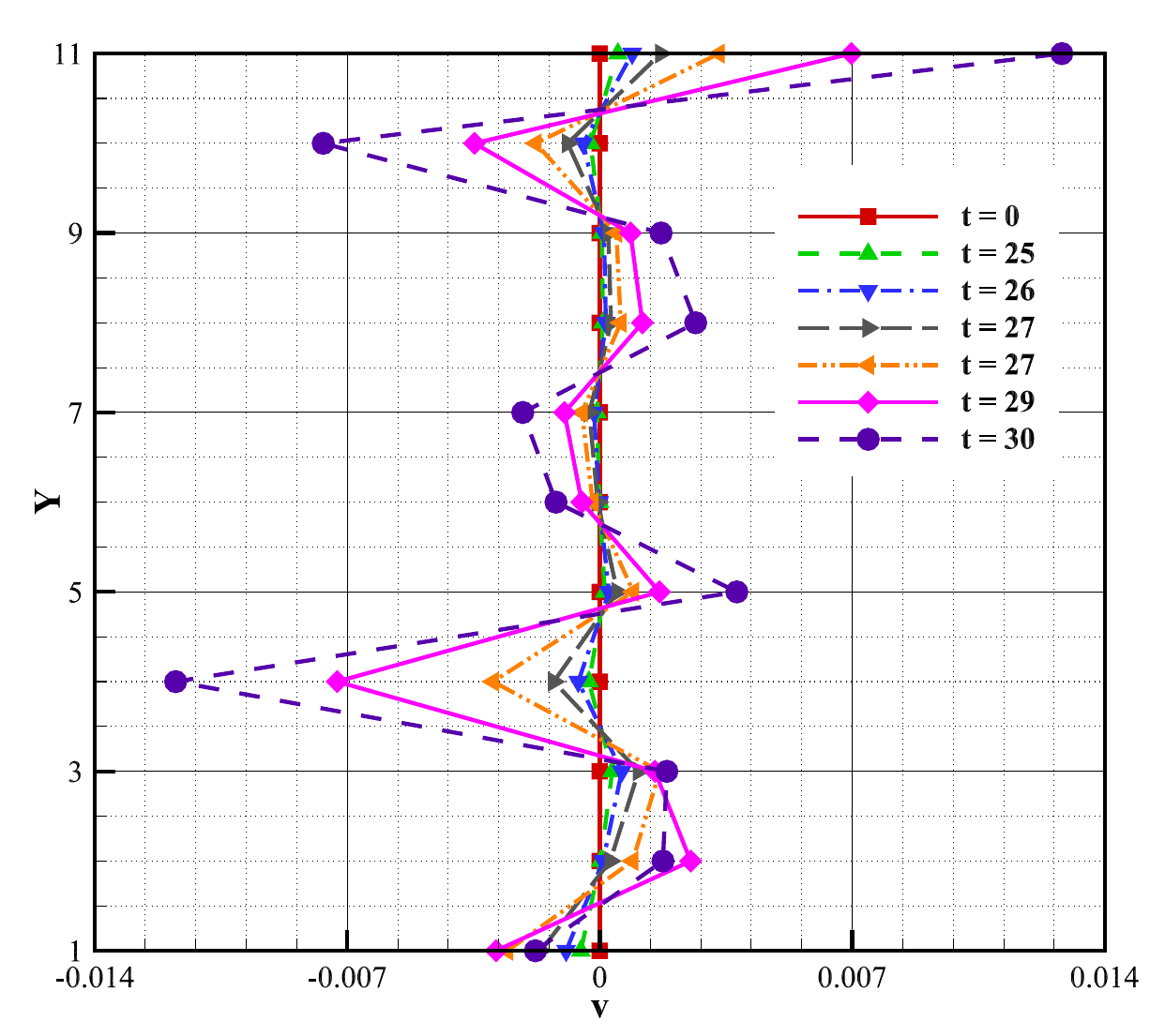}
	\end{minipage}
	}
	\subfigure[the profile of $ j=3 $]{
	\begin{minipage}[t]{0.46\linewidth}
	\centering
	\includegraphics[width=0.9\textwidth]{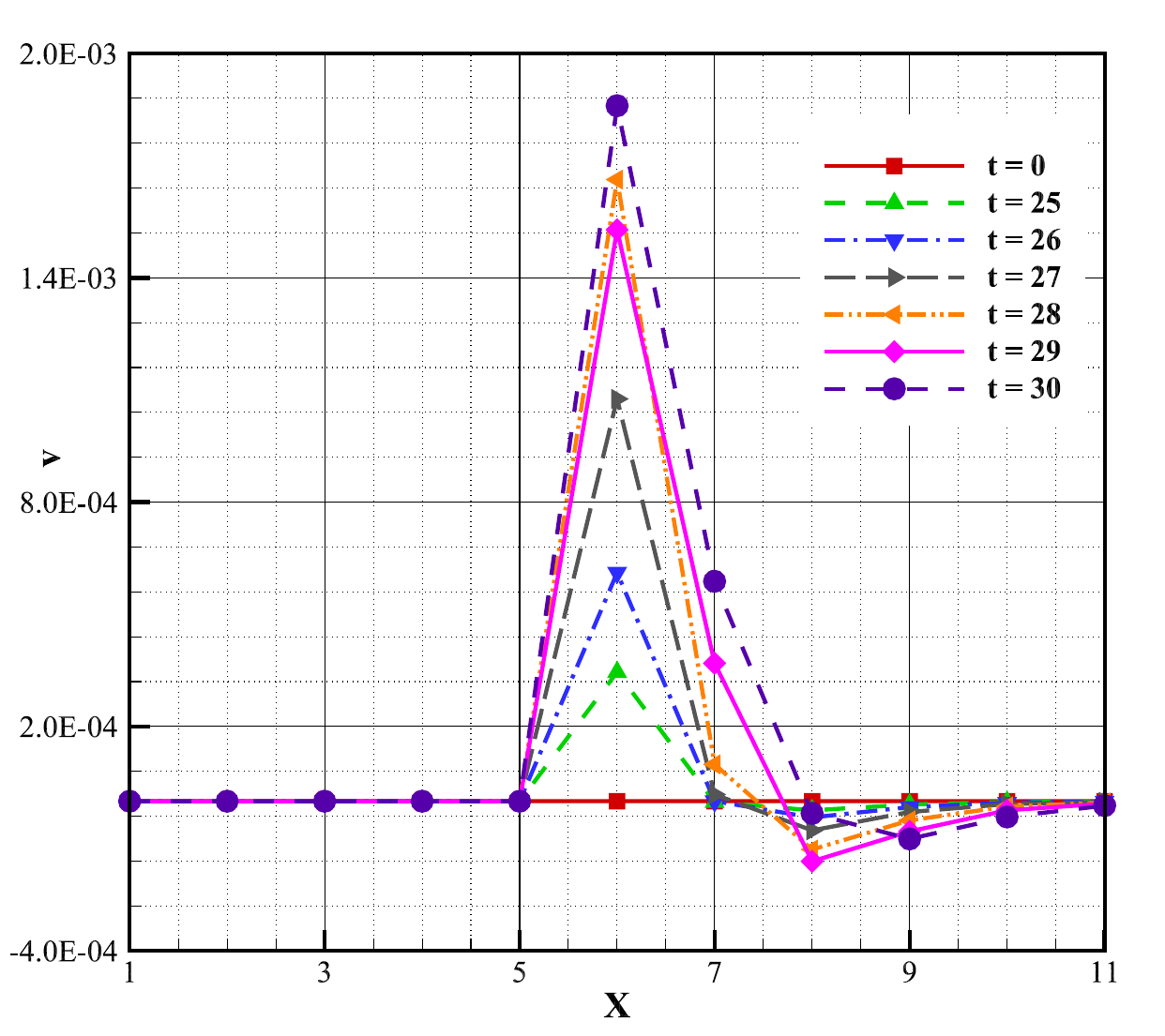}
	\end{minipage}
	}

	\subfigure[the profile of $ j=6 $]{
	\begin{minipage}[t]{0.46\linewidth}
	\centering
	\includegraphics[width=0.9\textwidth]{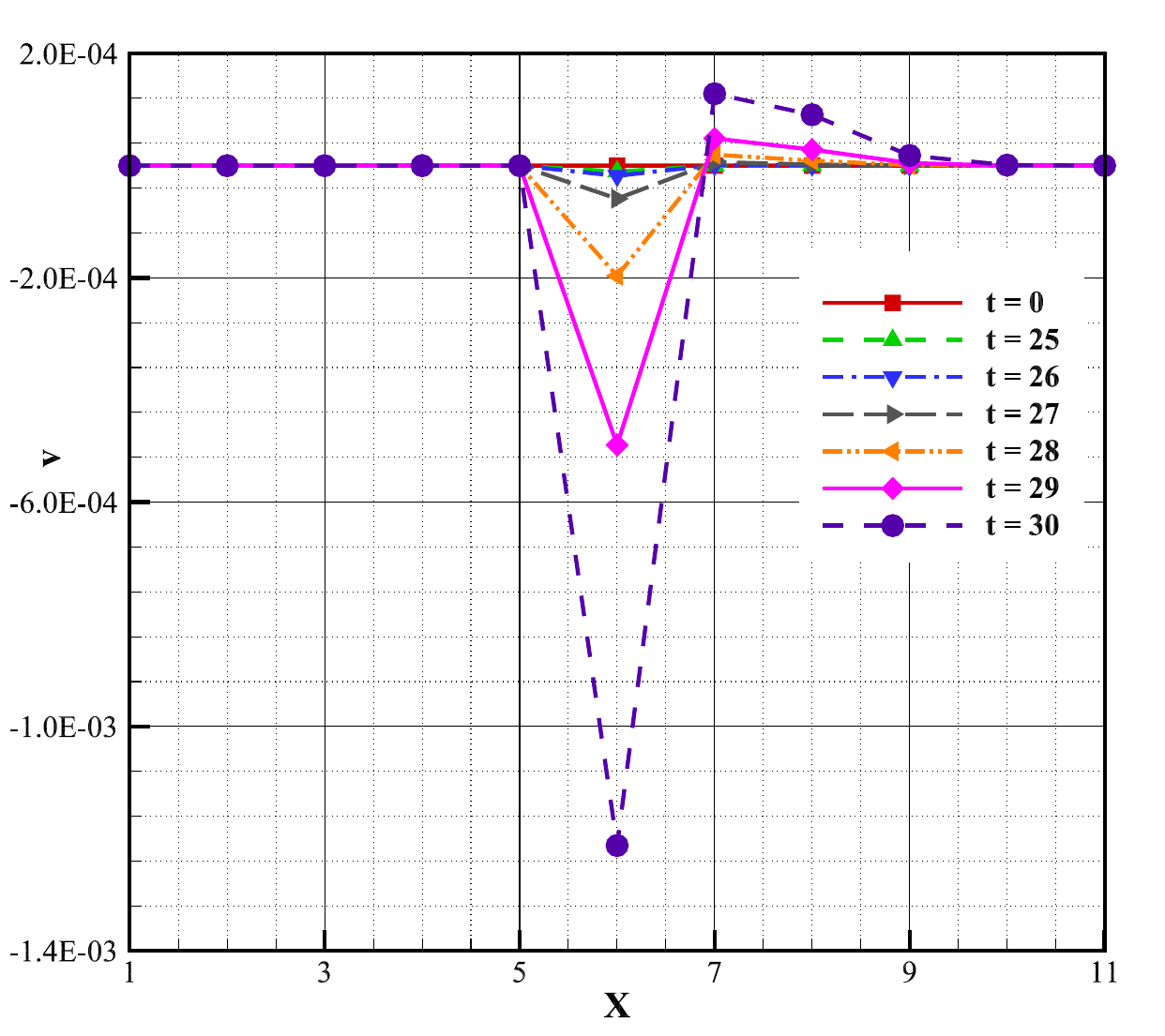}
	\end{minipage}
	}
	\subfigure[the profile of $ j=9 $]{
	\begin{minipage}[t]{0.46\linewidth}
	\centering
	\includegraphics[width=0.9\textwidth]{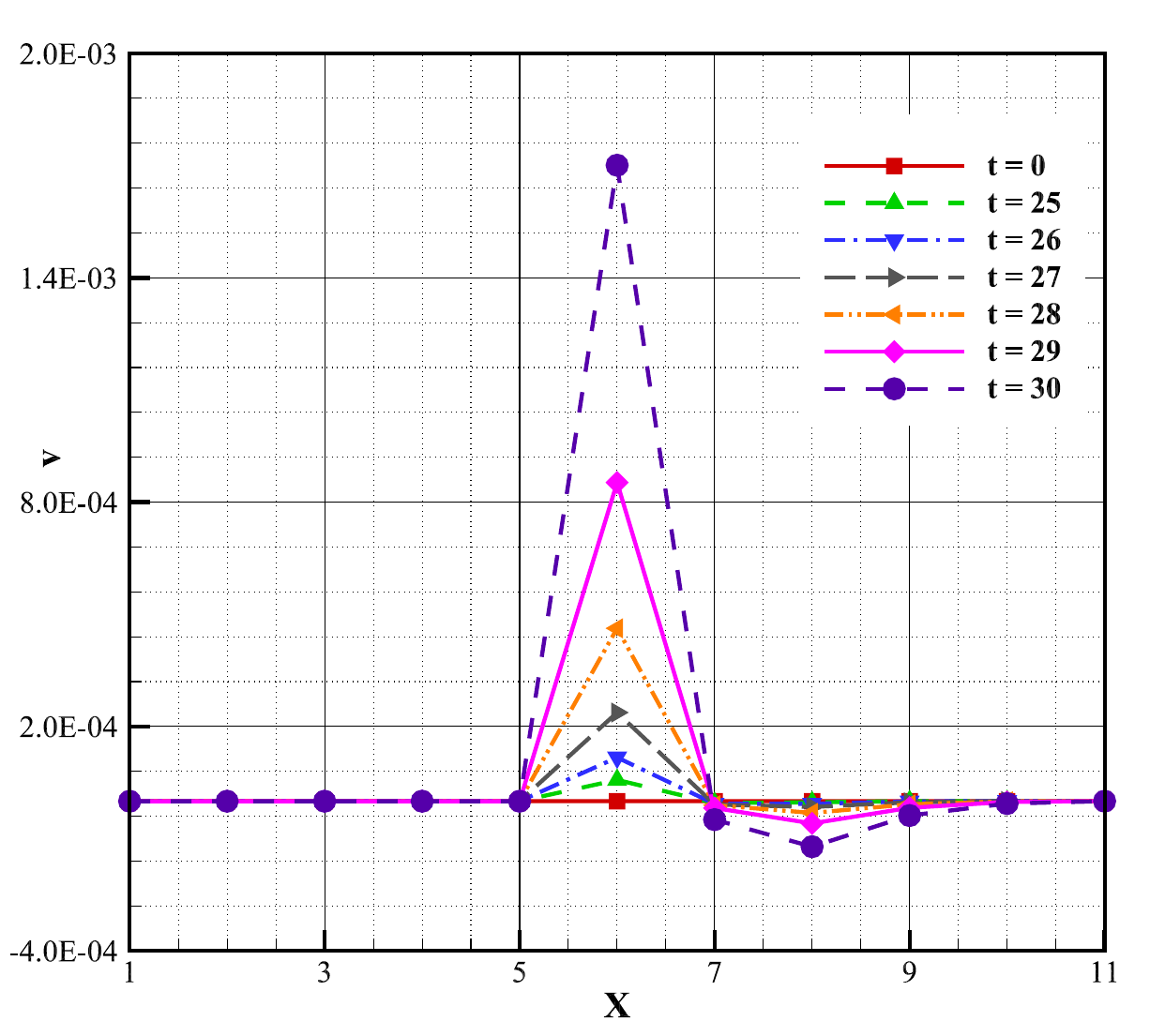}
	\end{minipage}
	}

	\centering
	\caption{The perturbation errors at different times and locations.}\label{fig v t}
\end{figure}

Fig.\ref{fig eigenvector} shows the eigenvectors of the primitive variables of the most unstable eigenvalue, which contains the information of the spatial-behavior of shock instability \cite{Dumbser_Matrix_2004}. As shown, the sawtooth distribution appears near $ i=6 $, which is the numerical shock structure. This indicates that the shock instability in the fifth-order case also originates from the numerical shock structure, just like the low-order case. To validate this and investigate the process of instability, we display the perturbation errors at different times in Fig.\ref{fig v t}. Fig.\ref{fig v t} (a) shows the perturbation error at the profile of $ i=6 $ (numerical shock structure), while Fig. ref{fig v t} (b)-(d) show the perturbation errors at $ j = 3,6,9 $ (perpendicular to the shock) respectively. As shown in Fig.\ref{fig v t} (a), the perturbation errors inside the numerical shock structure show a sawtooth profile with the time increasing. And it can also be found from Fig.\ref{fig v t} (b)-(d) that the perturbation error first appears inside the numerical shock structure, which is consistent with the results obtained by Fig.\ref{fig eigenvector}. Also, when the perturbation error inside the numerical shock structure develops to a certain value, it will propagate downstream, leading to instability in the downstream flow field. In contrast, the perturbation error in the upstream domain always remains zero.

\begin{figure}[htbp]
	\centering
	\includegraphics[width=0.7\textwidth]{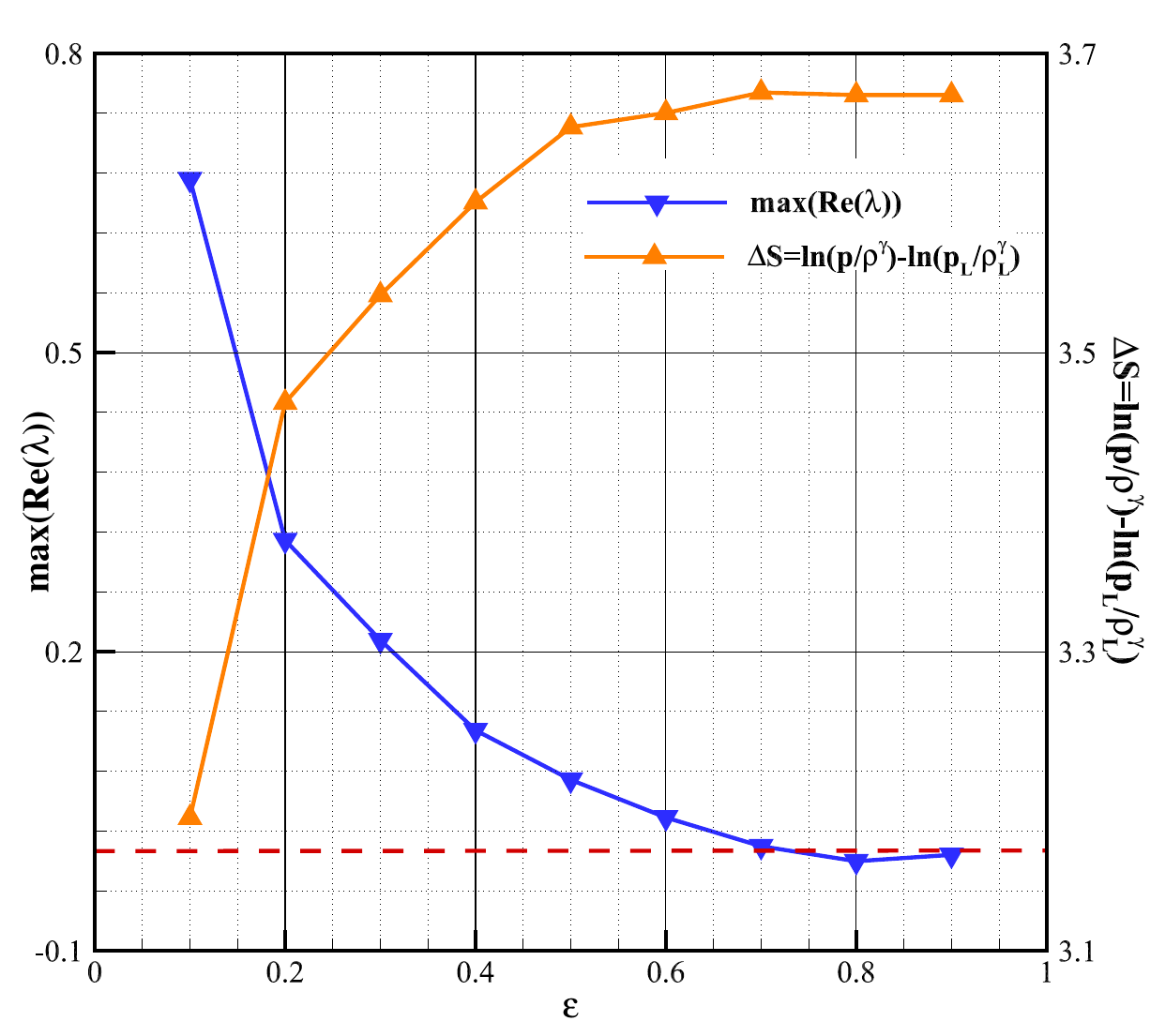}
	\caption{The entropy increase of the numerical shock structure. (Grid with 11$ \times $11 cells, fifthe order scheme with Roe solver, $ M_0 = 20 $.)}
	\label{fig entropy}
\end{figure}

According to the result of the above analysis, one can observe that for fifth-order schemes, the spatial location of shock instabilities also originates from the numerical shock structure, which is consistent with the first and second-order cases. Hence, it can be inferred that the state of the shock structure plays a vital role in ensuring shock-capturing stability, regardless of whether a low-order or high-order scheme is employed. This section utilizes the matrix stability analysis method proposed in this paper to explore the correlation between the state within the numerical shock structure and shock instabilities. The variations of $\max(\text{Re}(\lambda))$ with respect to $\varepsilon$ are depicted in Fig.\ref{fig validation} and Fig.\ref{fig entropy}. It is evident that, for fifth-order schemes, as $\varepsilon$ approaches 1, $\max(\text{Re}(\lambda))$ gradually decreases to a value below 0, signifying an improvement in the stability of shock-capturing. As indicated in (\ref{eq numerical shock structure}) and (\ref{eq epsilon}), as $\varepsilon$ approaches 1, the state within the numerical shock structure becomes more akin to the downstream region. Consequently, it can be concluded that a closer resemblance between the state in the numerical shock structure and the downstream region leads to greater computational stability, which is consistent with the observations made in the first and second-order cases. Therefore, this poses a question: why does the stability of shock-capturing improve as the state within the numerical shock structure approaches the downstream region? This question holds significant importance in understanding and curing the shock instability problem. It should be noted that, by theoretical analysis and numerical experiments, Xie et al.\cite{Xie_Further_2021} demonstrate that the stability of shock-capturing is closely related to the entropy state within the numerical shock structure. If enough entropy production is guaranteed, then the instability problem can be successfully eliminated. The orange curve in Fig.\ref{fig entropy} illustrates the relative increase in entropy from the upstream region to the numerical shock structure. It can be found from Fig.\ref{fig entropy} that as $\varepsilon$ increases, the entropy in the shock structure continuously amplifies, accompanied by a decrease in the corresponding $\max(\text{Re}(\lambda))$, indicating the improvement of the stability of shock-capturing. Such results can enhance the understanding of shock instabilities and offer significant aid towards curing the shock instability problem.

\subsection{Influence of the local characteristic decomposition}
The investigations in section \ref{subsection 4.1} - \ref{subsection 4.3} are all based on the reconstruction of primitive variables since it is widely employed in hypersonic flow simulations due to its lower cost \cite{Rider_Methods_1993, Zanotti_Efficient_2016}. In fact, it has been demonstrated that the local characteristic decomposition can significantly reduce oscillations when the order of accuracy is high, although it will increase the computational cost significantly \cite{Qiu_Construction_2002,Shu_Essentially_2020}. Previous studies have primarily focused on the effectiveness of local characteristic decomposition in reducing oscillations. However, there remains a question of whether local characteristic decomposition can also mitigate shock instabilities in high-order schemes. This problem will be addressed and examined in this section.

\begin{figure}[htbp]
	\centering

	\subfigure[numerical shock structure]{
	\begin{minipage}[t]{0.45\linewidth}
	\centering
	\includegraphics[width=0.9\textwidth]{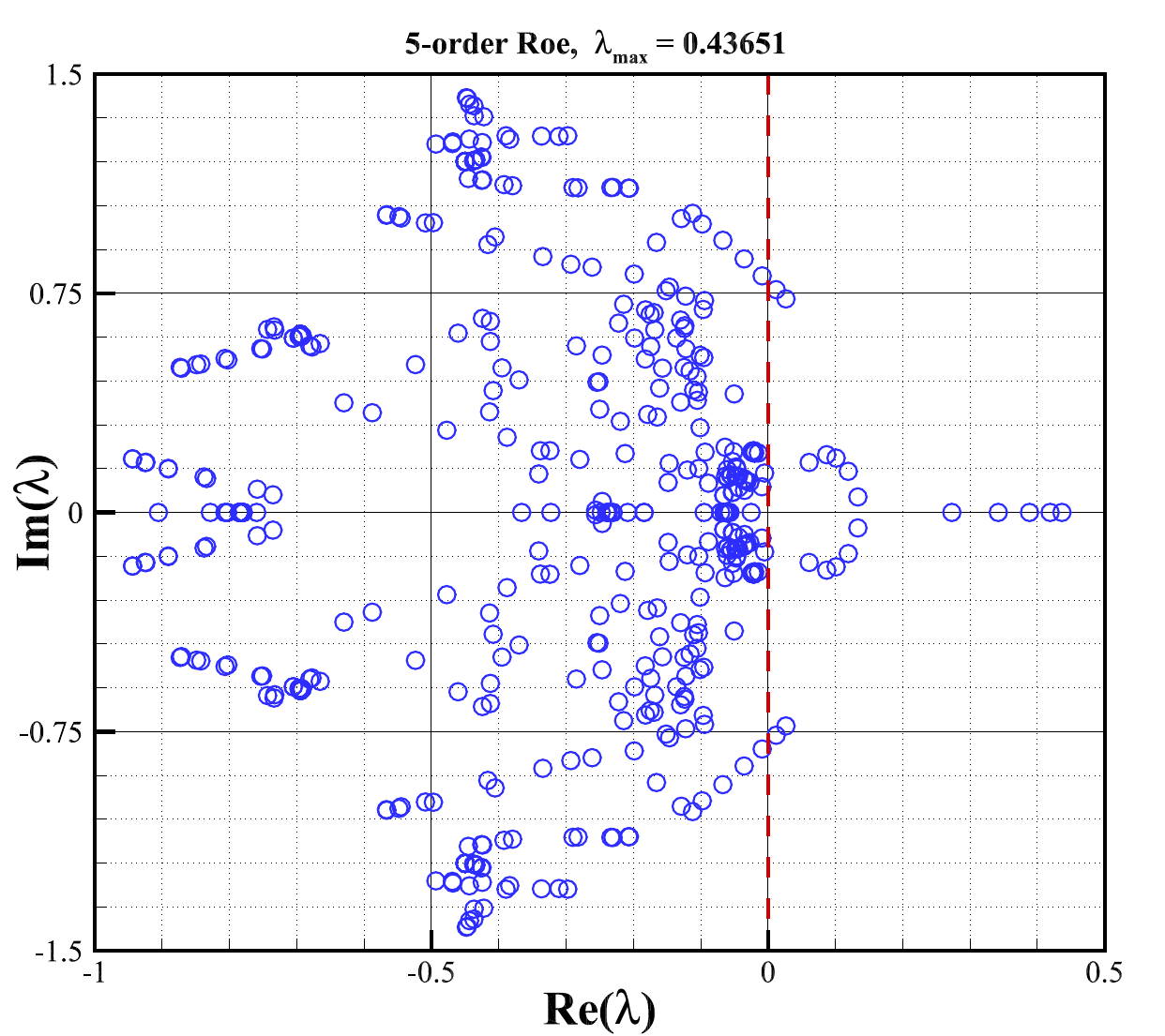}
	\end{minipage}
	}
	\subfigure[wave patterns of \textit{L,M}]{
	\begin{minipage}[t]{0.46\linewidth}
	\centering
	\includegraphics[width=0.9\textwidth]{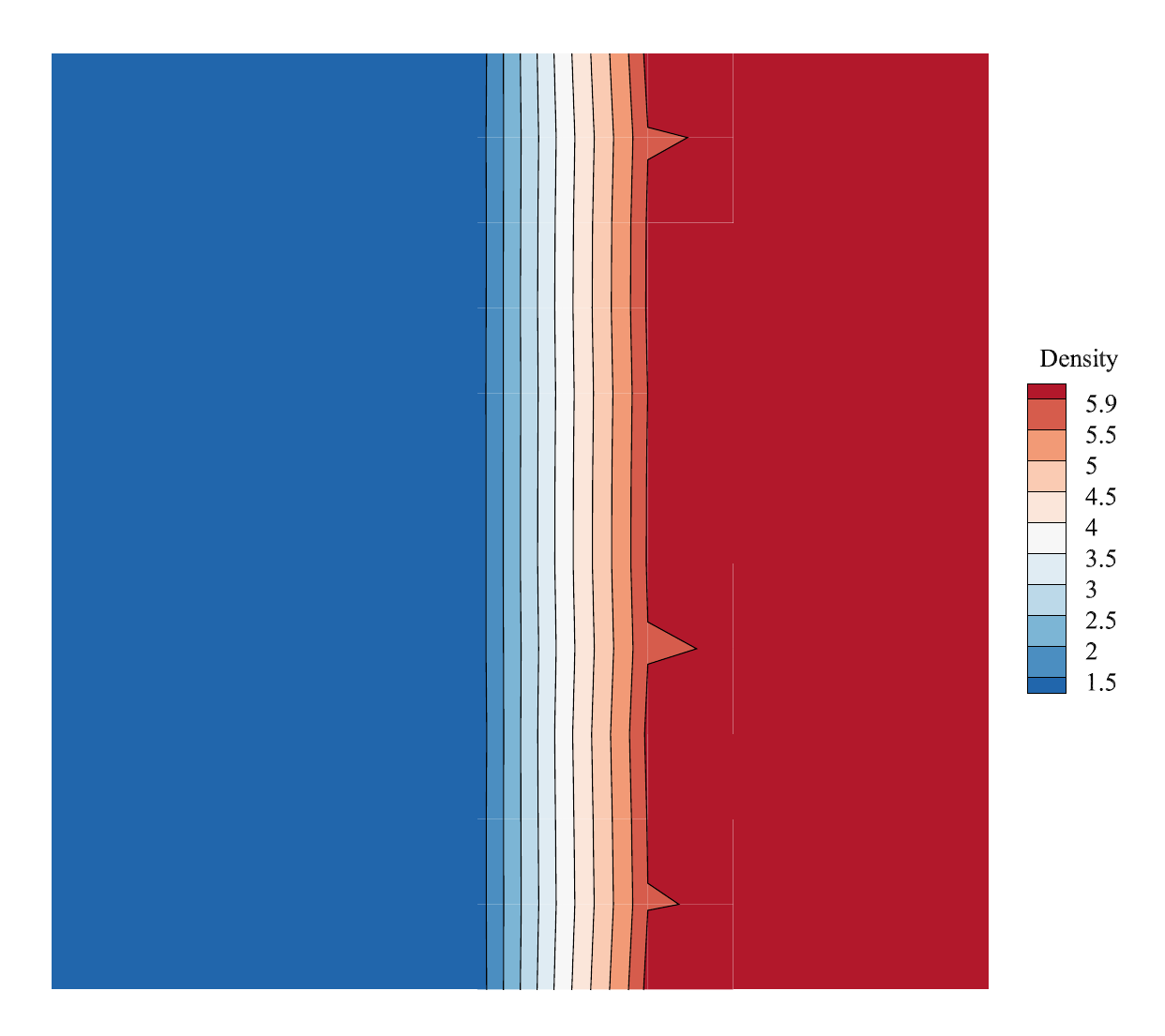}
	\end{minipage}
	}

	\centering
	\caption{Matrix stability analysis and the corresponding density with local characteristic decomposition.(Grid with 11$ \times $11 cells, fifthe order scheme with Roe solver, $ M_0 = 20 $, and $ \varepsilon =0.1 $.)}\label{fig characteristic}
\end{figure}

\begin{figure}[htbp]
	\centering
	\includegraphics[width=0.6\textwidth]{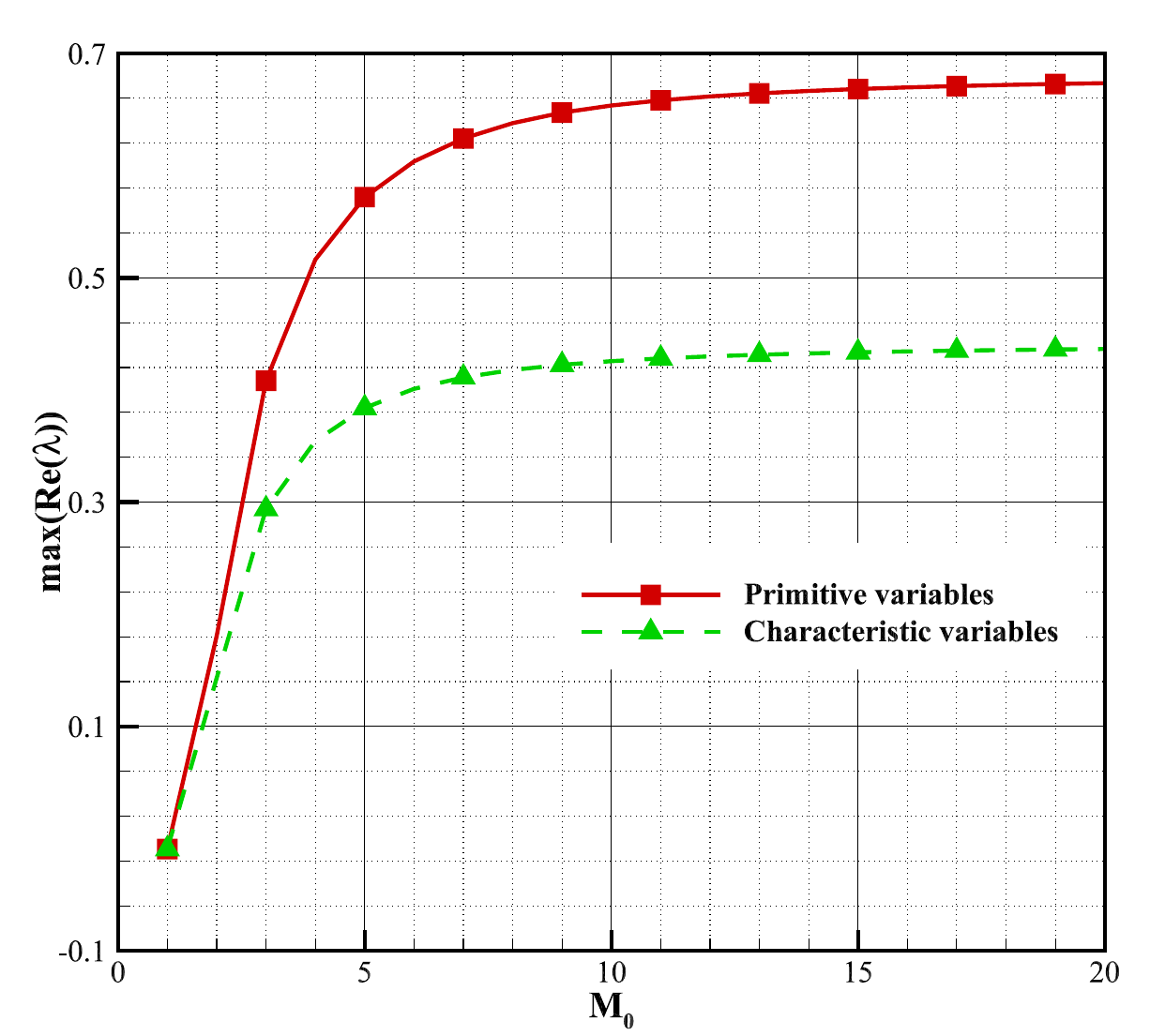}
	\caption{The comparison of $\max(\text{Re}(\lambda))$ between primitive and characteristic variables.}
	\label{fig variable Ma}
\end{figure}

To investigate the influence of the local characteristic decomposition on shock instabilities, the matrix stability analysis method is employed, the detail of which can be found in \ref{Appendix B}. Fig.\ref{fig characteristic} shows the result of the matrix analysis and the corresponding density contour. By comparing Fig.\ref{fig characteristic} with Fig.\ref{fig Roe matrix}, it is evident that there is a smaller maximum real part of $\lambda$ ($\max(\text{Re}(\lambda))$) when the local characteristic decomposition is employed. This indicates that the evolution of perturbation errors will occur at a slower rate with the use of local characteristic decomposition. So, the local characteristic decomposition is helpful to mitigate the shock instability. However, it is important to note that $\max(\text{Re}(\lambda))$ still exceeds zero, suggesting that the computation remains unstable. Such conclusion is validated by Fig.\ref{fig characteristic} (b), which displays the distorted shock profile. Furthermore, Fig.\ref{fig variable Ma} compares $\max(\text{Re}(\lambda))$ of primitive and characteristic variables with various Mach numbers. It can be observed that the maximal real part of the eigenvalues for the primitive variables always surpasses that of the characteristic variables. Despite utilizing the local characteristic decomposition, positive $\max(\text{Re}(\lambda))$ persists. Consequently, it can be concluded that the local characteristic decomposition can help to alleviate shock instabilities, although it does not fully eliminate them. As a result, in the pursuit of developing robust high-order schemes, employing the local characteristic decomposition can contribute to more stable capturing of strong shocks. 

\section{Conclusion}\label{section 5}
In this study, our primary focus is to explore the fundamental numerical characteristics and underlying mechanisms of shock instabilities in fifth-order finite-volume WENO schemes. To achieve this, we have established the matrix stability analysis method for the fifth-order WENO scheme. This method can predict the evolution of perturbation errors and links the ability of schemes to stably capture shocks with the eigenvalues of the stability matrix. Results reveal that when employing low-dissipative solvers like Roe and HLLC, fifth-order schemes can yield unstable results, similar to the first- and second-order cases. One of the most significant findings of this study is that dissipative solvers such as HLL and van Leer also suffer from shock instabilities. Further investigations demonstrate that this phenomenon occurs due to the excessively high spatial accuracy near the numerical shock structure. The computation remains stable if the spatial accuracy is first- or second-order near the shock structure. Additionally, it has been observed that the shock instability problem of fifth-order schemes is a multidimensional coupling problem. To stably capture strong shocks with fifth-order schemes, it is crucial to have sufficient dissipation on transverse faces and ensure at least two points within the numerical shock structure in the direction perpendicular to the shock. Through matrix analysis and numerical experiments, we have determined that the instability of fifth-order schemes originates from the numerical shock structure. Thus, the state of the shock structure plays a vital role in ensuring shock-capturing stability. When the state within the numerical shock structure becomes closer to the downstream region, the computation tends to be more stable. This paper also investigates the influence of local characteristic decomposition on shock instabilities in fifth-order schemes, which has been proven to be helpful in mitigating the instability. However, it is important to note that the instability is still exist.

This study paves the way for a better understanding of the shock instability problem of fifth-order schemes and provides several guidance for developing robust high-order schemes
\begin{itemize}
	\item [(1)]The spatial accuracy near the numerical shock structure should be limited to first- or second-order;

	\item [(2)]There should be sufficient dissipation on transverse faces and at least two points within the numerical shock structure in the direction perpendicular to the shock;
	
	\item [(3)]Targeted treatments can be implemented to the numerical shock structure to suppress the instability;
	
	\item [(4)]The employment of local characteristic decomposition can help to stably capture strong shocks.
\end{itemize}
Based on these findings, we will devote to developing robust high-order schemes in the follow-up work.

\section*{Acknowledgements}	
This work was supported by National Natural Science Foundation of China (Grant No.12202490), Natural Science Foundation of Hunan Province, China (Grant No. 11472004), the Scientific Research Foundation of NUDT (Grant No. ZK21-10).

\appendix
\section{The matrix stability analysis for primitive variables}\label{Appendix A}
\renewcommand\theequation{A.\arabic{equation}}

The fifth-order WENO reconstruction procedure can also be implemented on primitive variables. So, if we assume that 
\begin{equation}  \label{eq primitive variables decomposition}
	\mathbf{W}_{i,j}=\bar{\mathbf{W}}_{i,j} + \delta\mathbf{W}_{i,j},
\end{equation}
(\ref{eq fifth-order variables decomposition}) can also be written in the form of the primitive variables
\begin{equation}  \label{eq fifth-order variables decomposition for primitive variables}
	\begin{aligned}
		\delta \mathbf{W}_{i+1/2,j}^L &=\frac{1}{3} \boldsymbol{\omega}_{i+1/2,j}^{L,0} \delta\mathbf{W}_{i-2,j}-\frac{1}{6}\left(7\boldsymbol{\omega}_{i+1/2,j}^{L,0}+\boldsymbol{\omega}_{i+1/2,j}^{L,1}\right)\delta\mathbf{W}_{i-1,j}\\
		&+\frac{1}{6}\left(11\boldsymbol{\omega}_{i+1/2,j}^{L,0}+5\boldsymbol{\omega}_{i+1/2,j}^{L,1}+2\boldsymbol{\omega}_{i+1/2,j}^{L,2}\right)\delta\mathbf{W}_{i,j}\\
		&+\frac{1}{6}\left(2\boldsymbol{\omega}_{i+1/2,j}^{L,1}+5\boldsymbol{\omega}_{i+1/2,j}^{L,2} \right)\delta\mathbf{W}_{i+1,j}-\frac{1}{6}\boldsymbol{\omega}_{i+1/2,j}^{L,2}\delta\mathbf{W}_{i+2,j},\\
		\delta \mathbf{W}_{i+1/2}^R&=\frac{1}{3} \boldsymbol{\omega}_{i+1/2,j}^{R,0} \delta\mathbf{W}_{i+3,j}-\frac{1}{6}\left(7\boldsymbol{\omega}_{i+1/2,j}^{R,0}+\boldsymbol{\omega}_{i+1/2,j}^{R,1}\right)\delta\mathbf{W}_{i+2,j}\\
		&+\frac{1}{6}\left(11\boldsymbol{\omega}_{i+1/2,j}^{R,0}+5\boldsymbol{\omega}_{i+1/2,j}^{R,1}+2\boldsymbol{\omega}_{i+1/2,j}^{R,2}\right)\delta\mathbf{W}_{i+1,j}\\
		&+\frac{1}{6}\left(2\boldsymbol{\omega}_{i+1/2,j}^{R,1}+5\boldsymbol{\omega}_{i+1/2,j}^{R,2} \right)\delta\mathbf{W}_{i,j}-\frac{1}{6}\boldsymbol{\omega}_{i+1/2,j}^{R,2}\delta\mathbf{W}_{i-1,j},
	\end{aligned}
\end{equation}
where $ \mathbf{W}=(\rho,u,v,p)^T $ is the vector of primitive variables. Then the flux function $ \mathbf{F}_{i+1/2,j} $ can be linearized around the steady mean value as
\begin{equation}
	\begin{aligned}\label{eq simplify second-order flux linearized for primitive variables}
		\mathbf{F}_{i+1/2,j}& =\mathbf{F}_{i+1/2,j}\left(\bar{{\mathbf{W}}}_{i+1/2,j}^{L},\bar{{\mathbf{W}}}_{i+1/2,j}^{R}\right)  \\
		&+\boldsymbol{\alpha}_{i+1/2,j}^{L}\delta\mathbf{W}_{i-2,j}+\boldsymbol{\beta}_{i+1/2,j}^{L}\delta\mathbf{W}_{i-1,j}+\boldsymbol{\chi}_{i+1/2,j}^{L}\delta\mathbf{W}_{i,j} \\
		&+\boldsymbol{\chi}_{i+1/2,j}^{R}\delta\mathbf{W}_{i+1,j}+\boldsymbol{\beta}_{i+1/2,j}^{R}\delta\mathbf{W}_{i+2,j}+\boldsymbol{\alpha}_{i+1/2,j}^{R}\delta\mathbf{W}_{i+3,j}
	\end{aligned}.
\end{equation}
The expressions of $ \boldsymbol{\alpha} $, $ \boldsymbol{\beta} $, and $ \boldsymbol{\chi} $ can refer to ({eq alpha beta chi}), but they are based on the primitive variables here. Substituting (\ref{eq primitive variables decomposition}) and (\ref{eq simplify second-order flux linearized for primitive variables}) into (\ref{eq discrete Euler equations}), it can be obtained that
\begin{equation}\label{eq linear error evolution model for primitive variables}
	\begin{aligned}
	  \frac{\mathrm{d} \delta \mathbf{U}_{i, j}}{\mathrm{dt}}=&\left(\frac{\text{d}\mathbf{U}_{i,j}}{\text{d}\mathbf{W}_{i,j}}\right)^{-1}\left[\left(\boldsymbol{\zeta}_{i+1/2,j}^{L}+\boldsymbol{\zeta}_{i,j+1/2}^{L}+\boldsymbol{\zeta}_{i-1/2,j}^{R}+\boldsymbol{\zeta}_{i,j-1/2}^{R}\right)\delta\mathbf{U}_{i,j}\right.\\
	  +&\left(\boldsymbol{\zeta}_{i+1/2,j}^R+\boldsymbol{\eta}_{i-1/2,j}^R\right)\delta\mathbf{U}_{i+1,j}
	  +\left(\boldsymbol{\zeta}_{i,j+1/2}^R+\boldsymbol{\eta}_{i,j-1/2}^R\right)\delta\mathbf{U}_{i,j+1}\\
	  +&\left(\boldsymbol{\zeta}_{i-1/2,j}^L+\boldsymbol{\eta}_{i+1/2,j}^L\right)\delta\mathbf{U}_{i-1,j}
	  +\left(\boldsymbol{\zeta}_{i,j-1/2}^L+\boldsymbol{\eta}_{i,j+1/2}^L\right)\delta\mathbf{U}_{i,j-1}\\
	  +&\left(\boldsymbol{\eta}_{i+1/2,j}^R+\boldsymbol{\theta}_{i-1/2,j}^R\right)\delta\mathbf{U}_{i+2,j}
	  +\left(\boldsymbol{\eta}_{i,j+1/2}^R+\boldsymbol{\theta}_{i,j-1/2}^R\right)\delta\mathbf{U}_{i,j+2}\\
	  +&\left(\boldsymbol{\eta}_{i-1/2,j}^L+\boldsymbol{\theta}_{i+1/2,j}^L\right)\delta\mathbf{U}_{i-2,j}
	  +\left(\boldsymbol{\eta}_{i,j-1/2}^L+\boldsymbol{\theta}_{i,j+1/2}^L\right)\delta\mathbf{U}_{i,j-2}\\
	  +&\left.\boldsymbol{\theta}_{i+1/2,j}^R\delta\mathbf{U}_{i+3,j}+\boldsymbol{\theta}_{i,j+1/2}^R\mathbf{U}_{i,j+3}+\boldsymbol{\theta}_{i-1/2,j}^L\mathbf{U}_{i-3,j}+\boldsymbol{\theta}_{i,j-1/2}^L\delta\mathbf{U}_{i,j-3}\right]
	\end{aligned},
\end{equation}
where $ \boldsymbol{\zeta} $, $ \boldsymbol{\eta} $, and $ \boldsymbol{\theta} $ are computed by (\ref{eq zeta eta theta}). $ \frac{\text{d}\mathbf{U}_{i,j}}{\text{d}\mathbf{W}_{i,j}} $ is the transformation matrix between the conservative variables and the primitive variables
\begin{equation}  
	\frac{\mathrm{d} \mathbf{U}_{i, j}}{\mathrm{~d} \mathbf{W}_{i, j}}=\left[\begin{array}{cccc}
		1 & 0 & 0 & 0 \\
		u & \rho & 0 & 0 \\
		v & 0 & \rho & 0 \\
		\left(u^2+v^2\right) / 2 & \rho u & \rho v & 1 /(\gamma-1)
		\end{array}\right]_{i, j} .
\end{equation}
(\ref{eq linear error evolution model for primitive variables}) is the primitive variables form of the linear error evolution model in $ \Omega_{i,j} $. (\ref{eq linear error evolution model for all cells}) and (\ref{eq solution of linear error evolution model}) can also be written in the primitive variables forms and the same stability criterion is shown in (\ref{eq stability criterion}).

\section{The matrix stability analysis for characteristic variables}\label{Appendix B}
\renewcommand\theequation{B.\arabic{equation}}
The similar expressions like (\ref{eq variables decomposition}) and (\ref{eq primitive variables decomposition}) can also be obtained for characteristic variables
\begin{equation}  \label{eq fifth-order variables decomposition for characteristic variables}
	\begin{aligned}
		\delta \mathbf{V}_{i+1/2,j}^L &=\frac{1}{3} \boldsymbol{\omega}_{i+1/2,j}^{L,0} \delta\mathbf{V}_{i-2,j}-\frac{1}{6}\left(7\boldsymbol{\omega}_{i+1/2,j}^{L,0}+\boldsymbol{\omega}_{i+1/2,j}^{L,1}\right)\delta\mathbf{V}_{i-1,j}\\
		&+\frac{1}{6}\left(11\boldsymbol{\omega}_{i+1/2,j}^{L,0}+5\boldsymbol{\omega}_{i+1/2,j}^{L,1}+2\boldsymbol{\omega}_{i+1/2,j}^{L,2}\right)\delta\mathbf{V}_{i,j}\\
		&+\frac{1}{6}\left(2\boldsymbol{\omega}_{i+1/2,j}^{L,1}+5\boldsymbol{\omega}_{i+1/2,j}^{L,2} \right)\delta\mathbf{V}_{i+1,j}-\frac{1}{6}\boldsymbol{\omega}_{i+1/2,j}^{L,2}\delta\mathbf{V}_{i+2,j},\\
		\delta \mathbf{V}_{i+1/2}^R&=\frac{1}{3} \boldsymbol{\omega}_{i+1/2,j}^{R,0} \delta\mathbf{V}_{i+3,j}-\frac{1}{6}\left(7\boldsymbol{\omega}_{i+1/2,j}^{R,0}+\boldsymbol{\omega}_{i+1/2,j}^{R,1}\right)\delta\mathbf{V}_{i+2,j}\\
		&+\frac{1}{6}\left(11\boldsymbol{\omega}_{i+1/2,j}^{R,0}+5\boldsymbol{\omega}_{i+1/2,j}^{R,1}+2\boldsymbol{\omega}_{i+1/2,j}^{R,2}\right)\delta\mathbf{V}_{i+1,j}\\
		&+\frac{1}{6}\left(2\boldsymbol{\omega}_{i+1/2,j}^{R,1}+5\boldsymbol{\omega}_{i+1/2,j}^{R,2} \right)\delta\mathbf{V}_{i,j}-\frac{1}{6}\boldsymbol{\omega}_{i+1/2,j}^{R,2}\delta\mathbf{V}_{i-1,j}.
	\end{aligned}
\end{equation}
Here, $ \mathbf{V} $ denotes the vector of characteristic variables and the nonlinear weights are computed by the characteristic variables. The relation between conservative variables and characteristic variables is as follows
\begin{equation}  \label{eq conservative variables to characteristic variables}
	\mathbf{V}=\mathbf{L}\mathbf{U},
\end{equation}
where 
\begin{equation}  
	\mathbf{L}=\left[\begin{array}{cccc}\label{eq left eigenvector matrix}
		\frac{1}{2}\left(\frac{\gamma-1}{2 c^2} \mathbf{v}^2+\frac{q}{c}\right) & -\frac{1}{2}\left(\frac{\gamma-1}{c^2} u+\frac{n_x}{c}\right) & -\frac{1}{2}\left(\frac{\gamma-1}{c^2} v+\frac{n_y}{c}\right) & \frac{\gamma-1}{2 c^2} \\
		1-\frac{\gamma-1}{2 c^2} \mathbf{v}^2 & \frac{\gamma-1}{c^2} u & \frac{\gamma-1}{c^2} v & -\frac{\gamma-1}{c^2} \\
		\frac{1}{2}\left(\frac{\gamma-1}{2 c^2} \mathbf{v}^2-\frac{q}{c}\right) & -\frac{1}{2}\left(\frac{\gamma-1}{c^2} u-\frac{n_x}{c}\right) & -\frac{1}{2}\left(\frac{\gamma-1}{c^2} v-\frac{n_y}{c}\right) & \frac{\gamma-1}{2 c^2} \\
		-q_{\ell} & \ell_x & \ell_y & 0
		\end{array}\right]
\end{equation}
is the left eigenvector matrix. $ [\ell_x,\ell_y]^T = [-n_y,n_x]^T $ is the tangent vector (perpendicular to $ \mathbf{n} $), and $ q_{\ell} = u\ell_x+v\ell_y $ is the velocity component in this direction. Substituting (\ref{eq conservative variables to characteristic variables}) into (\ref{eq fifth-order variables decomposition for characteristic variables}), it can be obtained that
\begin{equation}  \label{eq fifth-order variables decomposition for primitive variables 2}
	\begin{aligned}
		\delta \mathbf{V}_{i+1/2,j}^L &=\frac{1}{3} \boldsymbol{\omega}_{i+1/2,j}^{L,0} \mathbf{L}_{i+1/2,j}\delta\mathbf{U}_{i-2,j}-\frac{1}{6}\left(7\boldsymbol{\omega}_{i+1/2,j}^{L,0}+\boldsymbol{\omega}_{i+1/2,j}^{L,1}\right)\mathbf{L}_{i+1/2,j}\delta\mathbf{U}_{i-1,j}\\
		&+\frac{1}{6}\left(11\boldsymbol{\omega}_{i+1/2,j}^{L,0}+5\boldsymbol{\omega}_{i+1/2,j}^{L,1}+2\boldsymbol{\omega}_{i+1/2,j}^{L,2}\right)\mathbf{L}_{i+1/2,j}\delta\mathbf{U}_{i,j}\\
		&+\frac{1}{6}\left(2\boldsymbol{\omega}_{i+1/2,j}^{L,1}+5\boldsymbol{\omega}_{i+1/2,j}^{L,2} \right)\mathbf{L}_{i+1/2,j}\delta\mathbf{U}_{i+1,j}-\frac{1}{6}\boldsymbol{\omega}_{i+1/2,j}^{L,2}\mathbf{L}_{i+1/2,j}\delta\mathbf{U}_{i+2,j},\\
		\delta \mathbf{V}_{i+1/2}^R&=\frac{1}{3} \boldsymbol{\omega}_{i+1/2,j}^{R,0} \mathbf{L}_{i+1/2,j}\delta\mathbf{U}_{i+3,j}-\frac{1}{6}\left(7\boldsymbol{\omega}_{i+1/2,j}^{R,0}+\boldsymbol{\omega}_{i+1/2,j}^{R,1}\right)\mathbf{L}_{i+1/2,j}\delta\mathbf{U}_{i+2,j}\\
		&+\frac{1}{6}\left(11\boldsymbol{\omega}_{i+1/2,j}^{R,0}+5\boldsymbol{\omega}_{i+1/2,j}^{R,1}+2\boldsymbol{\omega}_{i+1/2,j}^{R,2}\right)\mathbf{L}_{i+1/2,j}\delta\mathbf{U}_{i+1,j}\\
		&+\frac{1}{6}\left(2\boldsymbol{\omega}_{i+1/2,j}^{R,1}+5\boldsymbol{\omega}_{i+1/2,j}^{R,2} \right)\mathbf{L}_{i+1/2,j}\delta\mathbf{U}_{i,j}-\frac{1}{6}\boldsymbol{\omega}_{i+1/2,j}^{R,2}\mathbf{L}_{i+1/2,j}\delta\mathbf{U}_{i-1,j}.
	\end{aligned}
\end{equation}
Then, the flux fucntion $ \mathbf{F}_{i+1/2,j} $ can be linearized as
\begin{equation}
	\begin{aligned}\label{eq simplify second-order flux linearized for characteristic variables}
	  \mathbf{F}_{i+1/2,j}& =\mathbf{F}_{i+1/2,j}\left(\bar{{\mathbf{V}}}_{i+1/2,j}^{L},\bar{{\mathbf{V}}}_{i+1/2,j}^{R}\right)  \\
	  &+\boldsymbol{\alpha}_{i+1/2,j}^{L}\mathbf{L}_{i+1/2,j}\delta\mathbf{U}_{i-2,j}+\boldsymbol{\beta}_{i+1/2,j}^{L}\mathbf{L}_{i+1/2,j}\delta\mathbf{U}_{i-1,j}+\boldsymbol{\chi}_{i+1/2,j}^{L}\mathbf{L}_{i+1/2,j}\delta\mathbf{U}_{i,j} \\
	  &+\boldsymbol{\chi}_{i+1/2,j}^{R}\mathbf{L}_{i+1/2,j}\delta\mathbf{U}_{i+1,j}+\boldsymbol{\beta}_{i+1/2,j}^{R}\mathbf{L}_{i+1/2,j}\delta\mathbf{U}_{i+2,j}+\boldsymbol{\alpha}_{i+1/2,j}^{R}\mathbf{L}_{i+1/2,j}\delta\mathbf{U}_{i+3,j}
	\end{aligned},
 \end{equation}
Substituting (\ref{eq fifth-order variables decomposition for primitive variables 2}) into (\ref{eq discrete Euler equations}), it can be obtained that
\begin{equation}\label{eq linear error evolution model for characteristic variables}
	\begin{aligned}
	  \frac{\mathrm{d} \delta \mathbf{U}_{i, j}}{\mathrm{dt}}=&\left(\boldsymbol{\zeta}_{i+1/2,j}^{L}+\boldsymbol{\zeta}_{i,j+1/2}^{L}+\boldsymbol{\zeta}_{i-1/2,j}^{R}+\boldsymbol{\zeta}_{i,j-1/2}^{R}\right)\delta\mathbf{U}_{i,j}\\
	  +&\left(\boldsymbol{\zeta}_{i+1/2,j}^R+\boldsymbol{\eta}_{i-1/2,j}^R\right)\delta\mathbf{U}_{i+1,j}
	  +\left(\boldsymbol{\zeta}_{i,j+1/2}^R+\boldsymbol{\eta}_{i,j-1/2}^R\right)\delta\mathbf{U}_{i,j+1}\\
	  +&\left(\boldsymbol{\zeta}_{i-1/2,j}^L+\boldsymbol{\eta}_{i+1/2,j}^L\right)\delta\mathbf{U}_{i-1,j}
	  +\left(\boldsymbol{\zeta}_{i,j-1/2}^L+\boldsymbol{\eta}_{i,j+1/2}^L\right)\delta\mathbf{U}_{i,j-1}\\
	  +&\left(\boldsymbol{\eta}_{i+1/2,j}^R+\boldsymbol{\theta}_{i-1/2,j}^R\right)\delta\mathbf{U}_{i+2,j}
	  +\left(\boldsymbol{\eta}_{i,j+1/2}^R+\boldsymbol{\theta}_{i,j-1/2}^R\right)\delta\mathbf{U}_{i,j+2}\\
	  +&\left(\boldsymbol{\eta}_{i-1/2,j}^L+\boldsymbol{\theta}_{i+1/2,j}^L\right)\delta\mathbf{U}_{i-2,j}
	  +\left(\boldsymbol{\eta}_{i,j-1/2}^L+\boldsymbol{\theta}_{i,j+1/2}^L\right)\delta\mathbf{U}_{i,j-2}\\
	  +&\boldsymbol{\theta}_{i+1/2,j}^R\delta\mathbf{U}_{i+3,j}+\boldsymbol{\theta}_{i,j+1/2}^R\mathbf{U}_{i,j+3}+\boldsymbol{\theta}_{i-1/2,j}^L\mathbf{U}_{i-3,j}+\boldsymbol{\theta}_{i,j-1/2}^L\delta\mathbf{U}_{i,j-3}
	\end{aligned},
\end{equation}
where
\begin{equation}
  \begin{aligned}\label{eq zeta eta theta for characteristic variables}
    \boldsymbol{\zeta}_{i\pm1/2,j\pm1/2}^{L/R} &= -\frac{\mathcal{L}_{i\pm1/2,j\pm1/2}}{\Omega_{i,j}}\boldsymbol{\chi}_{i\pm1/2,j\pm1/2}^{L/R}\mathbf{L}_{i\pm1/2,j\pm1/2},\\
    \boldsymbol{\eta}_{i\pm1/2,j\pm1/2}^{L/R} &= -\frac{\mathcal{L}_{i\pm1/2,j\pm1/2}}{\Omega_{i,j}}\boldsymbol{\beta}_{i\pm1/2,j\pm1/2}^{L/R}\mathbf{L}_{i\pm1/2,j\pm1/2},\\
    \boldsymbol{\theta}_{i\pm1/2,j\pm1/2}^{L/R} &= -\frac{\mathcal{L}_{i\pm1/2,j\pm1/2}}{\Omega_{i,j}}\boldsymbol{\alpha}_{i\pm1/2,j\pm1/2}^{L/R}\mathbf{L}_{i\pm1/2,j\pm1/2}.\\
  \end{aligned}
\end{equation}
(\ref{eq linear error evolution model for characteristic variables}) is the perturbation error evolution in $ \Omega_{i,j} $ when reconstructing characteristic variables. (\ref{eq linear error evolution model for all cells}) and (\ref{eq solution of linear error evolution model}) can also be written in the primitive variables forms and the same stability criterion is shown in (\ref{eq stability criterion}).

\bibliography{mybibfile}

\end{document}